
\ifx\shlhetal\undefinedcontrolsequence\let\shlhetal\relax\fi
\def\fmtname{AmS-TeX}

\def\fmtversion{2.1}
\catcode`\@=11
\ifx\amstexloaded@\relax\catcode`\@=\active
  \endinput\else\let\amstexloaded@\relax\fi
\newlinechar=`\^^J
\def\W@{\immediate\write\sixt@@n}
\def\CR@{\W@{^^J\fmtname - Version \fmtversion^^J%
COPYRIGHT 1985, 1990, 1991 - AMERICAN MATHEMATICAL SOCIETY^^J%
Use of this macro package is not restricted provided^^J%
each use is acknowledged upon publication.^^J}}
\CR@ \everyjob{\CR@}
\message{Loading definitions for}
\message{misc utility macros,}
\toksdef\toks@@=2
\long\def\rightappend@#1\to#2{\toks@{\\{#1}}\toks@@
 =\expandafter{#2}\xdef#2{\the\toks@@\the\toks@}\toks@{}\toks@@{}}
\def\alloclist@{}
\newif\ifalloc@
\def\showallocations{{\def\\{\immediate\write\m@ne}\alloclist@}\alloc@true}
\def\alloc@#1#2#3#4#5{\global\advance\count1#1by\@ne
 \ch@ck#1#4#2\allocationnumber=\count1#1
 \global#3#5=\allocationnumber
 \edef\next@{\string#5=\string#2\the\allocationnumber}%
 \expandafter\rightappend@\next@\to\alloclist@}
\newcount\count@@
\newcount\count@@@
\def\FN@{\futurelet\next}
\def\DN@{\def\next@}
\def\DNii@{\def\nextii@}
\def\RIfM@{\relax\ifmmode}
\def\RIfMIfI@{\relax\ifmmode\ifinner}
\def\setboxz@h{\setbox\z@\hbox}
\def\wdz@{\wd\z@}
\def\boxz@{\box\z@}
\def\setbox@ne{\setbox\@ne}
\def\wd@ne{\wd\@ne}
\def\iterate{\body\expandafter\iterate\else\fi}
\def\err@#1{\errmessage{AmS-TeX error: #1}}
\newhelp\defaulthelp@{Sorry, I already gave what help I could...^^J
Maybe you should try asking a human?^^J
An error might have occurred before I noticed any problems.^^J
``If all else fails, read the instructions.''}
\def\Err@{\errhelp\defaulthelp@\err@}
\def\eat@#1{}
\def\in@#1#2{\def\in@@##1#1##2##3\in@@{\ifx\in@##2\in@false\else\in@true\fi}%
 \in@@#2#1\in@\in@@}
\newif\ifin@
\def\space@.{\futurelet\space@\relax}
\space@. %
\newhelp\athelp@
{Only certain combinations beginning with @ make sense to me.^^J
Perhaps you wanted \string\@\space for a printed @?^^J
I've ignored the character or group after @.}
{\catcode`\~=\active 
 \lccode`\~=`\@ \lowercase{\gdef~{\FN@\at@}}}
\def\at@{\let\next@\at@@
 \ifcat\noexpand\next a\else\ifcat\noexpand\next0\else
 \ifcat\noexpand\next\relax\else
   \let\next\at@@@\fi\fi\fi
 \next@}
\def\at@@#1{\expandafter
 \ifx\csname\space @\string#1\endcsname\relax
  \expandafter\at@@@ \else
  \csname\space @\string#1\expandafter\endcsname\fi}
\def\at@@@#1{\errhelp\athelp@ \err@{\Invalid@@ @}}
\def\atdef@#1{\expandafter\def\csname\space @\string#1\endcsname}
\newhelp\defahelp@{If you typed \string\define\space cs instead of
\string\define\string\cs\space^^J
I've substituted an inaccessible control sequence so that your^^J
definition will be completed without mixing me up too badly.^^J
If you typed \string\define{\string\cs} the inaccessible control sequence^^J
was defined to be \string\cs, and the rest of your^^J
definition appears as input.}
\newhelp\defbhelp@{I've ignored your definition, because it might^^J
conflict with other uses that are important to me.}
\def\define{\FN@\define@}
\def\define@{\ifcat\noexpand\next\relax
 \expandafter\define@@\else\errhelp\defahelp@                               
 \err@{\string\define\space must be followed by a control
 sequence}\expandafter\def\expandafter\nextii@\fi}                          
\def\undefined@@@@@@@@@@{}
\def\preloaded@@@@@@@@@@{}
\def\next@@@@@@@@@@{}
\def\define@@#1{\ifx#1\relax\errhelp\defbhelp@                              
 \err@{\string#1\space is already defined}\DN@{\DNii@}\else
 \expandafter\ifx\csname\expandafter\eat@\string                            
 #1@@@@@@@@@@\endcsname\undefined@@@@@@@@@@\errhelp\defbhelp@
 \err@{\string#1\space can't be defined}\DN@{\DNii@}\else
 \expandafter\ifx\csname\expandafter\eat@\string#1\endcsname\relax          
 \global\let#1\undefined\DN@{\def#1}\else\errhelp\defbhelp@
 \err@{\string#1\space is already defined}\DN@{\DNii@}\fi
 \fi\fi\next@}

\def\predefine#1#2{\let#1#2}
\def\undefine#1{\let#1\undefined}
\message{page layout,}
\newdimen\captionwidth@
\captionwidth@\hsize
\advance\captionwidth@-1.5in
\def\pagewidth#1{\hsize#1\relax
 \captionwidth@\hsize\advance\captionwidth@-1.5in}
\def\pageheight#1{\vsize#1\relax}
\def\hcorrection#1{\advance\hoffset#1\relax}
\def\vcorrection#1{\advance\voffset#1\relax}
\message{accents/punctuation,}

\let\graveaccent\`
\let\acuteaccent\'
\let\tildeaccent\~
\let\hataccent\^
\let\underscore\_
\let\B\=
\let\D\.
\let\ic@\/
\def\/{\unskip\ic@}
\def\textfonti{\the\textfont\@ne}
\def\t#1#2{{\edef\next@{\the\font}\textfonti\accent"7F \next@#1#2}}
\def~{\unskip\nobreak\ \ignorespaces}
\def\.{.\spacefactor\@m}
\atdef@;{\leavevmode\null;}
\atdef@:{\leavevmode\null:}
\atdef@?{\leavevmode\null?}
\edef\@{\string @}
\def\relaxnext@{\let\next\relax}
\atdef@-{\relaxnext@\leavevmode
 \DN@{\ifx\next-\DN@-{\FN@\nextii@}\else
  \DN@{\leavevmode\hbox{-}}\fi\next@}%
 \DNii@{\ifx\next-\DN@-{\leavevmode\hbox{---}}\else
  \DN@{\leavevmode\hbox{--}}\fi\next@}%
 \FN@\next@}
\def\srdr@{\kern.16667em}
\def\drsr@{\kern.02778em}
\def\sldl@{\drsr@}
\def\dlsl@{\srdr@}
\atdef@"{\unskip\relaxnext@
 \DN@{\ifx\next\space@\DN@. {\FN@\nextii@}\else
  \DN@.{\FN@\nextii@}\fi\next@.}%
 \DNii@{\ifx\next`\DN@`{\FN@\nextiii@}\else
  \ifx\next\lq\DN@\lq{\FN@\nextiii@}\else
  \DN@####1{\FN@\nextiv@}\fi\fi\next@}%
 \def\nextiii@{\ifx\next`\DN@`{\sldl@``}\else\ifx\next\lq
  \DN@\lq{\sldl@``}\else\DN@{\dlsl@`}\fi\fi\next@}%
 \def\nextiv@{\ifx\next'\DN@'{\srdr@''}\else
  \ifx\next\rq\DN@\rq{\srdr@''}\else\DN@{\drsr@'}\fi\fi\next@}%
 \FN@\next@}

\def\textfontii{\the\textfont\tw@}
\def\lbrace@{\delimiter"4266308 }
\def\rbrace@{\delimiter"5267309 }
\def\{{\RIfM@\lbrace@\else{\textfontii f}\spacefactor\@m\fi}
\def\}{\RIfM@\rbrace@\else
 \let\@sf\empty\ifhmode\edef\@sf{\spacefactor\the\spacefactor}\fi
 {\textfontii g}\@sf\relax\fi}
\let\lbrace\{
\let\rbrace\}
\def\AmSTeX{{\textfontii A\kern-.1667em%
  \lower.5ex\hbox{M}\kern-.125emS}-\TeX}
\message{line and page breaks,}
\def\vmodeerr@#1{\Err@{\string#1\space not allowed between paragraphs}}
\def\mathmodeerr@#1{\Err@{\string#1\space not allowed in math mode}}
\def\linebreak{\RIfM@\mathmodeerr@\linebreak\else
 \ifhmode\unskip\unkern\break\else\vmodeerr@\linebreak\fi\fi}

\newskip\saveskip@
\def\allowlinebreak{\RIfM@\mathmodeerr@\allowlinebreak\else
 \ifhmode\saveskip@\lastskip\unskip
 \allowbreak\ifdim\saveskip@>\z@\hskip\saveskip@\fi
 \else\vmodeerr@\allowlinebreak\fi\fi}
\def\nolinebreak{\RIfM@\mathmodeerr@\nolinebreak\else
 \ifhmode\saveskip@\lastskip\unskip
 \nobreak\ifdim\saveskip@>\z@\hskip\saveskip@\fi
 \else\vmodeerr@\nolinebreak\fi\fi}
\def\newline{\relaxnext@
 \DN@{\RIfM@\expandafter\mathmodeerr@\expandafter\newline\else
  \ifhmode\ifx\next\par\else
  \expandafter\unskip\expandafter\null\expandafter\hfill\expandafter\break\fi
  \else
  \expandafter\vmodeerr@\expandafter\newline\fi\fi}%
 \FN@\next@}
\def\dmatherr@#1{\Err@{\string#1\space not allowed in display math mode}}
\def\nondmatherr@#1{\Err@{\string#1\space not allowed in non-display math
 mode}}
\def\onlydmatherr@#1{\Err@{\string#1\space allowed only in display math mode}}
\def\nonmatherr@#1{\Err@{\string#1\space allowed only in math mode}}
\def\mathbreak{\RIfMIfI@\break\else
 \dmatherr@\mathbreak\fi\else\nonmatherr@\mathbreak\fi}
\def\nomathbreak{\RIfMIfI@\nobreak\else
 \dmatherr@\nomathbreak\fi\else\nonmatherr@\nomathbreak\fi}
\def\allowmathbreak{\RIfMIfI@\allowbreak\else
 \dmatherr@\allowmathbreak\fi\else\nonmatherr@\allowmathbreak\fi}
\def\pagebreak{\RIfM@
 \ifinner\nondmatherr@\pagebreak\else\postdisplaypenalty-\@M\fi
 \else\ifvmode\removelastskip\break\else\vadjust{\break}\fi\fi}
\def\nopagebreak{\RIfM@
 \ifinner\nondmatherr@\nopagebreak\else\postdisplaypenalty\@M\fi
 \else\ifvmode\nobreak\else\vadjust{\nobreak}\fi\fi}
\def\nonvmodeerr@#1{\Err@{\string#1\space not allowed within a paragraph
 or in math}}
\def\vnonvmode@#1#2{\relaxnext@\DNii@{\ifx\next\par\DN@{#1}\else
 \DN@{#2}\fi\next@}%
 \ifvmode\DN@{#1}\else
 \DN@{\FN@\nextii@}\fi\next@}
\def\newpage{\vnonvmode@{\vfill\break}{\nonvmodeerr@\newpage}}
\def\smallpagebreak{\vnonvmode@\smallbreak{\nonvmodeerr@\smallpagebreak}}
\def\medpagebreak{\vnonvmode@\medbreak{\nonvmodeerr@\medpagebreak}}
\def\bigpagebreak{\vnonvmode@\bigbreak{\nonvmodeerr@\bigpagebreak}}
\def\NoBlackBoxes{\global\overfullrule\z@}
\def\BlackBoxes{\global\overfullrule5\p@}
\def\Invalid@#1{\def#1{\Err@{\Invalid@@\string#1}}}
\def\Invalid@@{Invalid use of }
\message{figures,}
\Invalid@\caption
\Invalid@\captionwidth
\newdimen\smallcaptionwidth@
\def\topspace{\mid@false\ins@}
\def\midspace{\mid@true\ins@}
\newif\ifmid@
\def\captionfont@{}
\def\ins@#1{\relaxnext@\allowbreak
 \smallcaptionwidth@\captionwidth@\gdef\thespace@{#1}%
 \DN@{\ifx\next\space@\DN@. {\FN@\nextii@}\else
  \DN@.{\FN@\nextii@}\fi\next@.}%
 \DNii@{\ifx\next\caption\DN@\caption{\FN@\nextiii@}%
  \else\let\next@\nextiv@\fi\next@}%
 \def\nextiv@{\vnonvmode@
  {\ifmid@\expandafter\midinsert\else\expandafter\topinsert\fi
   \vbox to\thespace@{}\endinsert}
  {\ifmid@\nonvmodeerr@\midspace\else\nonvmodeerr@\topspace\fi}}%
 \def\nextiii@{\ifx\next\captionwidth\expandafter\nextv@
  \else\expandafter\nextvi@\fi}%
 \def\nextv@\captionwidth##1##2{\smallcaptionwidth@##1\relax\nextvi@{##2}}%
 \def\nextvi@##1{\def\thecaption@{\captionfont@##1}%
  \DN@{\ifx\next\space@\DN@. {\FN@\nextvii@}\else
   \DN@.{\FN@\nextvii@}\fi\next@.}%
  \FN@\next@}%
 \def\nextvii@{\vnonvmode@
  {\ifmid@\expandafter\midinsert\else
  \expandafter\topinsert\fi\vbox to\thespace@{}\nobreak\smallskip
  \setboxz@h{\noindent\ignorespaces\thecaption@\unskip}%
  \ifdim\wdz@>\smallcaptionwidth@\centerline{\vbox{\hsize\smallcaptionwidth@
   \noindent\ignorespaces\thecaption@\unskip}}%
  \else\centerline{\boxz@}\fi\endinsert}
  {\ifmid@\nonvmodeerr@\midspace
  \else\nonvmodeerr@\topspace\fi}}%
 \FN@\next@}
\message{comments,}
\def\newcodes@{\catcode`\\12\catcode`\{12\catcode`\}12\catcode`\#12%
 \catcode`\%12\relax}
\def\oldcodes@{\catcode`\\0\catcode`\{1\catcode`\}2\catcode`\#6%
 \catcode`\%14\relax}
\def\comment{\newcodes@\endlinechar=10 \comment@}
{\lccode`\0=`\\
\lowercase{\gdef\comment@#1^^J{\comment@@#10endcomment\comment@@@}%
\gdef\comment@@#10endcomment{\FN@\comment@@@}%
\gdef\comment@@@#1\comment@@@{\ifx\next\comment@@@\let\next\comment@
 \else\def\next{\oldcodes@\endlinechar=`\^^M\relax}%
 \fi\next}}}
\def\pr@m@s{\ifx'\next\DN@##1{\prim@s}\else\let\next@\egroup\fi\next@}
\def\prime{{\null\prime@\null}}
\mathchardef\prime@="0230
\let\dsize\displaystyle

\let\ssize\scriptstyle

\message{math spacing,}
\def\,{\RIfM@\mskip\thinmuskip\relax\else\kern.16667em\fi}
\def\!{\RIfM@\mskip-\thinmuskip\relax\else\kern-.16667em\fi}
\let\thinspace\,
\let\negthinspace\!
\def\medspace{\RIfM@\mskip\medmuskip\relax\else\kern.222222em\fi}
\def\negmedspace{\RIfM@\mskip-\medmuskip\relax\else\kern-.222222em\fi}
\def\thickspace{\RIfM@\mskip\thickmuskip\relax\else\kern.27777em\fi}
\let\;\thickspace
\def\negthickspace{\RIfM@\mskip-\thickmuskip\relax\else
 \kern-.27777em\fi}
\atdef@,{\RIfM@\mskip.1\thinmuskip\else\leavevmode\null,\fi}
\atdef@!{\RIfM@\mskip-.1\thinmuskip\else\leavevmode\null!\fi}
\atdef@.{\RIfM@&&\else\leavevmode.\spacefactor3000 \fi}
\def\and{\DOTSB\;\mathchar"3026 \;}

\message{fractions,}
\def\frac#1#2{{#1\over#2}}

\newdimen\ex@
\ex@.2326ex
\Invalid@\thickness
\def\thickfrac{\relaxnext@
 \DN@{\ifx\next\thickness\let\next@\nextii@\else
 \DN@{\nextii@\thickness1}\fi\next@}%
 \DNii@\thickness##1##2##3{{##2\above##1\ex@##3}}%
 \FN@\next@}

\def\thickfracwithdelims#1#2{\relaxnext@\def\ldelim@{#1}\def\rdelim@{#2}%
 \DN@{\ifx\next\thickness\let\next@\nextii@\else
 \DN@{\nextii@\thickness1}\fi\next@}%
 \DNii@\thickness##1##2##3{{##2\abovewithdelims
 \ldelim@\rdelim@##1\ex@##3}}%
 \FN@\next@}

\def\:{\nobreak\hskip.1111em\mathpunct{}\nonscript\mkern-\thinmuskip{:}\hskip
 .3333emplus.0555em\relax}
\def\snug{\unskip\kern-\mathsurround}
\message{smash commands,}
\def\topsmash{\top@true\bot@false\smash@}
\def\botsmash{\top@false\bot@true\smash@}
\newif\iftop@
\newif\ifbot@
\def\smash{\top@true\bot@true\smash@}
\def\smash@{\RIfM@\expandafter\mathpalette\expandafter\mathsm@sh\else
 \expandafter\makesm@sh\fi}
\def\finsm@sh{\iftop@\ht\z@\z@\fi\ifbot@\dp\z@\z@\fi\leavevmode\boxz@}
\message{large operator symbols,}
\def\LimitsOnSums{\global\let\slimits@\displaylimits}
\def\NoLimitsOnSums{\global\let\slimits@\nolimits}
\LimitsOnSums
\mathchardef\coprod@="1360       \def\coprod{\DOTSB\coprod@\slimits@}
\mathchardef\bigvee@="1357       \def\bigvee{\DOTSB\bigvee@\slimits@}
\mathchardef\bigwedge@="1356     \def\bigwedge{\DOTSB\bigwedge@\slimits@}
\mathchardef\biguplus@="1355     \def\biguplus{\DOTSB\biguplus@\slimits@}
\mathchardef\bigcap@="1354       \def\bigcap{\DOTSB\bigcap@\slimits@}
\mathchardef\bigcup@="1353       \def\bigcup{\DOTSB\bigcup@\slimits@}
\mathchardef\prod@="1351         \def\prod{\DOTSB\prod@\slimits@}
\mathchardef\sum@="1350          \def\sum{\DOTSB\sum@\slimits@}
\mathchardef\bigotimes@="134E    \def\bigotimes{\DOTSB\bigotimes@\slimits@}
\mathchardef\bigoplus@="134C     \def\bigoplus{\DOTSB\bigoplus@\slimits@}
\mathchardef\bigodot@="134A      \def\bigodot{\DOTSB\bigodot@\slimits@}
\mathchardef\bigsqcup@="1346     \def\bigsqcup{\DOTSB\bigsqcup@\slimits@}
\message{integrals,}
\def\LimitsOnInts{\global\let\ilimits@\displaylimits}
\def\NoLimitsOnInts{\global\let\ilimits@\nolimits}
\NoLimitsOnInts
\def\int{\DOTSI\intop\ilimits@}
\def\oint{\DOTSI\ointop\ilimits@}
\def\intic@{\mathchoice{\hskip.5em}{\hskip.4em}{\hskip.4em}{\hskip.4em}}
\def\negintic@{\mathchoice
 {\hskip-.5em}{\hskip-.4em}{\hskip-.4em}{\hskip-.4em}}
\def\intkern@{\mathchoice{\!\!\!}{\!\!}{\!\!}{\!\!}}
\def\intdots@{\mathchoice{\plaincdots@}
 {{\cdotp}\mkern1.5mu{\cdotp}\mkern1.5mu{\cdotp}}
 {{\cdotp}\mkern1mu{\cdotp}\mkern1mu{\cdotp}}
 {{\cdotp}\mkern1mu{\cdotp}\mkern1mu{\cdotp}}}
\newcount\intno@
\def\iint{\DOTSI\intno@\tw@\FN@\ints@}
\def\iiint{\DOTSI\intno@\thr@@\FN@\ints@}
\def\iiiint{\DOTSI\intno@4 \FN@\ints@}
\def\idotsint{\DOTSI\intno@\z@\FN@\ints@}
\def\ints@{\findlimits@\ints@@}
\newif\iflimtoken@
\newif\iflimits@
\def\findlimits@{\limtoken@true\ifx\next\limits\limits@true
 \else\ifx\next\nolimits\limits@false\else
 \limtoken@false\ifx\ilimits@\nolimits\limits@false\else
 \ifinner\limits@false\else\limits@true\fi\fi\fi\fi}
\def\multint@{\int\ifnum\intno@=\z@\intdots@                                
 \else\intkern@\fi                                                          
 \ifnum\intno@>\tw@\int\intkern@\fi                                         
 \ifnum\intno@>\thr@@\int\intkern@\fi                                       
 \int}                                                                      
\def\multintlimits@{\intop\ifnum\intno@=\z@\intdots@\else\intkern@\fi
 \ifnum\intno@>\tw@\intop\intkern@\fi
 \ifnum\intno@>\thr@@\intop\intkern@\fi\intop}
\def\ints@@{\iflimtoken@                                                    
 \def\ints@@@{\iflimits@\negintic@\mathop{\intic@\multintlimits@}\limits    
  \else\multint@\nolimits\fi                                                
  \eat@}                                                                    
 \else                                                                      
 \def\ints@@@{\iflimits@\negintic@
  \mathop{\intic@\multintlimits@}\limits\else
  \multint@\nolimits\fi}\fi\ints@@@}
\def\LimitsOnNames{\global\let\nlimits@\displaylimits}
\def\NoLimitsOnNames{\global\let\nlimits@\nolimits@}
\LimitsOnNames
\def\nolimits@{\relaxnext@
 \DN@{\ifx\next\limits\DN@\limits{\nolimits}\else
  \let\next@\nolimits\fi\next@}%
 \FN@\next@}
\message{operator names,}
\def\newmcodes@{\mathcode`\'"27\mathcode`\*"2A\mathcode`\."613A%
 \mathcode`\-"2D\mathcode`\/"2F\mathcode`\:"603A }
\def\operatorname#1{\mathop{\newmcodes@\kern\z@\fam\z@#1}\nolimits@}
\def\operatornamewithlimits#1{\mathop{\newmcodes@\kern\z@\fam\z@#1}\nlimits@}
\def\qopname@#1{\mathop{\fam\z@#1}\nolimits@}
\def\qopnamewl@#1{\mathop{\fam\z@#1}\nlimits@}
\def\arccos{\qopname@{arccos}}
\def\arcsin{\qopname@{arcsin}}
\def\arctan{\qopname@{arctan}}
\def\arg{\qopname@{arg}}
\def\cos{\qopname@{cos}}
\def\cosh{\qopname@{cosh}}
\def\cot{\qopname@{cot}}
\def\coth{\qopname@{coth}}
\def\csc{\qopname@{csc}}
\def\deg{\qopname@{deg}}
\def\det{\qopnamewl@{det}}
\def\dim{\qopname@{dim}}
\def\exp{\qopname@{exp}}
\def\gcd{\qopnamewl@{gcd}}
\def\hom{\qopname@{hom}}
\def\inf{\qopnamewl@{inf}}
\def\injlim{\qopnamewl@{inj\,lim}}
\def\ker{\qopname@{ker}}
\def\lg{\qopname@{lg}}
\def\lim{\qopnamewl@{lim}}
\def\liminf{\qopnamewl@{lim\,inf}}
\def\limsup{\qopnamewl@{lim\,sup}}
\def\ln{\qopname@{ln}}
\def\log{\qopname@{log}}
\def\max{\qopnamewl@{max}}
\def\min{\qopnamewl@{min}}
\def\Pr{\qopnamewl@{Pr}}
\def\projlim{\qopnamewl@{proj\,lim}}
\def\sec{\qopname@{sec}}
\def\sin{\qopname@{sin}}
\def\sinh{\qopname@{sinh}}
\def\sup{\qopnamewl@{sup}}
\def\tan{\qopname@{tan}}
\def\tanh{\qopname@{tanh}}
\def\varinjlim{\mathop{\vtop{\ialign{##\crcr
 \hfil\rm lim\hfil\crcr\noalign{\nointerlineskip}\rightarrowfill\crcr
 \noalign{\nointerlineskip\kern-\ex@}\crcr}}}}
\def\varprojlim{\mathop{\vtop{\ialign{##\crcr
 \hfil\rm lim\hfil\crcr\noalign{\nointerlineskip}\leftarrowfill\crcr
 \noalign{\nointerlineskip\kern-\ex@}\crcr}}}}
\def\varliminf{\mathop{\underline{\vrule height\z@ depth.2exwidth\z@
 \hbox{\rm lim}}}}

\newdimen\buffer@
\buffer@\fontdimen13 \tenex
\newdimen\buffer
\buffer\buffer@

\def\ResetBuffer{\fontdimen13 \tenex\buffer@\global\buffer\buffer@}
\def\shave#1{\mathop{\hbox{$\m@th\fontdimen13 \tenex\z@                     
 \displaystyle{#1}$}}\fontdimen13 \tenex\buffer}

\message{multilevel sub/superscripts,}
\Invalid@\\
\def\Let@{\relax\iffalse{\fi\let\\=\cr\iffalse}\fi}
\Invalid@\vspace
\def\vspace@{\def\vspace##1{\crcr\noalign{\vskip##1\relax}}}
\def\multilimits@{\bgroup\vspace@\Let@
 \baselineskip\fontdimen10 \scriptfont\tw@
 \advance\baselineskip\fontdimen12 \scriptfont\tw@
 \lineskip\thr@@\fontdimen8 \scriptfont\thr@@
 \lineskiplimit\lineskip
 \vbox\bgroup\ialign\bgroup\hfil$\m@th\scriptstyle{##}$\hfil\crcr}
\def\Sb{_\multilimits@}
\def\endSb{\crcr\egroup\egroup\egroup}
\def\Sp{^\multilimits@}

\def\spreadlines#1{\RIfMIfI@\onlydmatherr@\spreadlines\else
 \openup#1\relax\fi\else\onlydmatherr@\spreadlines\fi}
\def\Mathstrut@{\copy\Mathstrutbox@}
\newbox\Mathstrutbox@
\setbox\Mathstrutbox@\null
\setboxz@h{$\m@th($}
\ht\Mathstrutbox@\ht\z@
\dp\Mathstrutbox@\dp\z@
\message{matrices,}
\newdimen\spreadmlines@
\def\spreadmatrixlines#1{\RIfMIfI@
 \onlydmatherr@\spreadmatrixlines\else
 \spreadmlines@#1\relax\fi\else\onlydmatherr@\spreadmatrixlines\fi}
\def\matrix{\null\,\vcenter\bgroup\Let@\vspace@
 \normalbaselines\openup\spreadmlines@\ialign
 \bgroup\hfil$\m@th##$\hfil&&\quad\hfil$\m@th##$\hfil\crcr
 \Mathstrut@\crcr\noalign{\kern-\baselineskip}}
\def\endmatrix{\crcr\Mathstrut@\crcr\noalign{\kern-\baselineskip}\egroup
 \egroup\,}
\def\format{\crcr\egroup\iffalse{\fi\ifnum`}=0 \fi\format@}
\newtoks\hashtoks@
\hashtoks@{#}
\def\format@#1\\{\def\preamble@{#1}%
 \def\l{$\m@th\the\hashtoks@$\hfil}%
 \def\c{\hfil$\m@th\the\hashtoks@$\hfil}%
 \def\r{\hfil$\m@th\the\hashtoks@$}%
 \edef\preamble@@{\preamble@}\ifnum`{=0 \fi\iffalse}\fi
 \ialign\bgroup\span\preamble@@\crcr}
\def\smallmatrix{\null\,\vcenter\bgroup\vspace@\Let@
 \baselineskip9\ex@\lineskip\ex@
 \ialign\bgroup\hfil$\m@th\scriptstyle{##}$\hfil&&\thickspace\hfil
 $\m@th\scriptstyle{##}$\hfil\crcr}
\def\endsmallmatrix{\crcr\egroup\egroup\,}

\newmuskip\dotsspace@
\dotsspace@1.5mu
\def\strip@#1 {#1}
\def\spacehdots#1\for#2{\multispan{#2}\xleaders
 \hbox{$\m@th\mkern\strip@#1 \dotsspace@.\mkern\strip@#1 \dotsspace@$}\hfill}
\def\hdotsfor#1{\spacehdots\@ne\for{#1}}
\def\multispan@#1{\omit\mscount#1\unskip\loop\ifnum\mscount>\@ne\sp@n\repeat}
\def\spaceinnerhdots#1\for#2\after#3{\multispan@{\strip@#2 }#3\xleaders
 \hbox{$\m@th\mkern\strip@#1 \dotsspace@.\mkern\strip@#1 \dotsspace@$}\hfill}
\def\innerhdotsfor#1\after#2{\spaceinnerhdots\@ne\for#1\after{#2}}
\def\cases{\bgroup\spreadmlines@\jot\left\{\,\matrix\format\l&\quad\l\\}
\def\endcases{\endmatrix\right.\egroup}
\message{multiline displays,}
\newif\ifinany@
\newif\ifinalign@
\newif\ifingather@
\def\strut@{\copy\strutbox@}
\newbox\strutbox@
\setbox\strutbox@\hbox{\vrule height8\p@ depth3\p@ width\z@}
\def\topaligned{\null\,\vtop\aligned@}
\def\botaligned{\null\,\vbox\aligned@}
\def\aligned{\null\,\vcenter\aligned@}
\def\aligned@{\bgroup\vspace@\Let@
 \ifinany@\else\openup\jot\fi\ialign
 \bgroup\hfil\strut@$\m@th\displaystyle{##}$&
 $\m@th\displaystyle{{}##}$\hfil\crcr}
\def\endaligned{\crcr\egroup\egroup}

\def\alignedat#1{\null\,\vcenter\bgroup\doat@{#1}\vspace@\Let@
 \ifinany@\else\openup\jot\fi\ialign\bgroup\span\preamble@@\crcr}
\newcount\atcount@
\def\doat@#1{\toks@{\hfil\strut@$\m@th
 \displaystyle{\the\hashtoks@}$&$\m@th\displaystyle
 {{}\the\hashtoks@}$\hfil}
 \atcount@#1\relax\advance\atcount@\m@ne                                    
 \loop\ifnum\atcount@>\z@\toks@=\expandafter{\the\toks@&\hfil$\m@th
 \displaystyle{\the\hashtoks@}$&$\m@th
 \displaystyle{{}\the\hashtoks@}$\hfil}\advance
  \atcount@\m@ne\repeat                                                     
 \xdef\preamble@{\the\toks@}\xdef\preamble@@{\preamble@}}

\def\gathered{\null\,\vcenter\bgroup\vspace@\Let@
 \ifinany@\else\openup\jot\fi\ialign
 \bgroup\hfil\strut@$\m@th\displaystyle{##}$\hfil\crcr}
\def\endgathered{\crcr\egroup\egroup}
\newif\iftagsleft@
\def\TagsOnLeft{\global\tagsleft@true}
\def\TagsOnRight{\global\tagsleft@false}
\TagsOnLeft
\newif\ifmathtags@
\def\TagsAsMath{\global\mathtags@true}
\def\TagsAsText{\global\mathtags@false}
\TagsAsText
\def\tagform@#1{\hbox{\rm(\ignorespaces#1\unskip)}}
\def\thetag{\leavevmode\tagform@}
\def\tag#1$${\iftagsleft@\leqno\else\eqno\fi                                
 \maketag@#1\maketag@                                                       
 $$}                                                                        
\def\maketag@{\FN@\maketag@@}
\def\maketag@@{\ifx\next"\expandafter\maketag@@@\else\expandafter\maketag@@@@
 \fi}
\def\maketag@@@"#1"#2\maketag@{\hbox{\rm#1}}                                
\def\maketag@@@@#1\maketag@{\ifmathtags@\tagform@{$\m@th#1$}\else
 \tagform@{#1}\fi}
\interdisplaylinepenalty\@M
\def\allowdisplaybreaks{\RIfMIfI@
 \onlydmatherr@\allowdisplaybreaks\else
 \interdisplaylinepenalty\z@\fi\else\onlydmatherr@\allowdisplaybreaks\fi}
\Invalid@\allowdisplaybreak
\Invalid@\displaybreak
\Invalid@\intertext
\def\allowdisplaybreak@{\def\allowdisplaybreak{\crcr\noalign{\allowbreak}}}
\def\displaybreak@{\def\displaybreak{\crcr\noalign{\break}}}
\def\intertext@{\def\intertext##1{\crcr\noalign{%
 \penalty\postdisplaypenalty \vskip\belowdisplayskip
 \vbox{\normalbaselines\noindent##1}%
 \penalty\predisplaypenalty \vskip\abovedisplayskip}}}
\newskip\centering@
\centering@\z@ plus\@m\p@
\def\align{\relax\ifingather@\DN@{\csname align (in
  \string\gather)\endcsname}\else
 \ifmmode\ifinner\DN@{\onlydmatherr@\align}\else
  \let\next@\align@\fi
 \else\DN@{\onlydmatherr@\align}\fi\fi\next@}
\newhelp\andhelp@
{An extra & here is so disastrous that you should probably exit^^J
and fix things up.}
\newif\iftag@
\newcount\and@
\def\align@{\inalign@true\inany@true
 \vspace@\allowdisplaybreak@\displaybreak@\intertext@
 \def\tag{\global\tag@true\ifnum\and@=\z@\DN@{&&}\else
          \DN@{&}\fi\next@}%
 \iftagsleft@\DN@{\csname align \endcsname}\else
  \DN@{\csname align \space\endcsname}\fi\next@}
\def\Tag@{\iftag@\else\errhelp\andhelp@\err@{Extra & on this line}\fi}
\newdimen\lwidth@
\newdimen\rwidth@
\newdimen\maxlwidth@
\newdimen\maxrwidth@
\newdimen\totwidth@
\def\measure@#1\endalign{\lwidth@\z@\rwidth@\z@\maxlwidth@\z@\maxrwidth@\z@
 \global\and@\z@                                                            
 \setbox@ne\vbox                                                            
  {\everycr{\noalign{\global\tag@false\global\and@\z@}}\Let@                
  \halign{\setboxz@h{$\m@th\displaystyle{\@lign##}$}
   \global\lwidth@\wdz@                                                     
   \ifdim\lwidth@>\maxlwidth@\global\maxlwidth@\lwidth@\fi                  
   \global\advance\and@\@ne                                                 
   &\setboxz@h{$\m@th\displaystyle{{}\@lign##}$}\global\rwidth@\wdz@        
   \ifdim\rwidth@>\maxrwidth@\global\maxrwidth@\rwidth@\fi                  
   \global\advance\and@\@ne                                                
   &\Tag@
   \eat@{##}\crcr#1\crcr}}
 \totwidth@\maxlwidth@\advance\totwidth@\maxrwidth@}                       
\def\displ@y@{\global\dt@ptrue\openup\jot
 \everycr{\noalign{\global\tag@false\global\and@\z@\ifdt@p\global\dt@pfalse
 \vskip-\lineskiplimit\vskip\normallineskiplimit\else
 \penalty\interdisplaylinepenalty\fi}}}
\def\black@#1{\noalign{\ifdim#1>\displaywidth
 \dimen@\prevdepth\nointerlineskip                                          
 \vskip-\ht\strutbox@\vskip-\dp\strutbox@                                   
 \vbox{\noindent\hbox to#1{\strut@\hfill}}
 \prevdepth\dimen@                                                          
 \fi}}
\expandafter\def\csname align \space\endcsname#1\endalign
 {\measure@#1\endalign\global\and@\z@                                       
 \ifingather@\everycr{\noalign{\global\and@\z@}}\else\displ@y@\fi           
 \Let@\tabskip\centering@                                                   
 \halign to\displaywidth
  {\hfil\strut@\setboxz@h{$\m@th\displaystyle{\@lign##}$}
  \global\lwidth@\wdz@\boxz@\global\advance\and@\@ne                        
  \tabskip\z@skip                                                           
  &\setboxz@h{$\m@th\displaystyle{{}\@lign##}$}
  \global\rwidth@\wdz@\boxz@\hfill\global\advance\and@\@ne                  
  \tabskip\centering@                                                       
  &\setboxz@h{\@lign\strut@\maketag@##\maketag@}
  \dimen@\displaywidth\advance\dimen@-\totwidth@
  \divide\dimen@\tw@\advance\dimen@\maxrwidth@\advance\dimen@-\rwidth@     
  \ifdim\dimen@<\tw@\wdz@\llap{\vtop{\normalbaselines\null\boxz@}}
  \else\llap{\boxz@}\fi                                                    
  \tabskip\z@skip                                                          
  \crcr#1\crcr                                                             
  \black@\totwidth@}}                                                      
\newdimen\lineht@
\expandafter\def\csname align \endcsname#1\endalign{\measure@#1\endalign
 \global\and@\z@
 \ifdim\totwidth@>\displaywidth\let\displaywidth@\totwidth@\else
  \let\displaywidth@\displaywidth\fi                                        
 \ifingather@\everycr{\noalign{\global\and@\z@}}\else\displ@y@\fi
 \Let@\tabskip\centering@\halign to\displaywidth
  {\hfil\strut@\setboxz@h{$\m@th\displaystyle{\@lign##}$}%
  \global\lwidth@\wdz@\global\lineht@\ht\z@                                 
  \boxz@\global\advance\and@\@ne
  \tabskip\z@skip&\setboxz@h{$\m@th\displaystyle{{}\@lign##}$}%
  \global\rwidth@\wdz@\ifdim\ht\z@>\lineht@\global\lineht@\ht\z@\fi         
  \boxz@\hfil\global\advance\and@\@ne
  \tabskip\centering@&\kern-\displaywidth@                                  
  \setboxz@h{\@lign\strut@\maketag@##\maketag@}%
  \dimen@\displaywidth\advance\dimen@-\totwidth@
  \divide\dimen@\tw@\advance\dimen@\maxlwidth@\advance\dimen@-\lwidth@
  \ifdim\dimen@<\tw@\wdz@
   \rlap{\vbox{\normalbaselines\boxz@\vbox to\lineht@{}}}\else
   \rlap{\boxz@}\fi
  \tabskip\displaywidth@\crcr#1\crcr\black@\totwidth@}}
\expandafter\def\csname align (in \string\gather)\endcsname
  #1\endalign{\vcenter{\align@#1\endalign}}
\Invalid@\endalign
\newif\ifxat@
\def\alignat{\RIfMIfI@\DN@{\onlydmatherr@\alignat}\else
 \DN@{\csname alignat \endcsname}\fi\else
 \DN@{\onlydmatherr@\alignat}\fi\next@}
\newif\ifmeasuring@
\newbox\savealignat@
\expandafter\def\csname alignat \endcsname#1#2\endalignat                   
 {\inany@true\xat@false
 \def\tag{\global\tag@true\count@#1\relax\multiply\count@\tw@
  \xdef\tag@{}\loop\ifnum\count@>\and@\xdef\tag@{&\tag@}\advance\count@\m@ne
  \repeat\tag@}%
 \vspace@\allowdisplaybreak@\displaybreak@\intertext@
 \displ@y@\measuring@true                                                   
 \setbox\savealignat@\hbox{$\m@th\displaystyle\Let@
  \attag@{#1}
  \vbox{\halign{\span\preamble@@\crcr#2\crcr}}$}%
 \measuring@false                                                           
 \Let@\attag@{#1}
 \tabskip\centering@\halign to\displaywidth
  {\span\preamble@@\crcr#2\crcr                                             
  \black@{\wd\savealignat@}}}                                               
\Invalid@\endalignat
\def\xalignat{\RIfMIfI@
 \DN@{\onlydmatherr@\xalignat}\else
 \DN@{\csname xalignat \endcsname}\fi\else
 \DN@{\onlydmatherr@\xalignat}\fi\next@}
\expandafter\def\csname xalignat \endcsname#1#2\endxalignat
 {\inany@true\xat@true
 \def\tag{\global\tag@true\def\tag@{}\count@#1\relax\multiply\count@\tw@
  \loop\ifnum\count@>\and@\xdef\tag@{&\tag@}\advance\count@\m@ne\repeat\tag@}%
 \vspace@\allowdisplaybreak@\displaybreak@\intertext@
 \displ@y@\measuring@true\setbox\savealignat@\hbox{$\m@th\displaystyle\Let@
 \attag@{#1}\vbox{\halign{\span\preamble@@\crcr#2\crcr}}$}%
 \measuring@false\Let@
 \attag@{#1}\tabskip\centering@\halign to\displaywidth
 {\span\preamble@@\crcr#2\crcr\black@{\wd\savealignat@}}}
\def\attag@#1{\let\Maketag@\maketag@\let\TAG@\Tag@                          
 \let\Tag@=0\let\maketag@=0
 \ifmeasuring@\def\llap@##1{\setboxz@h{##1}\hbox to\tw@\wdz@{}}%
  \def\rlap@##1{\setboxz@h{##1}\hbox to\tw@\wdz@{}}\else
  \let\llap@\llap\let\rlap@\rlap\fi                                         
 \toks@{\hfil\strut@$\m@th\displaystyle{\@lign\the\hashtoks@}$\tabskip\z@skip
  \global\advance\and@\@ne&$\m@th\displaystyle{{}\@lign\the\hashtoks@}$\hfil
  \ifxat@\tabskip\centering@\fi\global\advance\and@\@ne}
 \iftagsleft@
  \toks@@{\tabskip\centering@&\Tag@\kern-\displaywidth
   \rlap@{\@lign\maketag@\the\hashtoks@\maketag@}%
   \global\advance\and@\@ne\tabskip\displaywidth}\else
  \toks@@{\tabskip\centering@&\Tag@\llap@{\@lign\maketag@
   \the\hashtoks@\maketag@}\global\advance\and@\@ne\tabskip\z@skip}\fi      
 \atcount@#1\relax\advance\atcount@\m@ne
 \loop\ifnum\atcount@>\z@
 \toks@=\expandafter{\the\toks@&\hfil$\m@th\displaystyle{\@lign
  \the\hashtoks@}$\global\advance\and@\@ne
  \tabskip\z@skip&$\m@th\displaystyle{{}\@lign\the\hashtoks@}$\hfil\ifxat@
  \tabskip\centering@\fi\global\advance\and@\@ne}\advance\atcount@\m@ne
 \repeat                                                                    
 \xdef\preamble@{\the\toks@\the\toks@@}
 \xdef\preamble@@{\preamble@}
 \let\maketag@\Maketag@\let\Tag@\TAG@}                                      
\Invalid@\endxalignat
\def\xxalignat{\RIfMIfI@
 \DN@{\onlydmatherr@\xxalignat}\else\DN@{\csname xxalignat
  \endcsname}\fi\else
 \DN@{\onlydmatherr@\xxalignat}\fi\next@}
\expandafter\def\csname xxalignat \endcsname#1#2\endxxalignat{\inany@true
 \vspace@\allowdisplaybreak@\displaybreak@\intertext@
 \displ@y\setbox\savealignat@\hbox{$\m@th\displaystyle\Let@
 \xxattag@{#1}\vbox{\halign{\span\preamble@@\crcr#2\crcr}}$}%
 \Let@\xxattag@{#1}\tabskip\z@skip\halign to\displaywidth
 {\span\preamble@@\crcr#2\crcr\black@{\wd\savealignat@}}}
\def\xxattag@#1{\toks@{\tabskip\z@skip\hfil\strut@
 $\m@th\displaystyle{\the\hashtoks@}$&%
 $\m@th\displaystyle{{}\the\hashtoks@}$\hfil\tabskip\centering@&}%
 \atcount@#1\relax\advance\atcount@\m@ne\loop\ifnum\atcount@>\z@
 \toks@=\expandafter{\the\toks@&\hfil$\m@th\displaystyle{\the\hashtoks@}$%
  \tabskip\z@skip&$\m@th\displaystyle{{}\the\hashtoks@}$\hfil
  \tabskip\centering@}\advance\atcount@\m@ne\repeat
 \xdef\preamble@{\the\toks@\tabskip\z@skip}\xdef\preamble@@{\preamble@}}
\Invalid@\endxxalignat
\newdimen\gwidth@
\newdimen\gmaxwidth@
\def\gmeasure@#1\endgather{\gwidth@\z@\gmaxwidth@\z@\setbox@ne\vbox{\Let@
 \halign{\setboxz@h{$\m@th\displaystyle{##}$}\global\gwidth@\wdz@
 \ifdim\gwidth@>\gmaxwidth@\global\gmaxwidth@\gwidth@\fi
 &\eat@{##}\crcr#1\crcr}}}
\def\gather{\RIfMIfI@\DN@{\onlydmatherr@\gather}\else
 \ingather@true\inany@true\def\tag{&}%
 \vspace@\allowdisplaybreak@\displaybreak@\intertext@
 \displ@y\Let@
 \iftagsleft@\DN@{\csname gather \endcsname}\else
  \DN@{\csname gather \space\endcsname}\fi\fi
 \else\DN@{\onlydmatherr@\gather}\fi\next@}
\expandafter\def\csname gather \space\endcsname#1\endgather
 {\gmeasure@#1\endgather\tabskip\centering@
 \halign to\displaywidth{\hfil\strut@\setboxz@h{$\m@th\displaystyle{##}$}%
 \global\gwidth@\wdz@\boxz@\hfil&
 \setboxz@h{\strut@{\maketag@##\maketag@}}%
 \dimen@\displaywidth\advance\dimen@-\gwidth@
 \ifdim\dimen@>\tw@\wdz@\llap{\boxz@}\else
 \llap{\vtop{\normalbaselines\null\boxz@}}\fi
 \tabskip\z@skip\crcr#1\crcr\black@\gmaxwidth@}}
\newdimen\glineht@
\expandafter\def\csname gather \endcsname#1\endgather{\gmeasure@#1\endgather
 \ifdim\gmaxwidth@>\displaywidth\let\gdisplaywidth@\gmaxwidth@\else
 \let\gdisplaywidth@\displaywidth\fi\tabskip\centering@\halign to\displaywidth
 {\hfil\strut@\setboxz@h{$\m@th\displaystyle{##}$}%
 \global\gwidth@\wdz@\global\glineht@\ht\z@\boxz@\hfil&\kern-\gdisplaywidth@
 \setboxz@h{\strut@{\maketag@##\maketag@}}%
 \dimen@\displaywidth\advance\dimen@-\gwidth@
 \ifdim\dimen@>\tw@\wdz@\rlap{\boxz@}\else
 \rlap{\vbox{\normalbaselines\boxz@\vbox to\glineht@{}}}\fi
 \tabskip\gdisplaywidth@\crcr#1\crcr\black@\gmaxwidth@}}
\newif\ifctagsplit@
\def\CenteredTagsOnSplits{\global\ctagsplit@true}
\def\TopOrBottomTagsOnSplits{\global\ctagsplit@false}
\TopOrBottomTagsOnSplits
\def\split{\relax\ifinany@\let\next@\insplit@\else
 \ifmmode\ifinner\def\next@{\onlydmatherr@\split}\else
 \let\next@\outsplit@\fi\else
 \def\next@{\onlydmatherr@\split}\fi\fi\next@}
\def\insplit@{\global\setbox\z@\vbox\bgroup\vspace@\Let@\ialign\bgroup
 \hfil\strut@$\m@th\displaystyle{##}$&$\m@th\displaystyle{{}##}$\hfill\crcr}
\def\endsplit{\crcr\egroup\egroup\iftagsleft@\expandafter\lendsplit@\else
 \expandafter\rendsplit@\fi}
\def\rendsplit@{\global\setbox9 \vbox
 {\unvcopy\z@\global\setbox8 \lastbox\unskip}
 \setbox@ne\hbox{\unhcopy8 \unskip\global\setbox\tw@\lastbox
 \unskip\global\setbox\thr@@\lastbox}
 \global\setbox7 \hbox{\unhbox\tw@\unskip}
 \ifinalign@\ifctagsplit@                                                   
  \gdef\split@{\hbox to\wd\thr@@{}&
   \vcenter{\vbox{\moveleft\wd\thr@@\boxz@}}}
 \else\gdef\split@{&\vbox{\moveleft\wd\thr@@\box9}\crcr
  \box\thr@@&\box7}\fi                                                      
 \else                                                                      
  \ifctagsplit@\gdef\split@{\vcenter{\boxz@}}\else
  \gdef\split@{\box9\crcr\hbox{\box\thr@@\box7}}\fi
 \fi
 \split@}                                                                   
\def\lendsplit@{\global\setbox9\vtop{\unvcopy\z@}
 \setbox@ne\vbox{\unvcopy\z@\global\setbox8\lastbox}
 \setbox@ne\hbox{\unhcopy8\unskip\setbox\tw@\lastbox
  \unskip\global\setbox\thr@@\lastbox}
 \ifinalign@\ifctagsplit@                                                   
  \gdef\split@{\hbox to\wd\thr@@{}&
  \vcenter{\vbox{\moveleft\wd\thr@@\box9}}}
  \else                                                                     
  \gdef\split@{\hbox to\wd\thr@@{}&\vbox{\moveleft\wd\thr@@\box9}}\fi
 \else
  \ifctagsplit@\gdef\split@{\vcenter{\box9}}\else
  \gdef\split@{\box9}\fi
 \fi\split@}
\def\outsplit@#1$${\align\insplit@#1\endalign$$}
\newdimen\multlinegap@
\multlinegap@1em
\newdimen\multlinetaggap@
\multlinetaggap@1em
\def\MultlineGap#1{\global\multlinegap@#1\relax}
\def\multlinegap#1{\RIfMIfI@\onlydmatherr@\multlinegap\else
 \multlinegap@#1\relax\fi\else\onlydmatherr@\multlinegap\fi}
\def\nomultlinegap{\multlinegap{\z@}}
\def\multline{\RIfMIfI@
 \DN@{\onlydmatherr@\multline}\else
 \DN@{\multline@}\fi\else
 \DN@{\onlydmatherr@\multline}\fi\next@}
\newif\iftagin@
\def\tagin@#1{\tagin@false\in@\tag{#1}\ifin@\tagin@true\fi}
\def\multline@#1$${\inany@true\vspace@\allowdisplaybreak@\displaybreak@
 \tagin@{#1}\iftagsleft@\DN@{\multline@l#1$$}\else
 \DN@{\multline@r#1$$}\fi\next@}
\newdimen\mwidth@
\def\rmmeasure@#1\endmultline{%
 \def\shoveleft##1{##1}\def\shoveright##1{##1}
 \setbox@ne\vbox{\Let@\halign{\setboxz@h
  {$\m@th\@lign\displaystyle{}##$}\global\mwidth@\wdz@
  \crcr#1\crcr}}}
\newdimen\mlineht@
\newif\ifzerocr@
\newif\ifonecr@
\def\lmmeasure@#1\endmultline{\global\zerocr@true\global\onecr@false
 \everycr{\noalign{\ifonecr@\global\onecr@false\fi
  \ifzerocr@\global\zerocr@false\global\onecr@true\fi}}
  \def\shoveleft##1{##1}\def\shoveright##1{##1}%
 \setbox@ne\vbox{\Let@\halign{\setboxz@h
  {$\m@th\@lign\displaystyle{}##$}\ifonecr@\global\mwidth@\wdz@
  \global\mlineht@\ht\z@\fi\crcr#1\crcr}}}
\newbox\mtagbox@
\newdimen\ltwidth@
\newdimen\rtwidth@
\def\multline@l#1$${\iftagin@\DN@{\lmultline@@#1$$}\else
 \DN@{\setbox\mtagbox@\null\ltwidth@\z@\rtwidth@\z@
  \lmultline@@@#1$$}\fi\next@}
\def\lmultline@@#1\endmultline\tag#2$${%
 \setbox\mtagbox@\hbox{\maketag@#2\maketag@}
 \lmmeasure@#1\endmultline\dimen@\mwidth@\advance\dimen@\wd\mtagbox@
 \advance\dimen@\multlinetaggap@                                            
 \ifdim\dimen@>\displaywidth\ltwidth@\z@\else\ltwidth@\wd\mtagbox@\fi       
 \lmultline@@@#1\endmultline$$}
\def\lmultline@@@{\displ@y
 \def\shoveright##1{##1\hfilneg\hskip\multlinegap@}%
 \def\shoveleft##1{\setboxz@h{$\m@th\displaystyle{}##1$}%
  \setbox@ne\hbox{$\m@th\displaystyle##1$}%
  \hfilneg
  \iftagin@
   \ifdim\ltwidth@>\z@\hskip\ltwidth@\hskip\multlinetaggap@\fi
  \else\hskip\multlinegap@\fi\hskip.5\wd@ne\hskip-.5\wdz@##1}
  \halign\bgroup\Let@\hbox to\displaywidth
   {\strut@$\m@th\displaystyle\hfil{}##\hfil$}\crcr
   \hfilneg                                                                 
   \iftagin@                                                                
    \ifdim\ltwidth@>\z@                                                     
     \box\mtagbox@\hskip\multlinetaggap@                                    
    \else
     \rlap{\vbox{\normalbaselines\hbox{\strut@\box\mtagbox@}%
     \vbox to\mlineht@{}}}\fi                                               
   \else\hskip\multlinegap@\fi}                                             
\def\multline@r#1$${\iftagin@\DN@{\rmultline@@#1$$}\else
 \DN@{\setbox\mtagbox@\null\ltwidth@\z@\rtwidth@\z@
  \rmultline@@@#1$$}\fi\next@}
\def\rmultline@@#1\endmultline\tag#2$${\ltwidth@\z@
 \setbox\mtagbox@\hbox{\maketag@#2\maketag@}%
 \rmmeasure@#1\endmultline\dimen@\mwidth@\advance\dimen@\wd\mtagbox@
 \advance\dimen@\multlinetaggap@
 \ifdim\dimen@>\displaywidth\rtwidth@\z@\else\rtwidth@\wd\mtagbox@\fi
 \rmultline@@@#1\endmultline$$}
\def\rmultline@@@{\displ@y
 \def\shoveright##1{##1\hfilneg\iftagin@\ifdim\rtwidth@>\z@
  \hskip\rtwidth@\hskip\multlinetaggap@\fi\else\hskip\multlinegap@\fi}%
 \def\shoveleft##1{\setboxz@h{$\m@th\displaystyle{}##1$}%
  \setbox@ne\hbox{$\m@th\displaystyle##1$}%
  \hfilneg\hskip\multlinegap@\hskip.5\wd@ne\hskip-.5\wdz@##1}%
 \halign\bgroup\Let@\hbox to\displaywidth
  {\strut@$\m@th\displaystyle\hfil{}##\hfil$}\crcr
 \hfilneg\hskip\multlinegap@}
\def\endmultline{\iftagsleft@\expandafter\lendmultline@\else
 \expandafter\rendmultline@\fi}
\def\lendmultline@{\hfilneg\hskip\multlinegap@\crcr\egroup}
\def\rendmultline@{\iftagin@                                                
 \ifdim\rtwidth@>\z@                                                        
  \hskip\multlinetaggap@\box\mtagbox@                                       
 \else\llap{\vtop{\normalbaselines\null\hbox{\strut@\box\mtagbox@}}}\fi     
 \else\hskip\multlinegap@\fi                                                
 \hfilneg\crcr\egroup}
\def\bmod{\mskip-\medmuskip\mkern5mu\mathbin{\fam\z@ mod}\penalty900
 \mkern5mu\mskip-\medmuskip}
\def\pmod#1{\allowbreak\ifinner\mkern8mu\else\mkern18mu\fi
 ({\fam\z@ mod}\,\,#1)}
\def\pod#1{\allowbreak\ifinner\mkern8mu\else\mkern18mu\fi(#1)}
\def\mod#1{\allowbreak\ifinner\mkern12mu\else\mkern18mu\fi{\fam\z@ mod}\,\,#1}
\message{continued fractions,}
\newcount\cfraccount@
\def\cfrac{\bgroup\bgroup\advance\cfraccount@\@ne\strut
 \iffalse{\fi\def\\{\over\displaystyle}\iffalse}\fi}
\def\lcfrac{\bgroup\bgroup\advance\cfraccount@\@ne\strut
 \iffalse{\fi\def\\{\hfill\over\displaystyle}\iffalse}\fi}
\def\rcfrac{\bgroup\bgroup\advance\cfraccount@\@ne\strut\hfill
 \iffalse{\fi\def\\{\over\displaystyle}\iffalse}\fi}
\def\gloop@#1\repeat{\gdef\body{#1}\iterate}
\def\endcfrac{\gloop@\ifnum\cfraccount@>\z@\global\advance\cfraccount@\m@ne
 \egroup\hskip-\nulldelimiterspace\egroup\repeat}
\message{compound symbols,}
\def\binrel@#1{\setboxz@h{\thinmuskip0mu
  \medmuskip\m@ne mu\thickmuskip\@ne mu$#1\m@th$}%
 \setbox@ne\hbox{\thinmuskip0mu\medmuskip\m@ne mu\thickmuskip
  \@ne mu${}#1{}\m@th$}%
 \setbox\tw@\hbox{\hskip\wd@ne\hskip-\wdz@}}
\def\overset#1\to#2{\binrel@{#2}\ifdim\wd\tw@<\z@
 \mathbin{\mathop{\kern\z@#2}\limits^{#1}}\else\ifdim\wd\tw@>\z@
 \mathrel{\mathop{\kern\z@#2}\limits^{#1}}\else
 {\mathop{\kern\z@#2}\limits^{#1}}{}\fi\fi}
\def\underset#1\to#2{\binrel@{#2}\ifdim\wd\tw@<\z@
 \mathbin{\mathop{\kern\z@#2}\limits_{#1}}\else\ifdim\wd\tw@>\z@
 \mathrel{\mathop{\kern\z@#2}\limits_{#1}}\else
 {\mathop{\kern\z@#2}\limits_{#1}}{}\fi\fi}
\def\oversetbrace#1\to#2{\overbrace{#2}^{#1}}
\def\undersetbrace#1\to#2{\underbrace{#2}_{#1}}
\def\sideset#1\and#2\to#3{%
 \setbox@ne\hbox{$\dsize{\vphantom{#3}}#1{#3}\m@th$}%
 \setbox\tw@\hbox{$\dsize{#3}#2\m@th$}%
 \hskip\wd@ne\hskip-\wd\tw@\mathop{\hskip\wd\tw@\hskip-\wd@ne
  {\vphantom{#3}}#1{#3}#2}}
\def\rightarrowfill@#1{\setboxz@h{$#1-\m@th$}\ht\z@\z@
  $#1\m@th\copy\z@\mkern-6mu\cleaders
  \hbox{$#1\mkern-2mu\box\z@\mkern-2mu$}\hfill
  \mkern-6mu\mathord\rightarrow$}
\def\leftarrowfill@#1{\setboxz@h{$#1-\m@th$}\ht\z@\z@
  $#1\m@th\mathord\leftarrow\mkern-6mu\cleaders
  \hbox{$#1\mkern-2mu\copy\z@\mkern-2mu$}\hfill
  \mkern-6mu\box\z@$}
\def\leftrightarrowfill@#1{\setboxz@h{$#1-\m@th$}\ht\z@\z@
  $#1\m@th\mathord\leftarrow\mkern-6mu\cleaders
  \hbox{$#1\mkern-2mu\box\z@\mkern-2mu$}\hfill
  \mkern-6mu\mathord\rightarrow$}
\def\overrightarrow{\mathpalette\overrightarrow@}
\def\overrightarrow@#1#2{\vbox{\ialign{##\crcr\rightarrowfill@#1\crcr
 \noalign{\kern-\ex@\nointerlineskip}$\m@th\hfil#1#2\hfil$\crcr}}}

\def\overleftarrow{\mathpalette\overleftarrow@}
\def\overleftarrow@#1#2{\vbox{\ialign{##\crcr\leftarrowfill@#1\crcr
 \noalign{\kern-\ex@\nointerlineskip}$\m@th\hfil#1#2\hfil$\crcr}}}
\def\overleftrightarrow{\mathpalette\overleftrightarrow@}
\def\overleftrightarrow@#1#2{\vbox{\ialign{##\crcr\leftrightarrowfill@#1\crcr
 \noalign{\kern-\ex@\nointerlineskip}$\m@th\hfil#1#2\hfil$\crcr}}}
\def\underrightarrow{\mathpalette\underrightarrow@}
\def\underrightarrow@#1#2{\vtop{\ialign{##\crcr$\m@th\hfil#1#2\hfil$\crcr
 \noalign{\nointerlineskip}\rightarrowfill@#1\crcr}}}

\def\underleftarrow{\mathpalette\underleftarrow@}
\def\underleftarrow@#1#2{\vtop{\ialign{##\crcr$\m@th\hfil#1#2\hfil$\crcr
 \noalign{\nointerlineskip}\leftarrowfill@#1\crcr}}}
\def\underleftrightarrow{\mathpalette\underleftrightarrow@}
\def\underleftrightarrow@#1#2{\vtop{\ialign{##\crcr$\m@th\hfil#1#2\hfil$\crcr
 \noalign{\nointerlineskip}\leftrightarrowfill@#1\crcr}}}
\message{various kinds of dots,}
\let\DOTSI\relax
\let\DOTSB\relax

\newif\ifmath@
{\uccode`7=`\\ \uccode`8=`m \uccode`9=`a \uccode`0=`t \uccode`!=`h
 \uppercase{\gdef\math@#1#2#3#4#5#6\math@{\global\math@false\ifx 7#1\ifx 8#2%
 \ifx 9#3\ifx 0#4\ifx !#5\xdef\meaning@{#6}\global\math@true\fi\fi\fi\fi\fi}}}
\newif\ifmathch@
{\uccode`7=`c \uccode`8=`h \uccode`9=`\"
 \uppercase{\gdef\mathch@#1#2#3#4#5#6\mathch@{\global\mathch@false
  \ifx 7#1\ifx 8#2\ifx 9#5\global\mathch@true\xdef\meaning@{9#6}\fi\fi\fi}}}
\newcount\classnum@
\def\getmathch@#1.#2\getmathch@{\classnum@#1 \divide\classnum@4096
 \ifcase\number\classnum@\or\or\gdef\thedots@{\dotsb@}\or
 \gdef\thedots@{\dotsb@}\fi}
\newif\ifmathbin@
{\uccode`4=`b \uccode`5=`i \uccode`6=`n
 \uppercase{\gdef\mathbin@#1#2#3{\relaxnext@
  \DNii@##1\mathbin@{\ifx\space@\next\global\mathbin@true\fi}%
 \global\mathbin@false\DN@##1\mathbin@{}%
 \ifx 4#1\ifx 5#2\ifx 6#3\DN@{\FN@\nextii@}\fi\fi\fi\next@}}}
\newif\ifmathrel@
{\uccode`4=`r \uccode`5=`e \uccode`6=`l
 \uppercase{\gdef\mathrel@#1#2#3{\relaxnext@
  \DNii@##1\mathrel@{\ifx\space@\next\global\mathrel@true\fi}%
 \global\mathrel@false\DN@##1\mathrel@{}%
 \ifx 4#1\ifx 5#2\ifx 6#3\DN@{\FN@\nextii@}\fi\fi\fi\next@}}}
\newif\ifmacro@
{\uccode`5=`m \uccode`6=`a \uccode`7=`c
 \uppercase{\gdef\macro@#1#2#3#4\macro@{\global\macro@false
  \ifx 5#1\ifx 6#2\ifx 7#3\global\macro@true
  \xdef\meaning@{\macro@@#4\macro@@}\fi\fi\fi}}}
\def\macro@@#1->#2\macro@@{#2}
\newif\ifDOTS@
\newcount\DOTSCASE@
{\uccode`6=`\\ \uccode`7=`D \uccode`8=`O \uccode`9=`T \uccode`0=`S
 \uppercase{\gdef\DOTS@#1#2#3#4#5{\global\DOTS@false\DN@##1\DOTS@{}%
  \ifx 6#1\ifx 7#2\ifx 8#3\ifx 9#4\ifx 0#5\let\next@\DOTS@@\fi\fi\fi\fi\fi
  \next@}}}
{\uccode`3=`B \uccode`4=`I \uccode`5=`X
 \uppercase{\gdef\DOTS@@#1{\relaxnext@
  \DNii@##1\DOTS@{\ifx\space@\next\global\DOTS@true\fi}%
  \DN@{\FN@\nextii@}%
  \ifx 3#1\global\DOTSCASE@\z@\else
  \ifx 4#1\global\DOTSCASE@\@ne\else
  \ifx 5#1\global\DOTSCASE@\tw@\else\DN@##1\DOTS@{}%
  \fi\fi\fi\next@}}}
\newif\ifnot@
{\uccode`5=`\\ \uccode`6=`n \uccode`7=`o \uccode`8=`t
 \uppercase{\gdef\not@#1#2#3#4{\relaxnext@
  \DNii@##1\not@{\ifx\space@\next\global\not@true\fi}%
 \global\not@false\DN@##1\not@{}%
 \ifx 5#1\ifx 6#2\ifx 7#3\ifx 8#4\DN@{\FN@\nextii@}\fi\fi\fi
 \fi\next@}}}
\newif\ifkeybin@
\def\keybin@{\keybin@true
 \ifx\next+\else\ifx\next=\else\ifx\next<\else\ifx\next>\else\ifx\next-\else
 \ifx\next*\else\ifx\next:\else\keybin@false\fi\fi\fi\fi\fi\fi\fi}
\def\dots{\RIfM@\expandafter\mdots@\else\expandafter\tdots@\fi}
\def\tdots@{\unskip\relaxnext@
 \DN@{$\m@th\mathinner{\ldotp\ldotp\ldotp}\,
   \ifx\next,\,$\else\ifx\next.\,$\else\ifx\next;\,$\else\ifx\next:\,$\else
   \ifx\next?\,$\else\ifx\next!\,$\else$ \fi\fi\fi\fi\fi\fi}%
 \ \FN@\next@}
\def\mdots@{\FN@\mdots@@}
\def\mdots@@{\gdef\thedots@{\dotso@}
 \ifx\next\boldkey\gdef\thedots@\boldkey{\boldkeydots@}\else                
 \ifx\next\boldsymbol\gdef\thedots@\boldsymbol{\boldsymboldots@}\else       
 \ifx,\next\gdef\thedots@{\dotsc}
 \else\ifx\not\next\gdef\thedots@{\dotsb@}
 \else\keybin@
 \ifkeybin@\gdef\thedots@{\dotsb@}
 \else\xdef\meaning@{\meaning\next..........}\xdef\meaning@@{\meaning@}
  \expandafter\math@\meaning@\math@
  \ifmath@
   \expandafter\mathch@\meaning@\mathch@
   \ifmathch@\expandafter\getmathch@\meaning@\getmathch@\fi                 
  \else\expandafter\macro@\meaning@@\macro@                                 
  \ifmacro@                                                                
   \expandafter\not@\meaning@\not@\ifnot@\gdef\thedots@{\dotsb@}
  \else\expandafter\DOTS@\meaning@\DOTS@
  \ifDOTS@
   \ifcase\number\DOTSCASE@\gdef\thedots@{\dotsb@}%
    \or\gdef\thedots@{\dotsi}\else\fi                                      
  \else\expandafter\math@\meaning@\math@                                   
  \ifmath@\expandafter\mathbin@\meaning@\mathbin@
  \ifmathbin@\gdef\thedots@{\dotsb@}
  \else\expandafter\mathrel@\meaning@\mathrel@
  \ifmathrel@\gdef\thedots@{\dotsb@}
  \fi\fi\fi\fi\fi\fi\fi\fi\fi\fi\fi\fi
 \thedots@}
\def\plainldots@{\mathinner{\ldotp\ldotp\ldotp}}
\def\plaincdots@{\mathinner{\cdotp\cdotp\cdotp}}
\def\dotsi{\!\plaincdots@}
\let\dotsb@\plaincdots@
\newif\ifextra@
\newif\ifrightdelim@
\def\rightdelim@{\global\rightdelim@true                                    
 \ifx\next)\else                                                            
 \ifx\next]\else
 \ifx\next\rbrack\else
 \ifx\next\}\else
 \ifx\next\rbrace\else
 \ifx\next\rangle\else
 \ifx\next\rceil\else
 \ifx\next\rfloor\else
 \ifx\next\rgroup\else
 \ifx\next\rmoustache\else
 \ifx\next\right\else
 \ifx\next\bigr\else
 \ifx\next\biggr\else
 \ifx\next\Bigr\else                                                        
 \ifx\next\Biggr\else\global\rightdelim@false
 \fi\fi\fi\fi\fi\fi\fi\fi\fi\fi\fi\fi\fi\fi\fi}
\def\extra@{%
 \global\extra@false\rightdelim@\ifrightdelim@\global\extra@true            
 \else\ifx\next$\global\extra@true                                          
 \else\xdef\meaning@{\meaning\next..........}
 \expandafter\macro@\meaning@\macro@\ifmacro@                               
 \expandafter\DOTS@\meaning@\DOTS@
 \ifDOTS@
 \ifnum\DOTSCASE@=\tw@\global\extra@true                                    
 \fi\fi\fi\fi\fi}
\newif\ifbold@
\def\dotso@{\relaxnext@
 \ifbold@
  \let\next\delayed@
  \DNii@{\extra@\plainldots@\ifextra@\,\fi}%
 \else
  \DNii@{\DN@{\extra@\plainldots@\ifextra@\,\fi}\FN@\next@}%
 \fi
 \nextii@}
\def\extrap@#1{%
 \ifx\next,\DN@{#1\,}\else
 \ifx\next;\DN@{#1\,}\else
 \ifx\next.\DN@{#1\,}\else\extra@
 \ifextra@\DN@{#1\,}\else
 \let\next@#1\fi\fi\fi\fi\next@}
\def\ldots{\DN@{\extrap@\plainldots@}%
 \FN@\next@}
\def\cdots{\DN@{\extrap@\plaincdots@}%
 \FN@\next@}

\def\dotsc{\relaxnext@
 \DN@{\ifx\next;\plainldots@\,\else
  \ifx\next.\plainldots@\,\else\extra@\plainldots@
  \ifextra@\,\fi\fi\fi}%
 \FN@\next@}
\def\cdot{\mathchar"2201 }

\def\mapsto{\DOTSB\mapstochar\rightarrow}

\message{special superscripts,}
\def\dddot#1{{\mathop{#1}\limits^{\vbox to-1.4\ex@{\kern-\tw@\ex@
 \hbox{\rm...}\vss}}}}
\def\ddddot#1{{\mathop{#1}\limits^{\vbox to-1.4\ex@{\kern-\tw@\ex@
 \hbox{\rm....}\vss}}}}
\def\sphat{^{\mathchoice{}{}%
 {\,\,\botsmash{\hbox{\lower4\ex@\hbox{$\m@th\widehat{\null}$}}}}%
 {\,\botsmash{\hbox{\lower3\ex@\hbox{$\m@th\hat{\null}$}}}}}}

\def\spacute{^{\!\botsmash{\hbox{\lower\@ne ex\hbox{\'{}}}}}}
\def\spgrave{^{\mathchoice{}{}{}{\!}%
 \botsmash{\hbox{\lower\@ne ex\hbox{\`{}}}}}}
\def\spdot{^{\hbox{\raise\ex@\hbox{\rm.}}}}
\def\spddot{^{\hbox{\raise\ex@\hbox{\rm..}}}}
\def\spdddot{^{\hbox{\raise\ex@\hbox{\rm...}}}}
\def\spddddot{^{\hbox{\raise\ex@\hbox{\rm....}}}}
\def\spbreve{^{\!\botsmash{\hbox{\lower4\ex@\hbox{\u{}}}}}}

\message{\string\text,}
\def\textonlyfont@#1#2{\def#1{\RIfM@
 \Err@{Use \string#1\space only in text}\else#2\fi}}
\textonlyfont@\rm\tenrm
\textonlyfont@\it\tenit
\textonlyfont@\sl\tensl
\textonlyfont@\bf\tenbf
\def\oldnos#1{\RIfM@{\mathcode`\,="013B \fam\@ne#1}\else
 \leavevmode\hbox{$\m@th\mathcode`\,="013B \fam\@ne#1$}\fi}
\def\text{\RIfM@\expandafter\text@\else\expandafter\text@@\fi}
\def\text@@#1{\leavevmode\hbox{#1}}
\def\mathhexbox@#1#2#3{\text{$\m@th\mathchar"#1#2#3$}}
\def\dag{{\mathhexbox@279}}
\def\ddag{{\mathhexbox@27A}}
\def\S{{\mathhexbox@278}}
\def\P{{\mathhexbox@27B}}
\newif\iffirstchoice@
\firstchoice@true
\def\text@#1{\mathchoice
 {\hbox{\everymath{\displaystyle}\def\textfonti{\the\textfont\@ne}%
  \def\textfontii{\the\textfont\tw@}\textdef@@ T#1}}
 {\hbox{\firstchoice@false
  \everymath{\textstyle}\def\textfonti{\the\textfont\@ne}%
  \def\textfontii{\the\textfont\tw@}\textdef@@ T#1}}
 {\hbox{\firstchoice@false
  \everymath{\scriptstyle}\def\textfonti{\the\scriptfont\@ne}%
  \def\textfontii{\the\scriptfont\tw@}\textdef@@ S\rm#1}}
 {\hbox{\firstchoice@false
  \everymath{\scriptscriptstyle}\def\textfonti
  {\the\scriptscriptfont\@ne}%
  \def\textfontii{\the\scriptscriptfont\tw@}\textdef@@ s\rm#1}}}
\def\textdef@@#1{\textdef@#1\rm\textdef@#1\bf\textdef@#1\sl\textdef@#1\it}
\def\rmfam{0}
\def\textdef@#1#2{%
 \DN@{\csname\expandafter\eat@\string#2fam\endcsname}%
 \if S#1\edef#2{\the\scriptfont\next@\relax}%
 \else\if s#1\edef#2{\the\scriptscriptfont\next@\relax}%
 \else\edef#2{\the\textfont\next@\relax}\fi\fi}
\scriptfont\itfam\tenit \scriptscriptfont\itfam\tenit
\scriptfont\slfam\tensl \scriptscriptfont\slfam\tensl
\newif\iftopfolded@
\newif\ifbotfolded@
\def\topfoldedtext{\topfolded@true\botfolded@false\foldedtext@}
\def\botfoldedtext{\botfolded@true\topfolded@false\foldedtext@}
\def\foldedtext{\topfolded@false\botfolded@false\foldedtext@}
\Invalid@\foldedwidth
\def\foldedtext@{\relaxnext@
 \DN@{\ifx\next\foldedwidth\let\next@\nextii@\else
  \DN@{\nextii@\foldedwidth{.3\hsize}}\fi\next@}%
 \DNii@\foldedwidth##1##2{\setbox\z@\vbox
  {\normalbaselines\hsize##1\relax
  \tolerance1600 \noindent\ignorespaces##2}\ifbotfolded@\boxz@\else
  \iftopfolded@\vtop{\unvbox\z@}\else\vcenter{\boxz@}\fi\fi}%
 \FN@\next@}
\message{math font commands,}
\def\bold{\RIfM@\expandafter\bold@\else
 \expandafter\nonmatherr@\expandafter\bold\fi}
\def\bold@#1{{\bold@@{#1}}}
\def\bold@@#1{\fam\bffam\relax#1}
\def\slanted{\RIfM@\expandafter\slanted@\else
 \expandafter\nonmatherr@\expandafter\slanted\fi}
\def\slanted@#1{{\slanted@@{#1}}}
\def\slanted@@#1{\fam\slfam\relax#1}
\def\roman{\RIfM@\expandafter\roman@\else
 \expandafter\nonmatherr@\expandafter\roman\fi}
\def\roman@#1{{\roman@@{#1}}}
\def\roman@@#1{\fam\rmfam\relax#1}
\def\italic{\RIfM@\expandafter\italic@\else
 \expandafter\nonmatherr@\expandafter\italic\fi}
\def\italic@#1{{\italic@@{#1}}}
\def\italic@@#1{\fam\itfam\relax#1}
\def\Cal{\RIfM@\expandafter\Cal@\else
 \expandafter\nonmatherr@\expandafter\Cal\fi}
\def\Cal@#1{{\Cal@@{#1}}}
\def\Cal@@#1{\noaccents@\fam\tw@#1}
\mathchardef\Gamma="0000
\mathchardef\Delta="0001
\mathchardef\Theta="0002
\mathchardef\Lambda="0003
\mathchardef\Xi="0004
\mathchardef\Pi="0005
\mathchardef\Sigma="0006
\mathchardef\Upsilon="0007
\mathchardef\Phi="0008
\mathchardef\Psi="0009
\mathchardef\Omega="000A
\mathchardef\varGamma="0100
\mathchardef\varDelta="0101
\mathchardef\varTheta="0102
\mathchardef\varLambda="0103
\mathchardef\varXi="0104
\mathchardef\varPi="0105
\mathchardef\varSigma="0106
\mathchardef\varUpsilon="0107
\mathchardef\varPhi="0108
\mathchardef\varPsi="0109
\mathchardef\varOmega="010A
\let\alloc@@\alloc@
\def\hexnumber@#1{\ifcase#1 0\or 1\or 2\or 3\or 4\or 5\or 6\or 7\or 8\or
 9\or A\or B\or C\or D\or E\or F\fi}
\def\loadmsam{%
 \font@\tenmsa=msam10
 \font@\sevenmsa=msam7
 \font@\fivemsa=msam5
 \alloc@@8\fam\chardef\sixt@@n\msafam
 \textfont\msafam=\tenmsa
 \scriptfont\msafam=\sevenmsa
 \scriptscriptfont\msafam=\fivemsa
 \edef\next{\hexnumber@\msafam}%
 \mathchardef\dabar@"0\next39
 \edef\dashrightarrow{\mathrel{\dabar@\dabar@\mathchar"0\next4B}}%
 \edef\dashleftarrow{\mathrel{\mathchar"0\next4C\dabar@\dabar@}}%
 \let\dasharrow\dashrightarrow
 \edef\ulcorner{\delimiter"4\next70\next70 }%
 \edef\urcorner{\delimiter"5\next71\next71 }%
 \edef\llcorner{\delimiter"4\next78\next78 }%
 \edef\lrcorner{\delimiter"5\next79\next79 }%
 \edef\yen{{\noexpand\mathhexbox@\next55}}%
 \edef\checkmark{{\noexpand\mathhexbox@\next58}}%
 \edef\circledR{{\noexpand\mathhexbox@\next72}}%
 \edef\maltese{{\noexpand\mathhexbox@\next7A}}%
 \global\let\loadmsam\empty}%
\def\loadmsbm{%
 \font@\tenmsb=msbm10 \font@\sevenmsb=msbm7 \font@\fivemsb=msbm5
 \alloc@@8\fam\chardef\sixt@@n\msbfam
 \textfont\msbfam=\tenmsb
 \scriptfont\msbfam=\sevenmsb \scriptscriptfont\msbfam=\fivemsb
 \global\let\loadmsbm\empty
 }
\def\widehat#1{\ifx\undefined\msbfam \DN@{362}%
  \else \setboxz@h{$\m@th#1$}%
    \edef\next@{\ifdim\wdz@>\tw@ em%
        \hexnumber@\msbfam 5B%
      \else 362\fi}\fi
  \mathaccent"0\next@{#1}}
\def\widetilde#1{\ifx\undefined\msbfam \DN@{365}%
  \else \setboxz@h{$\m@th#1$}%
    \edef\next@{\ifdim\wdz@>\tw@ em%
        \hexnumber@\msbfam 5D%
      \else 365\fi}\fi
  \mathaccent"0\next@{#1}}
\message{\string\newsymbol,}
\def\newsymbol#1#2#3#4#5{\define#1{}%
  \count@#2\relax \advance\count@\m@ne 
 \ifcase\count@
   \ifx\undefined\msafam\loadmsam\fi \let\next@\msafam
 \or \ifx\undefined\msbfam\loadmsbm\fi \let\next@\msbfam
 \else  \Err@{\Invalid@@\string\newsymbol}\let\next@\tw@\fi
 \mathchardef#1="#3\hexnumber@\next@#4#5\space}

\def\Bbb{\RIfM@\expandafter\Bbb@\else
 \expandafter\nonmatherr@\expandafter\Bbb\fi}
\def\Bbb@#1{{\Bbb@@{#1}}}
\def\Bbb@@#1{\noaccents@\fam\msbfam\relax#1}
\message{bold Greek and bold symbols,}
\def\loadbold{%
 \font@\tencmmib=cmmib10 \font@\sevencmmib=cmmib7 \font@\fivecmmib=cmmib5
 \skewchar\tencmmib'177 \skewchar\sevencmmib'177 \skewchar\fivecmmib'177
 \alloc@@8\fam\chardef\sixt@@n\cmmibfam
 \textfont\cmmibfam\tencmmib
 \scriptfont\cmmibfam\sevencmmib \scriptscriptfont\cmmibfam\fivecmmib
 \font@\tencmbsy=cmbsy10 \font@\sevencmbsy=cmbsy7 \font@\fivecmbsy=cmbsy5
 \skewchar\tencmbsy'60 \skewchar\sevencmbsy'60 \skewchar\fivecmbsy'60
 \alloc@@8\fam\chardef\sixt@@n\cmbsyfam
 \textfont\cmbsyfam\tencmbsy
 \scriptfont\cmbsyfam\sevencmbsy \scriptscriptfont\cmbsyfam\fivecmbsy
 \let\loadbold\empty
}
\def\boldnotloaded#1{\Err@{\ifcase#1\or First\else Second\fi
       bold symbol font not loaded}}
\def\mathchari@#1#2#3{\ifx\undefined\cmmibfam
    \boldnotloaded@\@ne
  \else\mathchar"#1\hexnumber@\cmmibfam#2#3\space \fi}
\def\mathcharii@#1#2#3{\ifx\undefined\cmbsyfam
    \boldnotloaded\tw@
  \else \mathchar"#1\hexnumber@\cmbsyfam#2#3\space\fi}
\edef\bffam@{\hexnumber@\bffam}
\def\boldkey#1{\ifcat\noexpand#1A%
  \ifx\undefined\cmmibfam \boldnotloaded\@ne
  \else {\fam\cmmibfam#1}\fi
 \else
 \ifx#1!\mathchar"5\bffam@21 \else
 \ifx#1(\mathchar"4\bffam@28 \else\ifx#1)\mathchar"5\bffam@29 \else
 \ifx#1+\mathchar"2\bffam@2B \else\ifx#1:\mathchar"3\bffam@3A \else
 \ifx#1;\mathchar"6\bffam@3B \else\ifx#1=\mathchar"3\bffam@3D \else
 \ifx#1?\mathchar"5\bffam@3F \else\ifx#1[\mathchar"4\bffam@5B \else
 \ifx#1]\mathchar"5\bffam@5D \else
 \ifx#1,\mathchari@63B \else
 \ifx#1-\mathcharii@200 \else
 \ifx#1.\mathchari@03A \else
 \ifx#1/\mathchari@03D \else
 \ifx#1<\mathchari@33C \else
 \ifx#1>\mathchari@33E \else
 \ifx#1*\mathcharii@203 \else
 \ifx#1|\mathcharii@06A \else
 \ifx#10\bold0\else\ifx#11\bold1\else\ifx#12\bold2\else\ifx#13\bold3\else
 \ifx#14\bold4\else\ifx#15\bold5\else\ifx#16\bold6\else\ifx#17\bold7\else
 \ifx#18\bold8\else\ifx#19\bold9\else
  \Err@{\string\boldkey\space can't be used with #1}%
 \fi\fi\fi\fi\fi\fi\fi\fi\fi\fi\fi\fi\fi\fi\fi
 \fi\fi\fi\fi\fi\fi\fi\fi\fi\fi\fi\fi\fi\fi}
\def\boldsymbol#1{%
 \DN@{\Err@{You can't use \string\boldsymbol\space with \string#1}#1}%
 \ifcat\noexpand#1A%
   \let\next@\relax
   \ifx\undefined\cmmibfam \boldnotloaded\@ne
   \else {\fam\cmmibfam#1}\fi
 \else
  \xdef\meaning@{\meaning#1.........}%
  \expandafter\math@\meaning@\math@
  \ifmath@
   \expandafter\mathch@\meaning@\mathch@
   \ifmathch@
    \expandafter\boldsymbol@@\meaning@\boldsymbol@@
   \fi
  \else
   \expandafter\macro@\meaning@\macro@
   \expandafter\delim@\meaning@\delim@
   \ifdelim@
    \expandafter\delim@@\meaning@\delim@@
   \else
    \boldsymbol@{#1}%
   \fi
  \fi
 \fi
 \next@}
\def\mathhexboxii@#1#2{\ifx\undefined\cmbsyfam
    \boldnotloaded\tw@
  \else \mathhexbox@{\hexnumber@\cmbsyfam}{#1}{#2}\fi}
\def\boldsymbol@#1{\let\next@\relax\let\next#1%
 \ifx\next\cdot\mathcharii@201 \else
 \ifx\next\prime{{\null\mathcharii@030 \null}}\else
 \ifx\next\lbrack\mathchar"4\bffam@5B \else
 \ifx\next\rbrack\mathchar"5\bffam@5D \else
 \ifx\next\{\mathcharii@466 \else
 \ifx\next\lbrace\mathcharii@466 \else
 \ifx\next\}\mathcharii@567 \else
 \ifx\next\rbrace\mathcharii@567 \else
 \ifx\next\surd{{\mathcharii@170}}\else
 \ifx\next\S{{\mathhexboxii@78}}\else
 \ifx\next\P{{\mathhexboxii@7B}}\else
 \ifx\next\dag{{\mathhexboxii@79}}\else
 \ifx\next\ddag{{\mathhexboxii@7A}}\else
 \DN@{\Err@{You can't use \string\boldsymbol\space with \string#1}#1}%
 \fi\fi\fi\fi\fi\fi\fi\fi\fi\fi\fi\fi\fi}
\def\boldsymbol@@#1.#2\boldsymbol@@{\classnum@#1 \count@@@\classnum@        
 \divide\classnum@4096 \count@\classnum@                                    
 \multiply\count@4096 \advance\count@@@-\count@ \count@@\count@@@           
 \divide\count@@@\@cclvi \count@\count@@                                    
 \multiply\count@@@\@cclvi \advance\count@@-\count@@@                       
 \divide\count@@@\@cclvi                                                    
 \multiply\classnum@4096 \advance\classnum@\count@@                         
 \ifnum\count@@@=\z@                                                        
  \count@"\bffam@ \multiply\count@\@cclvi
  \advance\classnum@\count@
  \DN@{\mathchar\number\classnum@}%
 \else
  \ifnum\count@@@=\@ne                                                      
   \ifx\undefined\cmmibfam \DN@{\boldnotloaded\@ne}%
   \else \count@\cmmibfam \multiply\count@\@cclvi
     \advance\classnum@\count@
     \DN@{\mathchar\number\classnum@}\fi
  \else
   \ifnum\count@@@=\tw@                                                    
     \ifx\undefined\cmbsyfam
       \DN@{\boldnotloaded\tw@}%
     \else
       \count@\cmbsyfam \multiply\count@\@cclvi
       \advance\classnum@\count@
       \DN@{\mathchar\number\classnum@}%
     \fi
  \fi
 \fi
\fi}
\newif\ifdelim@
\newcount\delimcount@
{\uccode`6=`\\ \uccode`7=`d \uccode`8=`e \uccode`9=`l
 \uppercase{\gdef\delim@#1#2#3#4#5\delim@
  {\delim@false\ifx 6#1\ifx 7#2\ifx 8#3\ifx 9#4\delim@true
   \xdef\meaning@{#5}\fi\fi\fi\fi}}}
\def\delim@@#1"#2#3#4#5#6\delim@@{\if#32%
\let\next@\relax
 \ifx\undefined\cmbsyfam \boldnotloaded\@ne
 \else \mathcharii@#2#4#5\space \fi\fi}
\def\vert{\delimiter"026A30C }
\def\Vert{\delimiter"026B30D }
\let\|\Vert
\def\backslash{\delimiter"026E30F }
\def\boldkeydots@#1{\bold@true\let\next=#1\let\delayed@=#1\mdots@@
 \boldkey#1\bold@false}  
\def\boldsymboldots@#1{\bold@true\let\next#1\let\delayed@#1\mdots@@
 \boldsymbol#1\bold@false}
\message{Euler fonts,}

\def\frak{\mathfont@\frak}

\def\loadmathfont#1{%
   \expandafter\font@\csname ten#1\endcsname=#110
   \expandafter\font@\csname seven#1\endcsname=#17
   \expandafter\font@\csname five#1\endcsname=#15
   \edef\next{\noexpand\alloc@@8\fam\chardef\sixt@@n
     \expandafter\noexpand\csname#1fam\endcsname}%
   \next
   \textfont\csname#1fam\endcsname \csname ten#1\endcsname
   \scriptfont\csname#1fam\endcsname \csname seven#1\endcsname
   \scriptscriptfont\csname#1fam\endcsname \csname five#1\endcsname
   \expandafter\def\csname #1\expandafter\endcsname\expandafter{%
      \expandafter\mathfont@\csname#1\endcsname}%
 \expandafter\gdef\csname load#1\endcsname{}%
}
\def\mathfont@#1{\RIfM@\expandafter\mathfont@@\expandafter#1\else
  \expandafter\nonmatherr@\expandafter#1\fi}
\def\mathfont@@#1#2{{\mathfont@@@#1{#2}}}
\def\mathfont@@@#1#2{\noaccents@
   \fam\csname\expandafter\eat@\string#1fam\endcsname
   \relax#2}
\message{math accents,}
\def\accentclass@{7}
\def\noaccents@{\def\accentclass@{0}}
\def\makeacc@#1#2{\def#1{\mathaccent"\accentclass@#2 }}
\makeacc@\hat{05E}
\makeacc@\check{014}
\makeacc@\tilde{07E}
\makeacc@\acute{013}
\makeacc@\grave{012}
\makeacc@\dot{05F}
\makeacc@\ddot{07F}
\makeacc@\breve{015}
\makeacc@\bar{016}

\newcount\skewcharcount@
\newcount\familycount@
\def\theskewchar@{\familycount@\@ne
 \global\skewcharcount@\the\skewchar\textfont\@ne                           
 \ifnum\fam>\m@ne\ifnum\fam<16
  \global\familycount@\the\fam\relax
  \global\skewcharcount@\the\skewchar\textfont\the\fam\relax\fi\fi          
 \ifnum\skewcharcount@>\m@ne
  \ifnum\skewcharcount@<128
  \multiply\familycount@256
  \global\advance\skewcharcount@\familycount@
  \global\advance\skewcharcount@28672
  \mathchar\skewcharcount@\else
  \global\skewcharcount@\m@ne\fi\else
 \global\skewcharcount@\m@ne\fi}                                            
\newcount\pointcount@
\def\getpoints@#1.#2\getpoints@{\pointcount@#1 }
\newdimen\accentdimen@
\newcount\accentmu@
\def\dimentomu@{\multiply\accentdimen@ 100
 \expandafter\getpoints@\the\accentdimen@\getpoints@
 \multiply\pointcount@18
 \divide\pointcount@\@m
 \global\accentmu@\pointcount@}
\def\Makeacc@#1#2{\def#1{\RIfM@\DN@{\mathaccent@
 {"\accentclass@#2 }}\else\DN@{\nonmatherr@{#1}}\fi\next@}}
\def\unbracefonts@{\let\Cal@\Cal@@\let\roman@\roman@@\let\bold@\bold@@
 \let\slanted@\slanted@@}
\def\mathaccent@#1#2{\ifnum\fam=\m@ne\xdef\thefam@{1}\else
 \xdef\thefam@{\the\fam}\fi                                                 
 \accentdimen@\z@                                                           
 \setboxz@h{\unbracefonts@$\m@th\fam\thefam@\relax#2$}
 \ifdim\accentdimen@=\z@\DN@{\mathaccent#1{#2}}
  \setbox@ne\hbox{\unbracefonts@$\m@th\fam\thefam@\relax#2\theskewchar@$}
  \setbox\tw@\hbox{$\m@th\ifnum\skewcharcount@=\m@ne\else
   \mathchar\skewcharcount@\fi$}
  \global\accentdimen@\wd@ne\global\advance\accentdimen@-\wdz@
  \global\advance\accentdimen@-\wd\tw@                                     
  \global\multiply\accentdimen@\tw@
  \dimentomu@\global\advance\accentmu@\@ne                                 
 \else\DN@{{\mathaccent#1{#2\mkern\accentmu@ mu}%
    \mkern-\accentmu@ mu}{}}\fi                                             
 \next@}\Makeacc@\Hat{05E}
\Makeacc@\Check{014}
\Makeacc@\Tilde{07E}
\Makeacc@\Acute{013}
\Makeacc@\Grave{012}
\Makeacc@\Dot{05F}
\Makeacc@\Ddot{07F}
\Makeacc@\Breve{015}
\Makeacc@\Bar{016}
\def\Vec{\RIfM@\DN@{\mathaccent@{"017E }}\else
 \DN@{\nonmatherr@\Vec}\fi\next@}
\def\accentedsymbol#1#2{\csname newbox\expandafter\endcsname
  \csname\expandafter\eat@\string#1@box\endcsname
 \expandafter\setbox\csname\expandafter\eat@
  \string#1@box\endcsname\hbox{$\m@th#2$}\define
  #1{\copy\csname\expandafter\eat@\string#1@box\endcsname{}}}
\message{roots,}
\def\sqrt#1{\radical"270370 {#1}}
\let\underline@\underline
\let\overline@\overline
\def\underline#1{\underline@{#1}}
\def\overline#1{\overline@{#1}}
\Invalid@\leftroot
\Invalid@\uproot
\newcount\uproot@
\newcount\leftroot@
\def\root{\relaxnext@
  \DN@{\ifx\next\uproot\let\next@\nextii@\else
   \ifx\next\leftroot\let\next@\nextiii@\else
   \let\next@\plainroot@\fi\fi\next@}%
  \DNii@\uproot##1{\uproot@##1\relax\FN@\nextiv@}%
  \def\nextiv@{\ifx\next\space@\DN@. {\FN@\nextv@}\else
   \DN@.{\FN@\nextv@}\fi\next@.}%
  \def\nextv@{\ifx\next\leftroot\let\next@\nextvi@\else
   \let\next@\plainroot@\fi\next@}%
  \def\nextvi@\leftroot##1{\leftroot@##1\relax\plainroot@}%
   \def\nextiii@\leftroot##1{\leftroot@##1\relax\FN@\nextvii@}%
  \def\nextvii@{\ifx\next\space@
   \DN@. {\FN@\nextviii@}\else
   \DN@.{\FN@\nextviii@}\fi\next@.}%
  \def\nextviii@{\ifx\next\uproot\let\next@\nextix@\else
   \let\next@\plainroot@\fi\next@}%
  \def\nextix@\uproot##1{\uproot@##1\relax\plainroot@}%
  \bgroup\uproot@\z@\leftroot@\z@\FN@\next@}
\def\plainroot@#1\of#2{\setbox\rootbox\hbox{$\m@th\scriptscriptstyle{#1}$}%
 \mathchoice{\r@@t\displaystyle{#2}}{\r@@t\textstyle{#2}}
 {\r@@t\scriptstyle{#2}}{\r@@t\scriptscriptstyle{#2}}\egroup}
\def\r@@t#1#2{\setboxz@h{$\m@th#1\sqrt{#2}$}%
 \dimen@\ht\z@\advance\dimen@-\dp\z@
 \setbox@ne\hbox{$\m@th#1\mskip\uproot@ mu$}\advance\dimen@ 1.667\wd@ne
 \mkern-\leftroot@ mu\mkern5mu\raise.6\dimen@\copy\rootbox
 \mkern-10mu\mkern\leftroot@ mu\boxz@}
\def\boxed#1{\setboxz@h{$\m@th\displaystyle{#1}$}\dimen@.4\ex@
 \advance\dimen@3\ex@\advance\dimen@\dp\z@
 \hbox{\lower\dimen@\hbox{%
 \vbox{\hrule height.4\ex@
 \hbox{\vrule width.4\ex@\hskip3\ex@\vbox{\vskip3\ex@\boxz@\vskip3\ex@}%
 \hskip3\ex@\vrule width.4\ex@}\hrule height.4\ex@}%
 }}}
\message{commutative diagrams,}
\let\ampersand@\relax
\newdimen\minaw@
\minaw@11.11128\ex@
\newdimen\minCDaw@
\minCDaw@2.5pc
\def\minCDarrowwidth#1{\RIfMIfI@\onlydmatherr@\minCDarrowwidth
 \else\minCDaw@#1\relax\fi\else\onlydmatherr@\minCDarrowwidth\fi}
\newif\ifCD@
\def\CD{\bgroup\vspace@\relax\let\ampersand@&\iffalse}\fi
 \CD@true\vcenter\bgroup\Let@\tabskip\z@skip\baselineskip20\ex@
 \lineskip3\ex@\lineskiplimit3\ex@\halign\bgroup
 &\hfill$\m@th##$\hfill\crcr}
\def\endCD{\crcr\egroup\egroup\egroup}
\newdimen\bigaw@
\atdef@>#1>#2>{\ampersand@                                                  
 \setboxz@h{$\m@th\ssize\;{#1}\;\;$}
 \setbox@ne\hbox{$\m@th\ssize\;{#2}\;\;$}
 \setbox\tw@\hbox{$\m@th#2$}
 \ifCD@\global\bigaw@\minCDaw@\else\global\bigaw@\minaw@\fi                 
 \ifdim\wdz@>\bigaw@\global\bigaw@\wdz@\fi
 \ifdim\wd@ne>\bigaw@\global\bigaw@\wd@ne\fi                                
 \ifCD@\enskip\fi                                                           
 \ifdim\wd\tw@>\z@
  \mathrel{\mathop{\hbox to\bigaw@{\rightarrowfill@\displaystyle}}%
    \limits^{#1}_{#2}}
 \else\mathrel{\mathop{\hbox to\bigaw@{\rightarrowfill@\displaystyle}}%
    \limits^{#1}}\fi                                                        
 \ifCD@\enskip\fi                                                          
 \ampersand@}                                                              
\atdef@<#1<#2<{\ampersand@\setboxz@h{$\m@th\ssize\;\;{#1}\;$}%
 \setbox@ne\hbox{$\m@th\ssize\;\;{#2}\;$}\setbox\tw@\hbox{$\m@th#2$}%
 \ifCD@\global\bigaw@\minCDaw@\else\global\bigaw@\minaw@\fi
 \ifdim\wdz@>\bigaw@\global\bigaw@\wdz@\fi
 \ifdim\wd@ne>\bigaw@\global\bigaw@\wd@ne\fi
 \ifCD@\enskip\fi
 \ifdim\wd\tw@>\z@
  \mathrel{\mathop{\hbox to\bigaw@{\leftarrowfill@\displaystyle}}%
       \limits^{#1}_{#2}}\else
  \mathrel{\mathop{\hbox to\bigaw@{\leftarrowfill@\displaystyle}}%
       \limits^{#1}}\fi
 \ifCD@\enskip\fi\ampersand@}
\begingroup
 \catcode`\~=\active \lccode`\~=`\@
 \lowercase{%
  \global\atdef@)#1)#2){~>#1>#2>}
  \global\atdef@(#1(#2({~<#1<#2<}}
\endgroup
\atdef@ A#1A#2A{\llap{$\m@th\vcenter{\hbox
 {$\ssize#1$}}$}\Big\uparrow\rlap{$\m@th\vcenter{\hbox{$\ssize#2$}}$}&&}
\atdef@ V#1V#2V{\llap{$\m@th\vcenter{\hbox
 {$\ssize#1$}}$}\Big\downarrow\rlap{$\m@th\vcenter{\hbox{$\ssize#2$}}$}&&}
\atdef@={&\enskip\mathrel
 {\vbox{\hrule width\minCDaw@\vskip3\ex@\hrule width
 \minCDaw@}}\enskip&}
\atdef@|{\Big\Vert&&}
\atdef@\vert{\Big\Vert&&}
\def\pretend#1\haswidth#2{\setboxz@h{$\m@th\scriptstyle{#2}$}\hbox
 to\wdz@{\hfill$\m@th\scriptstyle{#1}$\hfill}}
\message{poor man's bold,}
\def\pmb{\RIfM@\expandafter\mathpalette\expandafter\pmb@\else
 \expandafter\pmb@@\fi}
\def\pmb@@#1{\leavevmode\setboxz@h{#1}%
   \dimen@-\wdz@
   \kern-.5\ex@\copy\z@
   \kern\dimen@\kern.25\ex@\raise.4\ex@\copy\z@
   \kern\dimen@\kern.25\ex@\box\z@
}
\def\binrel@@#1{\ifdim\wd2<\z@\mathbin{#1}\else\ifdim\wd\tw@>\z@
 \mathrel{#1}\else{#1}\fi\fi}
\newdimen\pmbraise@
\def\pmb@#1#2{\setbox\thr@@\hbox{$\m@th#1{#2}$}%
 \setbox4\hbox{$\m@th#1\mkern.5mu$}\pmbraise@\wd4\relax
 \binrel@{#2}%
 \dimen@-\wd\thr@@
   \binrel@@{%
   \mkern-.8mu\copy\thr@@
   \kern\dimen@\mkern.4mu\raise\pmbraise@\copy\thr@@
   \kern\dimen@\mkern.4mu\box\thr@@
}}
\def\documentstyle#1{\W@{}\input #1.sty\relax}
\message{syntax check,}
\font\dummyft@=dummy
\fontdimen1 \dummyft@=\z@
\fontdimen2 \dummyft@=\z@
\fontdimen3 \dummyft@=\z@
\fontdimen4 \dummyft@=\z@
\fontdimen5 \dummyft@=\z@
\fontdimen6 \dummyft@=\z@
\fontdimen7 \dummyft@=\z@
\fontdimen8 \dummyft@=\z@
\fontdimen9 \dummyft@=\z@
\fontdimen10 \dummyft@=\z@
\fontdimen11 \dummyft@=\z@
\fontdimen12 \dummyft@=\z@
\fontdimen13 \dummyft@=\z@
\fontdimen14 \dummyft@=\z@
\fontdimen15 \dummyft@=\z@
\fontdimen16 \dummyft@=\z@
\fontdimen17 \dummyft@=\z@
\fontdimen18 \dummyft@=\z@
\fontdimen19 \dummyft@=\z@
\fontdimen20 \dummyft@=\z@
\fontdimen21 \dummyft@=\z@
\fontdimen22 \dummyft@=\z@
\def\fontlist@{\\{\tenrm}\\{\sevenrm}\\{\fiverm}\\{\teni}\\{\seveni}%
 \\{\fivei}\\{\tensy}\\{\sevensy}\\{\fivesy}\\{\tenex}\\{\tenbf}\\{\sevenbf}%
 \\{\fivebf}\\{\tensl}\\{\tenit}}
\def\font@#1=#2 {\rightappend@#1\to\fontlist@\font#1=#2 }
\def\dodummy@{{\def\\##1{\global\let##1\dummyft@}\fontlist@}}
\def\nopages@{\output{\setbox\z@\box\@cclv \deadcycles\z@}%
 \alloc@5\toks\toksdef\@cclvi\output}
\let\galleys\nopages@
\newif\ifsyntax@
\newcount\countxviii@
\def\syntax{\syntax@true\dodummy@\countxviii@\count18
 \loop\ifnum\countxviii@>\m@ne\textfont\countxviii@=\dummyft@
 \scriptfont\countxviii@=\dummyft@\scriptscriptfont\countxviii@=\dummyft@
 \advance\countxviii@\m@ne\repeat                                           
 \dummyft@\tracinglostchars\z@\nopages@\frenchspacing\hbadness\@M}
\def\first@#1#2\end{#1}
\def\printoptions{\W@{Do you want S(yntax check),
  G(alleys) or P(ages)?}%
 \message{Type S, G or P, followed by <return>: }%
 \begingroup 
 \endlinechar\m@ne 
 \read\m@ne to\ans@
 \edef\ans@{\uppercase{\def\noexpand\ans@{%
   \expandafter\first@\ans@ P\end}}}%
 \expandafter\endgroup\ans@
 \if\ans@ P
 \else \if\ans@ S\syntax
 \else \if\ans@ G\galleys
 \else\message{? Unknown option: \ans@; using the `pages' option.}%
 \fi\fi\fi}
\def\alloc@#1#2#3#4#5{\global\advance\count1#1by\@ne
 \ch@ck#1#4#2\allocationnumber=\count1#1
 \global#3#5=\allocationnumber
 \ifalloc@\wlog{\string#5=\string#2\the\allocationnumber}\fi}
\def\document{\def\alloclist@{}\def\fontlist@{}}
\let\enddocument\bye

\let\proclaim\undefined
\let\footnote\undefined
\let\=\undefined
\let\>\undefined

\catcode`\@=\active
\message{... finished}

\expandafter\ifx\csname mathdefs.tex\endcsname\relax
  \expandafter\gdef\csname mathdefs.tex\endcsname{}
\else \message{Hey!  Apparently you were trying to
  \string\input{mathdefs.tex} twice.   This does not make sense.} 
\errmessage{Please edit your file (probably \jobname.tex) and remove
any duplicate ``\string\input'' lines}\endinput\fi




\catcode`\X=12\catcode`\@=11

\def\n@wcount{\alloc@0\count\countdef\insc@unt}
\def\n@wwrite{\alloc@7\write\chardef\sixt@@n}
\def\n@wread{\alloc@6\read\chardef\sixt@@n}
\def\r@s@t{\relax}\def\v@idline{\par}\def\@mputate#1/{#1}
\def\l@c@l#1X{\firstpart.#1}\def\gl@b@l#1X{#1}\def\t@d@l#1X{{}}

\def\crossrefs#1{\ifx\all#1\let\tr@ce=\all\else\def\tr@ce{#1,}\fi
   \n@wwrite\cit@tionsout\openout\cit@tionsout=\jobname.cit 
   \write\cit@tionsout{\tr@ce}\expandafter\setfl@gs\tr@ce,}
\def\setfl@gs#1,{\def\@{#1}\ifx\@\empty\let\next=\relax
   \else\let\next=\setfl@gs\expandafter\xdef
   \csname#1tr@cetrue\endcsname{}\fi\next}
\def\m@ketag#1#2{\expandafter\n@wcount\csname#2tagno\endcsname
     \csname#2tagno\endcsname=0\let\tail=\all\xdef\all{\tail#2,}
   \ifx#1\l@c@l\let\tail=\r@s@t\xdef\r@s@t{\csname#2tagno\endcsname=0\tail}\fi
   \expandafter\gdef\csname#2cite\endcsname##1{\expandafter
     \ifx\csname#2tag##1\endcsname\relax?\else\csname#2tag##1\endcsname\fi
     \expandafter\ifx\csname#2tr@cetrue\endcsname\relax\else
     \write\cit@tionsout{#2tag ##1 cited on page \folio.}\fi}
   \expandafter\gdef\csname#2page\endcsname##1{\expandafter
     \ifx\csname#2page##1\endcsname\relax?\else\csname#2page##1\endcsname\fi
     \expandafter\ifx\csname#2tr@cetrue\endcsname\relax\else
     \write\cit@tionsout{#2tag ##1 cited on page \folio.}\fi}
   \expandafter\gdef\csname#2tag\endcsname##1{\expandafter
      \ifx\csname#2check##1\endcsname\relax
      \expandafter\xdef\csname#2check##1\endcsname{}%
      \else\immediate\write16{Warning: #2tag ##1 used more than once.}\fi
      \multit@g{#1}{#2}##1/X%
      \write\t@gsout{#2tag ##1 assigned number \csname#2tag##1\endcsname\space
      on page \number\count0.}%
   \csname#2tag##1\endcsname}}

\def\multit@g#1#2#3/#4X{\def\t@mp{#4}\ifx\t@mp\empty%
      \global\advance\csname#2tagno\endcsname by 1 
      \expandafter\xdef\csname#2tag#3\endcsname
      {#1\number\csname#2tagno\endcsnameX}%
   \else\expandafter\ifx\csname#2last#3\endcsname\relax
      \expandafter\n@wcount\csname#2last#3\endcsname
      \global\advance\csname#2tagno\endcsname by 1 
      \expandafter\xdef\csname#2tag#3\endcsname
      {#1\number\csname#2tagno\endcsnameX}
      \write\t@gsout{#2tag #3 assigned number \csname#2tag#3\endcsname\space
      on page \number\count0.}\fi
   \global\advance\csname#2last#3\endcsname by 1
   \def\t@mp{\expandafter\xdef\csname#2tag#3/}%
   \expandafter\t@mp\@mputate#4\endcsname
   {\csname#2tag#3\endcsname\lastpart{\csname#2last#3\endcsname}}\fi}
\def\t@gs#1{\def\all{}\m@ketag#1e\m@ketag#1s\m@ketag\t@d@l p
\let\realscite\scite
\let\realstag\stag
   \m@ketag\gl@b@l r \n@wread\t@gsin
   \openin\t@gsin=\jobname.tgs \re@der \closein\t@gsin
   \n@wwrite\t@gsout\openout\t@gsout=\jobname.tgs }
\outer\def\localtags{\t@gs\l@c@l}
\outer\def\globaltags{\t@gs\gl@b@l}
\outer\def\newlocaltag#1{\m@ketag\l@c@l{#1}}
\outer\def\newglobaltag#1{\m@ketag\gl@b@l{#1}}

\newif\ifpr@ 
\def\m@kecs #1tag #2 assigned number #3 on page #4.%
   {\expandafter\gdef\csname#1tag#2\endcsname{#3}
   \expandafter\gdef\csname#1page#2\endcsname{#4}
   \ifpr@\expandafter\xdef\csname#1check#2\endcsname{}\fi}
\def\re@der{\ifeof\t@gsin\let\next=\relax\else
   \read\t@gsin to\t@gline\ifx\t@gline\v@idline\else
   \expandafter\m@kecs \t@gline\fi\let \next=\re@der\fi\next}
\def\pretags#1{\pr@true\pret@gs#1,,}
\def\pret@gs#1,{\def\@{#1}\ifx\@\empty\let\n@xtfile=\relax
   \else\let\n@xtfile=\pret@gs \openin\t@gsin=#1.tgs \message{#1} \re@der 
   \closein\t@gsin\fi \n@xtfile}

\newcount\sectno\sectno=0\newcount\subsectno\subsectno=0
\newif\ifultr@local \def\ultralocal{\ultr@localtrue}
\def\firstpart{\number\sectno}
\def\lastpart#1{\ifcase#1 \or a\or b\or c\or d\or e\or f\or g\or h\or 
   i\or k\or l\or m\or n\or o\or p\or q\or r\or s\or t\or u\or v\or w\or 
   x\or y\or z \fi}

\def\resetall{\global\advance\sectno by 1\subsectno=0
   \gdef\firstpart{\number\sectno}\r@s@t}
\def\resetsub{\global\advance\subsectno by 1
   \gdef\firstpart{\number\sectno.\number\subsectno}\r@s@t}
\def\newsection#1\par{\resetall\vskip0pt plus.3\vsize\penalty-250
   \vskip0pt plus-.3\vsize\bigskip\bigskip
   \message{#1}\leftline{\bf#1}\nobreak\bigskip}
\def\subsection#1\par{\ifultr@local\resetsub\fi
   \vskip0pt plus.2\vsize\penalty-250\vskip0pt plus-.2\vsize
   \bigskip\smallskip\message{#1}\leftline{\bf#1}\nobreak\medskip}


\newdimen\marginshift

\newdimen\margindelta
\newdimen\marginmax
\newdimen\marginmin

\def\margininit{       
\marginmax=3 true cm                  
				      
\margindelta=0.1 true cm              
\marginmin=0.1true cm                 
\marginshift=\marginmin
}    

\def\t@gsjj#1,{\def\@{#1}\ifx\@\empty\let\next=\relax\else\let\next=\t@gsjj
   \def\@@{p}\ifx\@\@@\else
   \expandafter\gdef\csname#1cite\endcsname##1{\citejj{##1}}
   \expandafter\gdef\csname#1page\endcsname##1{?}
   \expandafter\gdef\csname#1tag\endcsname##1{\tagjj{##1}}\fi\fi\next}
\newif\ifshowstuffinmargin
\showstuffinmarginfalse
\def\jjtags{\ifx\shlhetal\relax 
  \else
\ifx\shlhetal\undefinedcontrolseq
\else
\showstuffinmargintrue
\ifx\all\relax\else\expandafter\t@gsjj\all,\fi\fi \fi
}

\def\tagjj#1{\realstag{#1}\mginpar{\zeigen{#1}}}
\def\citejj#1{\rechnen{#1}\mginpar{\zeigen{#1}}}     

\def\rechnen#1{\expandafter\ifx\csname stag#1\endcsname\relax ??\else
                           \csname stag#1\endcsname\fi}

\newdimen\theight

\def\marginfont{\sevenrm}

\def\trymarginbox#1{\setbox0=\hbox{\marginfont\hskip\marginshift #1}%
		\global\marginshift\wd0 
		\global\advance\marginshift\margindelta}

\def \mginpar#1{%
\ifvmode\setbox0\hbox to \hsize{\hfill\rlap{\marginfont\quad#1}}%
\ht0 0cm
\dp0 0cm
\box0\vskip-\baselineskip
\else 
             \vadjust{\trymarginbox{#1}%
		\ifdim\marginshift>\marginmax \global\marginshift\marginmin
			\trymarginbox{#1}%
                \fi
             \theight=\ht0
             \advance\theight by \dp0    \advance\theight by \lineskip
             \kern -\theight \vbox to \theight{\rightline{\rlap{\box0}}%
\vss}}\fi}


\def\t@gsoff#1,{\def\@{#1}\ifx\@\empty\let\next=\relax\else\let\next=\t@gsoff
   \def\@@{p}\ifx\@\@@\else
   \expandafter\gdef\csname#1cite\endcsname##1{\zeigen{##1}}
   \expandafter\gdef\csname#1page\endcsname##1{?}
   \expandafter\gdef\csname#1tag\endcsname##1{\zeigen{##1}}\fi\fi\next}
\def\verbatimtags{\showstuffinmarginfalse
\ifx\all\relax\else\expandafter\t@gsoff\all,\fi}
\def\zeigen#1{\hbox{$\langle$}#1\hbox{$\rangle$}}

\def\margincite#1{\ifshowstuffinmargin\mginpar{\zeigen{#1}}\fi}

\def\margintag#1{\ifshowstuffinmargin\mginpar{\zeigen{#1}}\fi}

\def\(#1){\edef\dot@g{\ifmmode\ifinner(\hbox{\noexpand\etag{#1}})
   \else\noexpand\eqno(\hbox{\noexpand\etag{#1}})\fi
   \else(\noexpand\ecite{#1})\fi}\dot@g}

\newif\ifbr@ck
\def\eat#1{}
\def\[#1]{\br@cktrue[\br@cket#1'X]}
\def\br@cket#1'#2X{\def\temp{#2}\ifx\temp\empty\let\next\eat
   \else\let\next\br@cket\fi
   \ifbr@ck\br@ckfalse\br@ck@t#1,X\else\br@cktrue#1\fi\next#2X}
\def\br@ck@t#1,#2X{\def\temp{#2}\ifx\temp\empty\let\neext\eat
   \else\let\neext\br@ck@t\def\temp{,}\fi
   \def\teemp{#1}\ifx\teemp\empty\else\rcite{#1}\fi\temp\neext#2X}
\def\resetbr@cket{\gdef\[##1]{[\rtag{##1}]}}
\def\references{\resetbr@cket\newsection References\par}

\newtoks\symb@ls\newtoks\s@mb@ls\newtoks\p@gelist\n@wcount\ftn@mber
    \ftn@mber=1\newif\ifftn@mbers\ftn@mbersfalse\newif\ifbyp@ge\byp@gefalse
\def\defm@rk{\ifftn@mbers\n@mberm@rk\else\symb@lm@rk\fi}
\def\n@mberm@rk{\xdef\m@rk{{\the\ftn@mber}}%
    \global\advance\ftn@mber by 1 }
\def\rot@te#1{\let\temp=#1\global#1=\expandafter\r@t@te\the\temp,X}
\def\r@t@te#1,#2X{{#2#1}\xdef\m@rk{{#1}}}
\def\b@@st#1{{$^{#1}$}}\def\str@p#1{#1}
\def\symb@lm@rk{\ifbyp@ge\rot@te\p@gelist\ifnum\expandafter\str@p\m@rk=1 
    \s@mb@ls=\symb@ls\fi\write\f@nsout{\number\count0}\fi \rot@te\s@mb@ls}
\def\byp@ge{\byp@getrue\n@wwrite\f@nsin\openin\f@nsin=\jobname.fns 
    \n@wcount\currentp@ge\currentp@ge=0\p@gelist={0}
    \re@dfns\closein\f@nsin\rot@te\p@gelist
    \n@wread\f@nsout\openout\f@nsout=\jobname.fns }
\def\m@kelist#1X#2{{#1,#2}}
\def\re@dfns{\ifeof\f@nsin\let\next=\relax\else\read\f@nsin to \f@nline
    \ifx\f@nline\v@idline\else\let\t@mplist=\p@gelist
    \ifnum\currentp@ge=\f@nline
    \global\p@gelist=\expandafter\m@kelist\the\t@mplistX0
    \else\currentp@ge=\f@nline
    \global\p@gelist=\expandafter\m@kelist\the\t@mplistX1\fi\fi
    \let\next=\re@dfns\fi\next}
\def\symbols#1{\symb@ls={#1}\s@mb@ls=\symb@ls} 
\def\bigsymbol{\textstyle}
\symbols{\bigsymbol\ast,\dagger,\ddagger,\sharp,\flat,\natural,\star}
\def\ftnumbers{\ftn@mberstrue} \def\ftsymbols{\ftn@mbersfalse}
\def\paginal{\byp@ge} \def\resetftnumbers{\ftn@mber=1}
\def\ftnote#1{\defm@rk\expandafter\expandafter\expandafter\footnote
    \expandafter\b@@st\m@rk{#1}}

\long\def\jump#1\endjump{}
\def\ssum{\mathop{\lower .1em\hbox{$\textstyle\Sigma$}}\nolimits}

\def\qed{\nobreak\kern 1em \vrule height .5em width .5em depth 0em}
\def\newneq{\hbox{\rlap{\hbox to 1\wd9{\hss$=$\hss}}\raise .1em 
   \hbox to 1\wd9{\hss$\scriptscriptstyle/$\hss}}}
\def\subsetne{\setbox9 = \hbox{$\subset$}\mathrel{\hbox{\rlap
   {\lower .4em \newneq}\raise .13em \hbox{$\subset$}}}}
\def\supsetne{\setbox9 = \hbox{$\subset$}\mathrel{\hbox{\rlap
   {\lower .4em \newneq}\raise .13em \hbox{$\supset$}}}}

\def\vbar{\mathchoice{\vrule height6.3ptdepth-.5ptwidth.8pt\kern-.8pt}
   {\vrule height6.3ptdepth-.5ptwidth.8pt\kern-.8pt}
   {\vrule height4.1ptdepth-.35ptwidth.6pt\kern-.6pt}
   {\vrule height3.1ptdepth-.25ptwidth.5pt\kern-.5pt}}
\def\f@dge{\mathchoice{}{}{\mkern.5mu}{\mkern.8mu}}
\def\b@c#1#2{{\rm \mkern#2mu\vbar\mkern-#2mu#1}}
\def\b@b#1{{\rm I\mkern-3.5mu #1}}
\def\b@a#1#2{{\rm #1\mkern-#2mu\f@dge #1}}
\def\bb#1{{\count4=`#1 \advance\count4by-64 \ifcase\count4\or\b@a A{11.5}\or
   \b@b B\or\b@c C{5}\or\b@b D\or\b@b E\or\b@b F \or\b@c G{5}\or\b@b H\or
   \b@b I\or\b@c J{3}\or\b@b K\or\b@b L \or\b@b M\or\b@b N\or\b@c O{5} \or
   \b@b P\or\b@c Q{5}\or\b@b R\or\b@a S{8}\or\b@a T{10.5}\or\b@c U{5}\or
   \b@a V{12}\or\b@a W{16.5}\or\b@a X{11}\or\b@a Y{11.7}\or\b@a Z{7.5}\fi}}

\catcode`\X=11 \catcode`\@=12




\let\thischap\jobname

\def\partof#1{\csname returnthe#1part\endcsname}
\def\chapof#1{\csname returnthe#1chap\endcsname}

\def\setchapter#1,#2,#3.{%
  \expandafter\def\csname returnthe#1part\endcsname{#2}%
  \expandafter\def\csname returnthe#1chap\endcsname{#3}%
}

\setchapter 300a,A,I.
\setchapter 300b,A,II.
\setchapter 300c,A,III.
\setchapter 300d,A,IV.
\setchapter 300e,A,V.
\setchapter 300f,A,VI.
\setchapter 300g,A,VII.
\setchapter F604,B,0.
\setchapter  88r,B,I.
\setchapter  600,B,II.
\setchapter  705,B,III.

\def\cprefix#1{
\edef\theotherpart{\partof{#1}}\edef\theotherchap{\chapof{#1}}%
\ifx\theotherpart\thispart
   \ifx\theotherchap\thischap 
    \else 
     \theotherchap%
    \fi
   \else 
     \theotherpart.\theotherchap\fi}

\def\sectioncite[#1]#2{%
     \cprefix{#2}#1}

\edef\thispart{\partof{\thischap}}
\edef\thischap{\chapof{\thischap}}


\def\spuriousreset{}


\expandafter\ifx\csname citeadd.tex\endcsname\relax
\expandafter\gdef\csname citeadd.tex\endcsname{}
\else \message{Hey!  Apparently you were trying to
\string\input{citeadd.tex} twice.   This does not make sense.} 
\errmessage{Please edit your file (probably \jobname.tex) and remove
any duplicate ``\string\input'' lines}\endinput\fi

\def\sciteu{\sciteerror{undefined}}

\def\sciteerror#1#2{{\mathortextbf{\scite{#2}}}\complainaboutcitation{#1}{#2}}
\def\mathortextbf#1{\hbox{\bf #1}}
\def\complainaboutcitation#1#2{%
\vadjust{\line{\llap{---$\!\!>$ }\qquad scite$\{$#2$\}$ #1\hfil}}}

\sectno=-1   
\localtags
\jjtags
\NoBlackBoxes
\define\mr{\medskip\roster}
\define\sn{\smallskip\noindent}
\define\mn{\medskip\noindent}
\define\bn{\bigskip\noindent}
\define\ub{\underbar}
\define\wilog{\text{without loss of generality}}
\define\ermn{\endroster\medskip\noindent}

\define\dbcu{\dsize\bigcup}
\define \nl{\newline}
\magnification=\magstep 1
\documentstyle{amsppt}

{    
\catcode`@11

\ifx\alicetwothousandloaded@\relax
  \endinput\else\global\let\alicetwothousandloaded@\relax\fi

\gdef\subjclass{\let\savedef@\subjclass
 \def\subjclass##1\endsubjclass{\let\subjclass\savedef@
   \toks@{\def\usualspace{{\rm\enspace}}\eightpoint}%
   \toks@@{##1\unskip.}%
   \edef\thesubjclass@{\the\toks@
     \frills@{{\noexpand\rm2000 {\noexpand\it Mathematics Subject
       Classification}.\noexpand\enspace}}%
     \the\toks@@}}%
  \nofrillscheck\subjclass}
} 


\expandafter\ifx\csname alice2jlem.tex\endcsname\relax
  \expandafter\xdef\csname alice2jlem.tex\endcsname{\the\catcode`@}
\else \message{Hey!  Apparently you were trying to
\string\input{alice2jlem.tex}  twice.   This does not make sense.}
\errmessage{Please edit your file (probably \jobname.tex) and remove
any duplicate ``\string\input'' lines}\endinput\fi

\expandafter\ifx\csname bib4plain.tex\endcsname\relax
  \expandafter\gdef\csname bib4plain.tex\endcsname{}
\else \message{Hey!  Apparently you were trying to \string\input
  bib4plain.tex twice.   This does not make sense.}
\errmessage{Please edit your file (probably \jobname.tex) and remove
any duplicate ``\string\input'' lines}\endinput\fi

\def\renewcommand{\newcommand}	       
\edef\cite{\the\catcode`@}%
\catcode`@ = 11
\let\@oldatcatcode = \cite
\chardef\@letter = 11
\chardef\@other = 12
%
%
%
%
\def\@innerdef#1#2{\edef#1{\expandafter\noexpand\csname #2\endcsname}}%
%
%
\@innerdef\@innernewcount{newcount}%
\@innerdef\@innernewdimen{newdimen}%
\@innerdef\@innernewif{newif}%
\@innerdef\@innernewwrite{newwrite}%
%
%
%
\def\@gobble#1{}%
%
%
%
\ifx\inputlineno\@undefined
   \let\@linenumber = \empty 
\else
   \def\@linenumber{\the\inputlineno:\space}%
\fi
%
%
%
\def\@futurenonspacelet#1{\def\cs{#1}%
   \afterassignment\@stepone\let\@nexttoken=
}%
\begingroup 
\def\\{\global\let\@stoken= }%
\\ 
\endgroup
\def\@stepone{\expandafter\futurelet\cs\@steptwo}%
\def\@steptwo{\expandafter\ifx\cs\@stoken\let\@@next=\@stepthree
   \else\let\@@next=\@nexttoken\fi \@@next}%
\def\@stepthree{\afterassignment\@stepone\let\@@next= }%
%
%
%
\def\@getoptionalarg#1{%
   \let\@optionaltemp = #1%
   \let\@optionalnext = \relax
   \@futurenonspacelet\@optionalnext\@bracketcheck
}%
%
%
\def\@bracketcheck{%
   \ifx [\@optionalnext
      \expandafter\@@getoptionalarg
   \else
      \let\@optionalarg = \empty
      \expandafter\@optionaltemp
   \fi
}%
\def\@@getoptionalarg[#1]{%
   \def\@optionalarg{#1}%
   \@optionaltemp
}%
%
%
%
\def\@nnil{\@nil}%
\def\@fornoop#1\@@#2#3{}%
\def\@for#1:=#2\do#3{%
   \edef\@fortmp{#2}%
   \ifx\@fortmp\empty \else
      \expandafter\@forloop#2,\@nil,\@nil\@@#1{#3}%
   \fi
}%
\def\@forloop#1,#2,#3\@@#4#5{\def#4{#1}\ifx #4\@nnil \else
       #5\def#4{#2}\ifx #4\@nnil \else#5\@iforloop #3\@@#4{#5}\fi\fi
}%
\def\@iforloop#1,#2\@@#3#4{\def#3{#1}\ifx #3\@nnil
       \let\@nextwhile=\@fornoop \else
      #4\relax\let\@nextwhile=\@iforloop\fi\@nextwhile#2\@@#3{#4}%
}%
%
%
%
\@innernewif\if@fileexists
\def\@testfileexistence{\@getoptionalarg\@finishtestfileexistence}%
\def\@finishtestfileexistence#1{%
   \begingroup
      \def\extension{#1}%
      \immediate\openin0 =
         \ifx\@optionalarg\empty\jobname\else\@optionalarg\fi
         \ifx\extension\empty \else .#1\fi
         \space
      \ifeof 0
         \global\@fileexistsfalse
      \else
         \global\@fileexiststrue
      \fi
      \immediate\closein0
   \endgroup
}%
%
%
%
%
\def\bibliographystyle#1{%
   \@readauxfile
   \@writeaux{\string\bibstyle{#1}}%
}%
\let\bibstyle = \@gobble
%
%
\let\bblfilebasename = \jobname
\def\bibliography#1{%
   \@readauxfile
   \@writeaux{\string\bibdata{#1}}%
   \@testfileexistence[\bblfilebasename]{bbl}%
   \if@fileexists
      \nobreak
      \@readbblfile
   \fi
}%
\let\bibdata = \@gobble
%
%
\def\nocite#1{%
   \@readauxfile
   \@writeaux{\string\citation{#1}}%
}%
\@innernewif\if@notfirstcitation
%
%
\def\cite{\@getoptionalarg\@cite}%
%
%
\def\@cite#1{%
   \let\@citenotetext = \@optionalarg
   \printcitestart
   \nocite{#1}%
   \@notfirstcitationfalse
   \@for \@citation :=#1\do
   {%
      \expandafter\@onecitation\@citation\@@
   }%
   \ifx\empty\@citenotetext\else
      \printcitenote{\@citenotetext}%
   \fi
   \printcitefinish
}%
\newif\ifweareinprivate
\weareinprivatetrue
\ifx\shlhetal\undefinedcontrolseq\weareinprivatefalse\fi
\ifx\shlhetal\relax\weareinprivatefalse\fi
\def\@onecitation#1\@@{%
   \if@notfirstcitation
      \printbetweencitations
   \fi
   \expandafter \ifx \csname\@citelabel{#1}\endcsname \relax
      \if@citewarning
         \message{\@linenumber Undefined citation `#1'.}%
      \fi
     \ifweareinprivate
      \expandafter\gdef\csname\@citelabel{#1}\endcsname{%
\strut 
\vadjust{\vskip-\dp\strutbox
\vbox to 0pt{\vss\parindent0cm \leftskip=\hsize 
\advance\leftskip3mm
\advance\hsize 4cm\strut\openup-4pt 
\rightskip 0cm plus 1cm minus 0.5cm ?  #1 ?\strut}}
         {\tt
            \escapechar = -1
            \nobreak\hskip0pt\pfeilsw
            \expandafter\string\csname#1\endcsname
             \pfeilso
            \nobreak\hskip0pt
         }%
      }%
     \else  
      \expandafter\gdef\csname\@citelabel{#1}\endcsname{%
            {\tt\expandafter\string\csname#1\endcsname}
      }%
     \fi  
   \fi
   \csname\@citelabel{#1}\endcsname
   \@notfirstcitationtrue
}%
%
%
\def\@citelabel#1{b@#1}%
%
%
\def\@citedef#1#2{\expandafter\gdef\csname\@citelabel{#1}\endcsname{#2}}%
%
%
%
\def\@readbblfile{%
   \ifx\@itemnum\@undefined
      \@innernewcount\@itemnum
   \fi
   \begingroup
      \def\begin##1##2{%
         \setbox0 = \hbox{\biblabelcontents{##2}}%
         \biblabelwidth = \wd0
      }%
      \def\end##1{}
      %
      %
      \@itemnum = 0
      \def\bibitem{\@getoptionalarg\@bibitem}%
      \def\@bibitem{%
         \ifx\@optionalarg\empty
            \expandafter\@numberedbibitem
         \else
            \expandafter\@alphabibitem
         \fi
      }%
      \def\@alphabibitem##1{%
         \expandafter \xdef\csname\@citelabel{##1}\endcsname {\@optionalarg}%
         \ifx\biblabelprecontents\@undefined
            \let\biblabelprecontents = \relax
         \fi
         \ifx\biblabelpostcontents\@undefined
            \let\biblabelpostcontents = \hss
         \fi
         \@finishbibitem{##1}%
      }%
      \def\@numberedbibitem##1{%
         \advance\@itemnum by 1
         \expandafter \xdef\csname\@citelabel{##1}\endcsname{\number\@itemnum}%
         \ifx\biblabelprecontents\@undefined
            \let\biblabelprecontents = \hss
         \fi
         \ifx\biblabelpostcontents\@undefined
            \let\biblabelpostcontents = \relax
         \fi
         \@finishbibitem{##1}%
      }%
      \def\@finishbibitem##1{%
         \biblabelprint{\csname\@citelabel{##1}\endcsname}%
         \@writeaux{\string\@citedef{##1}{\csname\@citelabel{##1}\endcsname}}%
         \ignorespaces
      }%
      %
      %
      \let\em = \bblem
      \let\newblock = \bblnewblock
      \let\sc = \bblsc
      \frenchspacing
      \clubpenalty = 4000 \widowpenalty = 4000
      \tolerance = 10000 \hfuzz = .5pt
      \everypar = {\hangindent = \biblabelwidth
                      \advance\hangindent by \biblabelextraspace}%
      \bblrm
      \parskip = 1.5ex plus .5ex minus .5ex
      \biblabelextraspace = .5em
      \bblhook
      \input \bblfilebasename.bbl
   \endgroup
}%
%
%
\@innernewdimen\biblabelwidth
\@innernewdimen\biblabelextraspace
%
%
%
\def\biblabelprint#1{%
   \noindent
   \hbox to \biblabelwidth{%
      \biblabelprecontents
      \biblabelcontents{#1}%
      \biblabelpostcontents
   }%
   \kern\biblabelextraspace
}%
%
%
%
\def\biblabelcontents#1{{\bblrm [#1]}}%
%
%
\def\bblrm{\rm}%
%
%
\def\bblem{\it}%
%
%
\def\bblsc{\ifx\@scfont\@undefined
              \font\@scfont = cmcsc10
           \fi
           \@scfont
}%
%
%
\def\bblnewblock{\hskip .11em plus .33em minus .07em }%
%
%
\let\bblhook = \empty
%
%
%
\def\printcitestart{[}
\def\printcitefinish{]}
\def\printbetweencitations{, }
\def\printcitenote#1{, #1}
%
%
%
\let\citation = \@gobble
%
%
%
\@innernewcount\@numparams
%
%
\def\newcommand#1{%
   \def\@commandname{#1}%
   \@getoptionalarg\@continuenewcommand
}%
%
%
\def\@continuenewcommand{%
   \@numparams = \ifx\@optionalarg\empty 0\else\@optionalarg \fi \relax
   \@newcommand
}%
%
%
\def\@newcommand#1{%
   \def\@startdef{\expandafter\edef\@commandname}%
   \ifnum\@numparams=0
      \let\@paramdef = \empty
   \else
      \ifnum\@numparams>9
         \errmessage{\the\@numparams\space is too many parameters}%
      \else
         \ifnum\@numparams<0
            \errmessage{\the\@numparams\space is too few parameters}%
         \else
            \edef\@paramdef{%
               \ifcase\@numparams
                  \empty  No arguments.
               \or ####1%
               \or ####1####2%
               \or ####1####2####3%
               \or ####1####2####3####4%
               \or ####1####2####3####4####5%
               \or ####1####2####3####4####5####6%
               \or ####1####2####3####4####5####6####7%
               \or ####1####2####3####4####5####6####7####8%
               \or ####1####2####3####4####5####6####7####8####9%
               \fi
            }%
         \fi
      \fi
   \fi
   \expandafter\@startdef\@paramdef{#1}%
}%
%
%
%
%
\def\@readauxfile{%
   \if@auxfiledone \else 
      \global\@auxfiledonetrue
      \@testfileexistence{aux}%
      \if@fileexists
         \begingroup
            \endlinechar = -1
            \catcode`@ = 11
            \input \jobname.aux
         \endgroup
      \else
         \message{\@undefinedmessage}%
         \global\@citewarningfalse
      \fi
      \immediate\openout\@auxfile = \jobname.aux
   \fi
}%
%
%
\newif\if@auxfiledone
\ifx\noauxfile\@undefined \else \@auxfiledonetrue\fi
%
%
%
%
\@innernewwrite\@auxfile
\def\@writeaux#1{\ifx\noauxfile\@undefined \write\@auxfile{#1}\fi}%
%
%
%
\ifx\@undefinedmessage\@undefined
   \def\@undefinedmessage{No .aux file; I won't give you warnings about
                          undefined citations.}%
\fi
%
%
\@innernewif\if@citewarning
\ifx\noauxfile\@undefined \@citewarningtrue\fi
%
%
%
\catcode`@ = \@oldatcatcode

\def\pfeilso{\leavevmode
            \vrule width 1pt height9pt depth 0pt\relax
           \vrule width 1pt height8.7pt depth 0pt\relax
           \vrule width 1pt height8.3pt depth 0pt\relax
           \vrule width 1pt height8.0pt depth 0pt\relax
           \vrule width 1pt height7.7pt depth 0pt\relax
            \vrule width 1pt height7.3pt depth 0pt\relax
            \vrule width 1pt height7.0pt depth 0pt\relax
            \vrule width 1pt height6.7pt depth 0pt\relax
            \vrule width 1pt height6.3pt depth 0pt\relax
            \vrule width 1pt height6.0pt depth 0pt\relax
            \vrule width 1pt height5.7pt depth 0pt\relax
            \vrule width 1pt height5.3pt depth 0pt\relax
            \vrule width 1pt height5.0pt depth 0pt\relax
            \vrule width 1pt height4.7pt depth 0pt\relax
            \vrule width 1pt height4.3pt depth 0pt\relax
            \vrule width 1pt height4.0pt depth 0pt\relax
            \vrule width 1pt height3.7pt depth 0pt\relax
            \vrule width 1pt height3.3pt depth 0pt\relax
            \vrule width 1pt height3.0pt depth 0pt\relax
            \vrule width 1pt height2.7pt depth 0pt\relax
            \vrule width 1pt height2.3pt depth 0pt\relax
            \vrule width 1pt height2.0pt depth 0pt\relax
            \vrule width 1pt height1.7pt depth 0pt\relax
            \vrule width 1pt height1.3pt depth 0pt\relax
            \vrule width 1pt height1.0pt depth 0pt\relax
            \vrule width 1pt height0.7pt depth 0pt\relax
            \vrule width 1pt height0.3pt depth 0pt\relax}

\def\pfeilsw{ \leavevmode 
            \vrule width 1pt height0.3pt depth 0pt\relax
            \vrule width 1pt height0.7pt depth 0pt\relax
            \vrule width 1pt height1.0pt depth 0pt\relax
            \vrule width 1pt height1.3pt depth 0pt\relax
            \vrule width 1pt height1.7pt depth 0pt\relax
            \vrule width 1pt height2.0pt depth 0pt\relax
            \vrule width 1pt height2.3pt depth 0pt\relax
            \vrule width 1pt height2.7pt depth 0pt\relax
            \vrule width 1pt height3.0pt depth 0pt\relax
            \vrule width 1pt height3.3pt depth 0pt\relax
            \vrule width 1pt height3.7pt depth 0pt\relax
            \vrule width 1pt height4.0pt depth 0pt\relax
            \vrule width 1pt height4.3pt depth 0pt\relax
            \vrule width 1pt height4.7pt depth 0pt\relax
            \vrule width 1pt height5.0pt depth 0pt\relax
            \vrule width 1pt height5.3pt depth 0pt\relax
            \vrule width 1pt height5.7pt depth 0pt\relax
            \vrule width 1pt height6.0pt depth 0pt\relax
            \vrule width 1pt height6.3pt depth 0pt\relax
            \vrule width 1pt height6.7pt depth 0pt\relax
            \vrule width 1pt height7.0pt depth 0pt\relax
            \vrule width 1pt height7.3pt depth 0pt\relax
            \vrule width 1pt height7.7pt depth 0pt\relax
            \vrule width 1pt height8.0pt depth 0pt\relax
            \vrule width 1pt height8.3pt depth 0pt\relax
            \vrule width 1pt height8.7pt depth 0pt\relax
            \vrule width 1pt height9pt depth 0pt\relax
      }


\def\widestnumber#1#2{}

\def\citewarning#1{\ifx\shlhetal\relax 
    \else
    \par{#1}\par
    \fi
}

\def\rm{\fam0 \tenrm}

\def\fakesubhead#1\endsubhead{\bigskip\noindent{\bf#1}\par}



%
%
%

%

\font\textrsfs=rsfs10
\font\scriptrsfs=rsfs7
\font\scriptscriptrsfs=rsfs5

\newfam\rsfsfam
\textfont\rsfsfam=\textrsfs
\scriptfont\rsfsfam=\scriptrsfs
\scriptscriptfont\rsfsfam=\scriptscriptrsfs

\edef\oldcatcodeofat{\the\catcode`\@}
\catcode`\@11

\def\Cal@@#1{\noaccents@ \fam \rsfsfam #1}

\catcode`\@\oldcatcodeofat


\expandafter\ifx \csname margininit\endcsname \relax\else\margininit\fi

\long\def\red#1\endred{}
\long\def\green#1\endgreen{}
\long\def\blue#1\endblue{}

\def\endred{ \unmatched endred! }
\def\endgreen{ \unmatched endgreen! }
\def\endblue{ \unmatched endblue! }

\ifx\latexcolors\undefinedcs\def\latexcolors{}\fi

\def\emptycs{}
\def\evaluatelatexcolors{%
        \ifx\latexcolors\emptycs\else
        \expandafter\xxevaluate\latexcolors\xxfertig\evaluatelatexcolors\fi}
\def\xxevaluate#1,#2\xxfertig{\setupthiscolor{#1}%
        \def\latexcolors{#2}}

\font\smallfont=cmsl7
\def\rutgerscolor{\ifmmode\else\endgraf\fi\smallfont
\advance\leftskip0.5cm\relax}
\def\setupthiscolor#1{\edef\tmptmpcs{\noexpand\bgroup\noexpand\rutgerscolor
\noexpand\def\noexpand\currentcolor{#1}%
\noexpand}%
\expandafter\let\csname#1\endcsname\tmptmpcs
\def\tmptmpcs{\checkColorUnmatched{#1}\popthecolor}
\expandafter\let\csname end#1\endcsname\tmptmpcs}

\def\checkColorUnmatched#1{\def\expectcolor{#1}%
    \ifx\expectcolor\currentcolor   
    \else \edef\failhere{\noexpand\tryingToClose '\currentcolor' with end\expectcolor}\failhere\fi}

\def\currentcolor{???}

\def\popthecolor{\ifmmode\else\endgraf\fi\egroup}

\expandafter\def\csname#1\endcsname{}

\evaluatelatexcolors

 \let\outerhead\head
 \def\head{\innerhead}
 \let\innerhead\outerhead

 \let\outersubhead\subhead
 \def\subhead{\innersubhead}
 \let\innersubhead\outersubhead

 \let\outersubsubhead\subsubhead
 \def\subsubhead{\innersubsubhead}
 \let\innersubsubhead\outersubsubhead

 \def\proclaim{\innerproclaim}
 \let\innerproclaim\outerproclaim

 %
 %
 %
 %

\def\demo#1{\medskip\noindent{\it #1.\/}}
\def\enddemo{\smallskip}

\def\remark#1{\medskip\noindent{\it #1.\/}}
\def\endremark{\smallskip}

\pageheight{8.5truein}
\topmatter
\title{More on: The revised GCH and middle diamond} \endtitle
\author {Saharon Shelah \thanks {\null\newline 
The author would like to thank the United States-Israel Binational
Science Foundation for partial support of this research and
Alice Leonhardt for the beautiful typing. Publication 829} \endthanks}
\endauthor  

\affil{The Hebrew University of Jerusalem \\
Einstein Institute of Mathematics \\
Edmond J. Safra Campus, Givat Ram \\
Jerusalem 91904, Israel
 \medskip
 Department of Mathematics \\
 Hill Center-Busch Campus \\
  Rutgers, The State University of New Jersey \\
 110 Frelinghuysen Road \\
 Piscataway, NJ 08854-8019 USA} \endaffil

\abstract  We strengthen the revised GCH theorem by showing, e.g.,
that for $\lambda = \text{cf}(\lambda) > \beth_\omega$, for all but
finitely many regular $\kappa <\beth_\omega,\lambda$ is accessible on
cofinality $\kappa$ in a weak version of it holds.
In particular, $\lambda = 2^\mu = \mu^+ > \beth_\omega$ implies the
diamond on $\lambda$ is restricted to cofinality $\kappa$ for all but
finitely many $\kappa \in \text{ Reg} \cap \beth_\omega$ and we
strengthen the results on the middle diamond.  Moreover, we get
stronger results on the middle diamond.
\endabstract
\endtopmatter
\document

\newpage

\head {\S0 Introduction} \endhead  \resetall \sectno=0
 \spuriousreset
\bigskip

The main result of this paper is to define for any cardinal $\lambda$
a set ${\frak d}_0(\lambda)$ of regular cardinals $< \lambda$ such that for
strong limit $\theta < \lambda$ we prove that
$\theta \cap {\frak d}_0(\lambda)$ is
finite, and for every $\kappa \in \text{ Reg } \cap \theta \backslash {\frak
d}_0(\lambda)$, in some sense $\lambda$ has not too many subsets of
cardinality $\kappa$.  This serves as our main aim here to use this to
show: if cf$(\lambda) > \mu$ and $\kappa \in \text{ Reg } \cap \mu$ satisfies
$\lambda = \sup\{\alpha:\kappa \notin {\frak d}_0(\alpha)\}$
\ub{then} $\lambda$ has a ``good" sequence $\langle {\Cal
P}_\alpha:\alpha < \lambda \rangle,{\Cal P}_\alpha \subseteq
[\alpha]^{\le \kappa}$ and if $\lambda = \lambda^\mu$ more (see
\scite{md.2}, \scite{md.5}).  \nl
This gives as a main consequence that: if $\mu \ge \theta,\lambda = \text{
cf}(2^\mu)$ \ub{then} $(\lambda,\kappa)$ has middle diamond for all
but finitely many $\kappa \in \text{ Reg}$ satisfying
$\beth_\omega(\kappa) \le \mu$.  Also $\lambda =
2^\mu = \mu^+ > \theta \Rightarrow \lambda$ has diamond on cofinality
$\kappa$ for all regular $\kappa$ for which $\beth_\omega(\kappa) < \lambda$ 
except finitely many.  We also strengthen the results on the middle
diamond (in ZFC, see below).

So this is part of pcf theory (\cite{Sh:g}) continuing in particular
\cite{Sh:460}.  As I have been very happy about \cite{Sh:460} 
and the proofs here give a shorter 
proof of the main theorem there, we have here (\S1) gives a
self-contained proof of the revised GCH, the main theorem of
\cite{Sh:460} and discuss it.

By pcf theory (\cite{Sh:g},\cite{Sh:460}) a worthwhile choice of power
is (for $\kappa < \lambda$ regular) $\lambda^{[\kappa]}$ (or 
$\lambda^{<\kappa>}$), the minimal
cardinality of a family of subsets of $\lambda$ each of cardinality
$\le \kappa$ such that any other subset of $\lambda$ of cardinality
$\kappa$ is equal to (or is contained in) 
the union of $< \kappa$ members of the family (see
Definition \scite{g.2}).  Let ${\frak d}^+(\lambda) = \{\kappa:\kappa$
is regular $< \lambda$ and $\lambda < \lambda^{<\kappa>}\}$.

This gives a good partition of the exponentiation as $\lambda^\kappa =
\lambda \Leftrightarrow 
2^\kappa \le \lambda \and (\forall \sigma)(\sigma = \text{
cf}(\sigma) \le \kappa \Rightarrow \lambda^{<\sigma>} = \lambda)$.  So
G.C.H. is equivalent to: $\kappa$ regular $\Rightarrow 2^\kappa =
\kappa^+$ and [$\kappa < \lambda$ are regular $\Rightarrow
\lambda^{<\kappa>} = \lambda$].  See more in \cite{Sh:460} where it is
proved that (the revised G.C.H.):
\mr
\item "{$\circledast$}"  if $\lambda > \beth_\omega$ then ${\frak
d}^+(\lambda) \cap \beth_\omega$ is bounded, i.e., $\lambda =
\lambda^{<\kappa>}$ for every regular $\kappa < \beth_\omega$ large
enough.
\ermn
We can replace $\beth_\omega$ by any strong limit singular cardinal $\theta$.

The advances in pcf theory reveal several natural hypotheses.  The
Strongest Hypothesis (pp$(\mu)=\mu^+$ for every singular $\mu$) is
very nice but it implies the SCH hence it does not follow from ZFC.
The status of the Weak Hypothesis (somewhat more than $\{\mu:\text{cf}(\mu) <
\mu < \lambda \le \text{ pp}(\mu)\}$ is at most countable) is not
known but I am sure is consistent though it has large consistency
strength.  I am not sure about $(\forall {\frak a})(|{\frak a}| \ge
|\text{pcf}({\frak a})|)$.  
Still better then $\circledast$ would be 
(we believe but do not know it, particularly (2)).
\bigskip

\demo{\stag{0.1} Conjecture}  1) For every $\lambda,{\frak
d}^+(\lambda)$ is finite, or at least
\nl
2) For every strong limit $\mu,\lambda \ge \mu \Rightarrow 
{\frak d}^+(\lambda) \cap \mu$ is finite.
\nl
Now ${\frak d}_0(\lambda) \cap \theta$ being finite is a step in the
right direction and is enough to improve the results of
\cite{Sh:775}.  In particular being able to use $\kappa = \aleph_n$
for some $n$ rather than $\kappa$ regular $< \beth_\omega$ seems
crucial in abelian group theory (as there are non-free almost
$\kappa$-free abelian group of cardinality $\kappa$ when 
$\kappa = \aleph_n$).

So we can hope to get the right objects in each cardinality $\aleph_n$
whereas consistently they may not exist for arbitrary $\kappa = 
\text{ cf}(\kappa) < \beth_\omega$, this is the case in abelian group theory.

This work also continues ones on $I[\lambda]$.  By \cite{Sh:108} if
$\theta$ is strong limit singular, for some $A \in I[\lambda]$ for
some $\bold c:[\mu^+]^2 \rightarrow \text{ cf}(\mu)$, if $B \subseteq
\mu,\bold c \restriction [B]^2$ constant (or just has bounded range),
$\delta = \sup(B)$, cf$(\delta) \ne \text{ cf}(\mu)$ then $\delta \in
A$.

By Dzamanja and Shelah \cite{DjSh:562}, using \cite{Sh:460}, if
$\lambda = \text{ cf}(\lambda) > \theta$, a strong limit singular, for
some $\kappa <\theta$, for some $A \in I[\lambda]$, if for every $A'
\subseteq A,|A'| < \theta$ for some $\bold c:[A'] \rightarrow \kappa$,
we have: if $B \subseteq A',\bold c \restriction [B]$ is constant,
$\delta = \sup[B]$, cf$(\delta) > \kappa$ then $\delta \in A$.  By
\cite[5.20]{Sh:589}, $T_D$ helps.
\nl
We try to make this paper as self-contained as is reasonably possible.
\enddemo
\bigskip

\definition{\stag{0.A} Definition}  1) For an ideal $J$ on a set $X$:
\mr
\item "{$(a)$}"  $J^+ = {\Cal P}(X) \backslash J$; we agree that $J$
determines $X$ so $X = \text{ Dom}(J)$, this is an abuse of notation
when $\cup\{A:A \in J\} \subset X$ but usually clear in the context
\sn
\item "{$(b)$}"  for a two-place relation $R$ on $Y$ and for 
$f,g \in {}^X Y$, let $f R_J g$ means $\{t \in X:\neg f(t) R g(t)\} \in
J$; the cases we shall use are $=,\ne,<,\le$.
\ermn
2) If $D$ is a filter on $X,J$ the dual ideal on $X$ (i.e., $J = \{X
\backslash A:A \in D\}$) we may replace $J$ by $D$.
\nl
3) Let $(\forall^J t)\varphi(t)$ mean $\{t:\neg \varphi(t)\} \in J$
similarly $\exists^J,\forall^D,\exists^D$.
\nl
4) Let $S^\lambda_\kappa = \{\delta < \lambda:\text{cf}(\delta) =
\kappa\}$ and $S^\lambda_{< \kappa} = \{\delta <
\lambda:\text{cf}(\delta) < \kappa\}$.
\enddefinition
\bn
We quote (see \cite{Sh:71}, \cite[\S3]{Sh:506}, in \scite{3.13},
\scite{0.Q}, \scite{0.B})).
\definition{\stag{3.13} Definition}  1) Let $\bar A = \langle A_i:i
\in X \rangle,D$ a filter on $X$ and for simplicity first assume
$i \in X \Rightarrow A_i \ne \emptyset$.  We let
\mr
\item "{$(a)$}"   $T_D^0(\bar A) = \sup\{|{\Cal F}|:{\Cal F} \subseteq
\Pi(\bar A)$ and $f_1 \ne 
f_2 \in {\Cal F} \Rightarrow f_1\not=_D f_2\}$
\sn
\item "{$(b)$}"   
$$
\align
T^1_D(\bar A)=\text{ Min}\{|{\Cal F}|:&(i) \quad {\Cal F} 
\subseteq \Pi(\bar A) \\
  &(ii) \quad f_1 \ne f_2\in {\Cal F} \Rightarrow f_1 \ne_D f_2 \\
  &(iii) \quad {\Cal F} \text{ maximal under } (i)+(ii)\}
\endalign
$$
\item "{$(c)$}"   $T^2_D(\bar A) = \text{ Min}\{|{\Cal F}|:{\Cal F}
\subseteq \Pi \bar A$ and for every $f_1\in \Pi \bar A$, for some
$f_2 \in {\Cal F}$ we have $\neg(f_1 \ne_D f_2)\}$.
\ermn
2) If $\{i:A_i = \emptyset\} \in J$ then we let $T^\ell_D(\bar A) = T^\ell_{D
\restriction Y}(\bar A \restriction Y)$ where 
$Y = \{i:A_i \ne \emptyset\}$; if $\{i:A_i \ne \emptyset\} \in J$ 
then $T^\ell_D(\bar A) = 0$. 
\nl
3) For $f\in {}^\kappa\text{Ord}$ and 
$\ell<3$ let $T^l_D(f)$ means $T^l_D(\langle f(\alpha):\alpha<\kappa
\rangle)$.
\nl
4) If $T^0_D(\bar A) = T^1_D(\bar A)=T^2_D(\bar A)$ then let 
$T_D(\bar A) = T^l_D(\bar A)$ for $l < 3$; similarly $T_D(\bar A)$; 
we say that ${\Cal F}$ witness $T_D(\bar A) = \lambda$ if it is as in the
definition of $T^1_D(\bar A) =\lambda$;  similarly $T_D(f)$.
\enddefinition
\bigskip

\remark{Remark}  Actually the case $\bar A = \langle
\lambda_\alpha:\alpha < \lambda \rangle$ is enough so we concentrate
on it.
\endremark
\bigskip

\proclaim{\stag{0.Q} Claim}  0)  If $D_0\subseteq D_1$ are 
filters on $\kappa$ then $T^{\ell}_{D_0}(\bar \lambda)\le 
T^{\ell}_{D_1}(\bar \lambda)$ for $\ell=0,2$. 
\nl
1) $T^2_D(\bar \lambda) \le T^1_D(\bar \lambda) \le T^0_D(\bar \lambda)$,
in particular $T^1_D(\bar \lambda)$ is well defined.
\nl
2) If $(\forall i)\lambda_i > 2^\kappa$
\ub{then} $T^0_D(\bar\lambda)=T^1_D(\bar\lambda)=T^2_D(\bar\lambda)$ so
the supremum in \scite{3.13}(a) is obtained (so, e.g.,
$T^0_D(\bar\lambda)>2^\kappa$ suffice; also $(\forall i)\lambda_i \ge
2^\kappa$ suffice). 
\endproclaim
\bigskip

\demo{Proof}  0) Check. \nl
1)  If ${\Cal F}$ exemplifies the value of $T^1_D(\bar\lambda)$, it also
exemplifies $T^2_D(\bar\lambda)\le |{\Cal F}|$ hence but easily 
$T^2_D(\bar\lambda)\le T^1_D(\bar\lambda)$ is as in the definition it
can be extended.
In the definition of $T^0_D(\bar \lambda)$ the Min is taken over a
non-empty set (as maximal such ${\Cal F}$ exists), so $T^0_D(\bar
\lambda)$ is well defined as well as $T^1_D(\bar \lambda)$.

Lastly, if ${\Cal F}$ exemplifies the value of $T^1_D(\bar \lambda)$ it also
exemplifies $T^0_D(\bar\lambda)\ge |{\Cal F}|$, so
$T^1_D(\bar\lambda)\le T^0_D(\bar\lambda)$. \nl
2) Let $\mu$ be $2^\kappa$.
Assume that the desired conclusion fails so $T^2_D(\bar\lambda)<
T^0_D(\bar\lambda)$, so there is ${\Cal F}_0 \subseteq \Pi 
\bar\lambda$, such that $[f_1\not= f_2\in {\Cal F}_0 \Rightarrow 
f_1\not=_D f_2]$, and $|{\Cal F}_0| > T^2_D(\bar\lambda) + \mu$ 
(by the definition of $T^0_D(\bar\lambda)$).
Also there is ${\Cal F}_2 \subseteq \Pi \bar\lambda$ exemplifying the value
of $T^2_D(\bar\lambda)$.
For every $f \in {\Cal F}_0$ there is $g_f\in {\Cal F}_2$ such that
$\neg(f \ne_D g_f)$ (by the choice of ${\Cal F}_2$).
As $|{\Cal F}_0|> T^2_D(\bar\lambda) + \mu$ for some $g\in
{\Cal F}_2$, ${\Cal F}^*=:\{ f\in {\Cal F}_0:g_f = g\}$ has 
cardinality $>T^2_D(f)+\mu$.
Now for each $f\in {\Cal F}^*$ let $A_f=\{ i<\kappa:f(i)=g(i)\}$
clearly $A_f\in D^+$. Now $f\mapsto A_f/D$ is a function from 
${\Cal F}^*$ into ${\Cal P}(\kappa)/D$, hence, if as $\mu \ge
|{\Cal P}(\kappa)/D|$, it is not one to one
(by cardinality consideration) so for some $f' \ne f''$ from 
${\Cal F}^*$ (hence form ${\Cal F}_0$) we have $A_{f'} /D= A_{f''}/D$; but so

$$
\{i<\kappa:f'(i) = f''(i)\}\supseteq \{ i<\kappa:
f'(i) =g(i)\} \cap \{ i<\kappa:f''(i) = g(i)\} = A_{f'}/D
$$
\mn
hence is $\ne \emptyset \mod D$, so $\neg(f' \ne_D f'')$,
contradition the choice of ${\Cal F}_0$.  \hfill$\square_{\scite{0.Q}}$\margincite{0.Q}
\enddemo
\bigskip

\proclaim{\stag{0.B} Claim}  Let $J$ be a $\sigma$-complete ideal on
$\kappa$.
\nl
1) If $\bar A = \langle A_i:i < \kappa \rangle,\bar \lambda = \langle
\lambda_i:i < \kappa \rangle,\lambda_i = |A_i|$ \ub{then} $T^\ell_J(\bar A) =
T^\ell_J(\bar \lambda)$ and if $A \in J,B = \kappa \backslash A$ then
$T^\ell_J(\bar \lambda) = T^\ell_{J \restriction B}(\bar \lambda \restriction
B)$.
\nl
2) $T_J(\bar \lambda) > 2^\kappa$ iff $(\forall_J t)(\lambda_t > 2^\kappa)$.
\nl
3) $T^\ell_J(\bar \lambda^1) \le T^\ell_J(\bar \lambda^2)$ if $(\forall_J t)
(\lambda^1_t \le \lambda^2_t)$.
\nl
4) If {\rm Dom}$(J) = \cup\{A_\varepsilon:\varepsilon < \zeta\},\zeta <
\sigma$ and $\lambda_i \ge 2^\kappa$ for $i < \kappa$ 
\ub{then} $T^0_J(\bar \lambda) = { \text{\rm Min\/}}\{T^0_{J \restriction
A_\varepsilon}(\bar \lambda \restriction A_\varepsilon):
\varepsilon < \zeta$ and $A_\varepsilon \in J^+\}$.
\endproclaim
\bigskip

\demo{Proof}  E.g. (and the one we use): 
\nl
4) Let $A'_\varepsilon =
A_\varepsilon \backslash \cup \{A_\xi:\xi < \varepsilon\}$ for
$\varepsilon < \zeta$.
\nl
First assume that ${\Cal F} \subseteq \Pi \bar \lambda$ and $f_1
\ne f_2 \in {\Cal F}_2 \Rightarrow f_1 \ne_J f_2$.  Then for each
$\varepsilon < \zeta$ satisfying $A_\varepsilon \in J^+$ clearly
${\Cal F}^{[\varepsilon]} = \{f \restriction A_\varepsilon:f \in {\Cal
F}\}$ satisfies $|{\Cal F}^{[\varepsilon]}| = |{\Cal F}|$ as $f
\mapsto f \restriction A_\zeta$ is one to one by the assumption on
${\Cal F}$ and ${\Cal F}^{[\varepsilon]} \subseteq \dsize \prod_{i \in
A_\zeta} \lambda_i$; so $|{\Cal F}| = |{\Cal F}^{[\varepsilon]}| \le
T_{J \restriction A_\varepsilon}(\bar \lambda \restriction
A_\varepsilon)$.  As this holds for every $\varepsilon < \zeta$ for
which $A_\varepsilon \in J^+$ we get $|{\Cal F}| \le 
\text{ Min}\{T_{J \restriction A_\varepsilon}(\bar \lambda \restriction
A_\varepsilon):\varepsilon < \zeta,A_\varepsilon \in J^+\}$. By the
demand on ${\Cal F}$ we get the inequality $\le$ in part (4).  Second,
assume $\mu < \text{ Min}\{T^0_{J \restriction A_\varepsilon}(\bar
\lambda \restriction A_\varepsilon):\varepsilon < \zeta,A_\varepsilon
\in J^+\}$.  So for each such $\varepsilon$ there is 
${\Cal F}_\varepsilon \subseteq \dsize \prod_{i \in A_\varepsilon}
\lambda_i$ such that $f \ne g \in {\Cal F}_\varepsilon \Rightarrow f
\ne_{J \restriction A_\varepsilon} g,|{\Cal F}_\varepsilon| \ge
\mu^+$.  Let $f^\varepsilon_\alpha \in {\Cal F}_\varepsilon$ be
pairwise disinct, and define $f_\alpha \in \Pi \bar \lambda$ for
$\alpha < \lambda^+$ as follows $f_\alpha \restriction A'_\varepsilon =
f^\varepsilon_\alpha$ when $A_\varepsilon \in J^+,f_\alpha
\restriction A'_\xi$ is zero otherwise.  
\nl
Now check.
\enddemo
\bigskip

\definition{\stag{0.R} Definition}  For regular $\lambda$ and
stationary $S \subseteq \lambda$ let $(D\ell)_{\lambda,S}$ mean that
we can find $\bar{\Cal P} = \langle {\Cal P}_\alpha:\alpha \in S
\rangle,{\Cal P}_\alpha \subseteq {\Cal P}(\alpha)$ of cardinality $<
\lambda$ such that for every $A \subseteq \lambda$ the set $\{\alpha \in S:A
\cap \alpha \in {\Cal P}_\alpha\}$ is stationary.
\enddefinition
\bigskip

\definition{\stag{0.S} Definition}   For $\lambda$ regular uncountable
let $I[\lambda]$ be the family of sets $S \subseteq \lambda$ which have
a witness $(E,\bar{\Cal P})$ for $S \in I[\lambda]$, which means
\mr
\item "{$(*)$}"  $E$ is a club of $\lambda,\bar{\Cal P} = \langle {\Cal
P}_\alpha:\alpha < \lambda \rangle,{\Cal P}_\alpha \subseteq {\Cal
P}(\alpha),|{\Cal P}_\alpha| < \lambda$, and for every $\delta \in E
\cap S$ there is an unbounded subset $C$ of $\delta$ of order type $<
\delta$ satisfying 
$\alpha \in C \Rightarrow C \cap \alpha \in {\Cal P}_\alpha$.
\endroster
\enddefinition
\bigskip

\proclaim{\stag{0.51} Claim}  (\cite{Sh:420}): 1) For $\lambda$ regular
uncountable, $S \in I[\lambda]$ iff there is a 
pair $(E,\bar{\Cal P}),E$ a club of 
$\lambda,\bar a = \langle a_\alpha:\alpha < \lambda \rangle,
a_\alpha \subseteq \alpha$ such that $\beta \in a_\alpha \Rightarrow
 a_\beta = a_\alpha \cap \beta$ and $\delta \in E \cap S \Rightarrow
 \delta = \sup(a_\delta) > { \text{\rm otp\/}}(a_\delta)$ (or even $\delta \in
E \cap S \Rightarrow \delta = \sup(a_\delta))$, {\rm otp}$(a_\delta) 
= { \text{\rm cf\/}}(\delta) < \delta$.
\nl
2) If $\kappa^+ < \lambda$ and $\lambda,\kappa$ are regular, \ub{then}
for some stationary $S \in I[\lambda]$ we have $\delta \in S
\Rightarrow { \text{\rm cf\/}}(\delta) = \kappa$.
\endproclaim
\bigskip

\proclaim{\stag{0.K} Claim}  1) Assume that $f_\alpha \in
{}^\kappa${\rm Ord} for $\alpha < \lambda,\lambda = (2^\kappa)^+$ or
just $\lambda =$ {\rm cf}$(\lambda)$ and $(\forall \alpha <
\lambda)(|\alpha|^\kappa < \lambda)$ and $S_1 \subseteq \{\delta <
\lambda$: {\rm cf}$(\delta) > \kappa\}$ is stationary.  \ub{Then} for some
stationary $S_2 \subseteq S_1$ we have: for each $i < \kappa$ the
sequence $\langle f_\alpha(i):\alpha \in S_2 \rangle$ is increasing or
is constant.
\nl
2) If $D$ is a filter on $\kappa$ and
$f_\alpha \in {}^\kappa${\rm Ord} for $\alpha < \delta$ is
$<_D$-increasing and {\rm cf}$(\delta) > 2^\kappa$ then $\langle
f_\alpha:\alpha < \delta \rangle$ has a $<_D$-eub, i.e.,
\mr
\item "{$(i)$}"  $\alpha < \delta \Rightarrow f_\alpha \le_D f_\delta$
\sn
\item "{$(ii)$}"  $f' \in {}^\kappa${\rm Ord} $\and f' <_D$ {\rm
Max}$\{f,1_\kappa\}$ then $(\exists \alpha < \delta)(f' <_D f_\alpha)$.
\endroster
\endproclaim
\bigskip

\demo{Proof}  See \cite{Sh:110}, \cite{Sh:g}.
\enddemo
\newpage

\head {\S1 The revised G.C.H. Revisited} \endhead  \resetall \sectno=1
 \spuriousreset
\bigskip

Here we give a proof of the RGCH which requires little knowledge; this
is the main theorem of \cite{Sh:460}, see also \cite[\S1]{Sh:513}.
The presentation is self-contained; in particular, do not use the
pcf-theorem (hence repeat some proofs, in weak forms).
\bigskip

\definition{\stag{g.1} Definition}  1) For $\lambda > \theta \ge
\sigma = \text{ cf}(\sigma)$ let 
$\lambda^{[\sigma,\theta]} = \text{ Min}\{|{\Cal P}|:{\Cal P} 
\subseteq [\lambda]^{\le \theta}$ and every $u \in [\lambda]^{\le
\theta}$ is the union of $< \sigma$ members of ${\Cal P}\}$.
\nl
2) Let $\lambda^{[\sigma]} = \lambda^{[\sigma,\sigma]}$.
\nl
3) For $\lambda \ge \theta^{[\sigma,\kappa]}$ let
$\lambda^{[\sigma,\kappa,\theta]} = \text{ Min}\{|{\Cal P}|:{\Cal P} 
\subseteq [\lambda]^{\le \kappa}$ such that for every $u \subseteq
\lambda$ of cardinality $\le \theta$ we can find $i^* < \sigma$ and
$u_i \subseteq u$ for $i<i^*$ such that $u= \cup\{u_i:i <i^*\}$ and
$[u_i]^\kappa \subseteq {\Cal P}\}$.
\nl
4) We may replace $\theta$ by $< \theta$ with the obvious meaning
(also $< \kappa$).  
\enddefinition
\bigskip

\definition{\stag{g.2} Definition}  1) For $\lambda > \theta \ge 
\text{ cf}(\sigma) = \sigma$ let 
$\lambda^{<\sigma,\theta>} = \text{ Min}\{|{\Cal P}|:{\Cal P} 
\subseteq [\lambda]^\theta$ and every $u \in [\lambda]^{\le
\theta}$ is included in the union of $< \sigma$ members of ${\Cal P}\}$.
\nl
2) Let $\lambda^{<\sigma>} = \lambda^{<\sigma,\sigma>}$.
\nl
3) For $\lambda \ge \theta^{<\sigma,\kappa>}$ let
$\lambda^{<\sigma,\kappa,\theta>} = \text{ Min}\{|{\Cal P}|:{\Cal P} 
\subseteq [\lambda]^{\le \kappa}$ such that for every $u \subseteq
\lambda$ of cardinality $\le \theta$ we can find $i^* < \sigma$ and
$u_i \subseteq u$ for $i<i^*$ such that $u \subseteq \cup
\{u_i:i <i^*\}$ satisfying 
$(\forall v \in [u_i]^{\le \kappa})
(\exists w \in {\Cal P})(v \subseteq w)\}$.
\nl
4) We may replace $\theta$ by $< \theta$ with the obvious meaning
(also $< \kappa$).  
\enddefinition
\bigskip

\demo{\stag{g.3} Observation}  Let $\lambda > \theta \ge \kappa \ge
\sigma = \text{ cf}(\sigma)$.
\nl
1) $\lambda^{<\kappa>} \le \lambda^{[\kappa]} \le 
\lambda^{<\kappa>} + 2^\kappa$.
\nl
2) $\lambda^{<\sigma,\theta>} \le \lambda^{[\sigma,\theta]} 
\le \lambda^{<\sigma,\theta>} + 2^\theta$
(but see (3)).
\nl
3) If cf$(\theta) < \sigma$ then $\lambda^{<\sigma,\theta>} =
 \Sigma\{\lambda^{<\sigma,\theta'>}:\sigma \le \theta' < \theta\}$
 and $\lambda^{[\sigma,\theta]} =
 \Sigma\{\lambda^{[\sigma,\theta']}:\sigma \le \theta' < \theta\}$.
\nl
4) $\lambda^{<\sigma,\kappa,\theta>} \le \lambda^{[\sigma,\kappa,\theta]} \le
\lambda^{<\sigma,\kappa,\theta>} + 2^\kappa$.
\enddemo
\bigskip

\demo{Proof}  Easy.
\enddemo
\bn
The main claim of this section is
\proclaim{\stag{g.4} Claim}  Assume
\mr
\item "{$(a)$}"  $\aleph_0 < \sigma = { \text{\rm cf\/}}(\sigma) \le
\kappa < \partial \le \theta$
\sn
\item "{$(b)$}"  $J$ is a $\sigma$-complete ideal on $\kappa$
\sn
\item "{$(c)$}"  $\bar \lambda = \langle \lambda_i:i < \kappa \rangle$
\sn
\item "{$(d)$}"  $T_J(\bar \lambda) = \lambda$ 
\sn
\item "{$(e)$}"  $\lambda^{[\partial,\theta]}_i =
\lambda_i$ for $i < \kappa$ (yes $\partial$ not $\partial_i$!)
\sn
\item "{$(f)$}"   if $\partial_i < \partial$
for $i < \kappa$ then $\dsize \prod_{i < \kappa} \partial_i < 
\partial$
\sn
\item "{$(g)$}"  $\theta = \theta^\kappa$ and $2^\theta \le \lambda$.
\ermn
\ub{Then} $\lambda^{[\partial,\theta]} = \lambda$.
\endproclaim
\bigskip

\remark{Remark}  1) We may consider using a $\mu^+$-free family $\bar f$ 
(see \S2). \nl
2) Actually we use less than $T^1_J(\bar \lambda) = \lambda$, we
just use
\mr
\item "{$(a)$}"  there are $f_\alpha \in \dsize \prod_{i < \kappa}
\lambda_i$ for $\alpha < \lambda$ such that $\alpha < \beta
\Rightarrow f_\alpha \ne_J f_\beta$
\sn
\item "{$(b)$}"  there are $f_\alpha \in \dsize \prod_{i < \kappa}
\lambda_i$ for $\alpha < \lambda$ such that for every $f \in \dsize
\prod_{i < \kappa} \lambda_i$ for some $\alpha,\neg(f \ne_J
f_\alpha)$.
\endroster
\endremark
\bigskip

\demo{Proof}  Let $\bar f = \langle f_\alpha:\alpha < \lambda \rangle$
be pairwise $J$-different, $f_\alpha \in \dsize \prod_{i < \kappa}
\lambda_i$ (i.e. $\alpha \ne \beta \Rightarrow \{i:f_\alpha(i) =
f_\beta(i)\} \in J$).

For each $i < \kappa$ let ${\Cal P}_i \subseteq [\lambda_i]^{\le
\theta}$ be of cardinality $\lambda_i$ and witness
$\lambda^{[\partial,\theta]}_i = \lambda_i$ that is:  every $u \in
[\lambda_i]^{\le \theta}$ is the union of $< \partial$ members of 
${\Cal P}_i$; such family exists by assumption (e).  
Let $M \prec ({\Cal H}(\chi),\in)$ be
of cardinality $\lambda$ such that $\lambda +1 \subseteq M$ and $\bar f,
\langle \lambda_i:i < \kappa \rangle,\langle {\Cal P}_i:i <
\kappa \rangle,J,{\Cal P}(\kappa)$ belong to $M$.
\nl
Let ${\Cal P} = M \cap [\lambda]^{\le \theta}$.  We shall show that
${\Cal P}$ exemplifies the desired conclusion.  Now ${\Cal P}$ is a
family of $\le \|M\| = \lambda$ subsets of $\lambda$ each of
cardinality $\le \theta$ hence it is enough to show
\mr
\item "{$(*)$}"  if $u \in [\lambda]^{\le \theta}$ then $u$ is
included in the union of $< \partial$ members of ${\Cal P}$ (or equal to; 
equivalent here as $2^\theta \le \lambda$ hence $u_1 \subseteq u_2 \in {\Cal
P} \Rightarrow u_1 \in {\Cal P}$).
\endroster
\enddemo
\bigskip

\demo{Proof of $(*)$}  Let $u_i = \{f_\alpha(i):\alpha \in u\}$; so
$u_i \in [\lambda_i]^{\le \theta}$, hence we can find $\langle
v_{i,j}:j < j_i \rangle$ such that $v_{i,j} \in {\Cal P}_i$ and $u_i =
\cup\{v_{i,j}:j < j_i\}$ and $0 < j_i < \partial$.  For each $\eta \in
\dsize \prod_{i < \kappa} j_i$ let 

$$
w_\eta = \{\alpha \in u:i <\kappa \Rightarrow f_\alpha(i) \in
v_{i,\eta(i)}\}.
$$
\mn
Clearly $u = \cup\{w_\eta:\eta \in \dsize \prod_{i<\kappa} j_i\}$ as
for any $\alpha \in u$ for each $i < \kappa$ 
we can define $\varepsilon_i(\alpha) < j_i$ such that $f_\alpha(i) \in
v_{i,\varepsilon_i(\alpha)}$ and let $\eta_\alpha = \langle
\varepsilon_i(\alpha):i < \kappa \rangle$, clearly $\eta_\alpha \in
\dsize \prod_{i < \kappa} j_i$ and so $\alpha \in w_{\eta_\alpha}$.  By
the assumption (f) as $i < \kappa \Rightarrow j_i < \partial$ clearly
$|\dsize \prod_{i < \kappa} j_i| < \partial$ hence it is enough to
prove that $\eta \in \dsize \prod_{i < \kappa} j_i \Rightarrow w_\eta
\in {\Cal P}$.  As $u \in M \wedge |u| \le \theta \Rightarrow {\Cal P}(u)
\subseteq M$ it is
enough to prove for $\eta\in \dsize \prod_{i < \kappa} j_i$, that
\mr
\item "{$\circledast$}"  $w_\eta$ is included in some $w \in M \cap
[\lambda]^{\le \theta}$.
\ermn
\ub{Proof of $\circledast$}:  As $i < \kappa \Rightarrow |{\Cal P}_i|
= \lambda_i$ and $T_J(\bar \lambda) = \lambda$ by \scite{0.B} 
there is ${\Cal G}
\subseteq \{g:g \in \dsize \prod_{i < \kappa} {\Cal P}_i\}$ satisfying
$|{\Cal G}| = \lambda$ and $(\forall g \in \dsize \prod_{i < \kappa} {\Cal
P}_i)(\exists g' \in {\Cal G})(\{i:g(i) = g'(i)\} \in J^+)$.  As
$\langle {\Cal P}_i:i < \kappa \rangle \in M$ \wilog \, ${\Cal G}
\in M$ and as $\lambda +1 \subseteq M$ we have ${\Cal G}
\subseteq M$.  Apply the choice of ${\Cal G}$ to $\langle
v_{i,\eta(i)}:i < \kappa \rangle \in \dsize \prod_{i < \kappa} {\Cal
P}_i$, so for some $g \in {\Cal G} \subseteq M$ the
set $B =: \{i < \kappa:
v_{i,\eta(i)} = g(i)\}$ belongs to $J^+$.  Clearly $B \in M$ (as
$B \subseteq \kappa,{\Cal P}(\kappa) \in M$ and $|{\Cal P}(\kappa)|
\le 2^\kappa \le \theta^\kappa \le \lambda \subseteq M$) 
hence $\langle v_{i,\eta(i)}:i \in B \rangle
\in M$ hence $w = \{\alpha < \lambda$: for every $i \in B$ we
have $f_\alpha(i) \in v_{i,\eta(i)}\}$ belongs to $M$.  Now $|w| 
\le \dsize \prod_{i \in B} |v_{i,\eta(i)}| \le 
\theta^\kappa = \theta$ because $\alpha < \beta < \lambda
\Rightarrow f_\alpha \ne_J f_\beta \Rightarrow f_\alpha \restriction B
\ne f_\beta \restriction B$.  Lastly $w_\eta \subseteq w$
as $\alpha \in w_\eta \and i < \kappa \Rightarrow f_\alpha(i) \in
v_{i,\eta(i)}$, so we are done.  \hfill$\square_{\scite{g.4}}$\margincite{g.4}
\enddemo
\bigskip

\remark{Remark}  We could have used ${\Cal F} = \{\alpha <
\lambda:\{i:f_\alpha(i) \in v_{i,\eta(i)}\} \in J^+\}$.
\endremark
\bn
To make this section free of quoting the pcf theorem we use the
following definition.
\definition{\stag{g.5.Z} Definition/Observation}  1) For a set ${\frak
a}$ of regular cardinals and $\sigma = \text{ cf}(\sigma) \le \text{
cf}(\lambda)$ let

$$
\align
J^\sigma_\lambda[{\frak a}] = \{{\frak b} \subseteq {\frak a}:&\text{
there is a set } {\Cal F} \subseteq \pi {\frak b} \text{ of
cardinality } < \lambda \\
  &\text{ such that for every } g \in \pi {\frak b}
\text{ we can find } j < \sigma \text{ and} \\
  &\, f_i \in {\Cal F} \text{ for } i < j \text{ satisfying } \theta
  \in {\frak b} \Rightarrow (\exists i < j)(g(\theta) <
  f_i(\theta))\}.
\endalign
$$
\mn
2) Clearly $J^\sigma_\lambda[{\frak a}]$ is a $\sigma$-complete
ideal on ${\frak a}$ but possibly ${\frak a} \in
J^\sigma_\lambda[{\frak a}]$.
\enddefinition
\bigskip

\remark{Remark}  In fact, if Min$({\frak a}) > |{\frak
a}|,J^\sigma_\lambda[{\frak a}] = \{{\frak b} \subseteq {\frak a}$:
pcf$_{\sigma\text{-complete}}({\frak b}) \subseteq \lambda\} = \{b
\subseteq {\frak a}:{\frak b}$ is the union of $< \sigma$ members of
$J_\lambda[{\frak a}]\}$ can be proved but this is irrelevant here.  
\endremark
\bn
For completeness we recall and prove then we know
\proclaim{\stag{g.5} Claim}  $\lambda = \lambda^{[\sigma,<\theta]}$
when
\mr
\item "{$(a)$}"  $\lambda \ge 2^{< \theta} \ge \sigma = 
{ \text{\rm cf\/}}(\sigma) > \aleph_0$ and {\rm cf}$(\theta) \notin
[\sigma,\theta)$ 
\sn
\item "{$(b)$}"  for every set ${\frak a} \subseteq { \text{\rm Reg\/}} \cap
\lambda \backslash \theta$ of cardinality $< \theta$ we have
${\frak a} \in J^\sigma_\lambda[{\frak a}]$.
\endroster  
\endproclaim
\bigskip

\demo{Proof}  Let $\chi$ be large enough, choose $M \prec ({\Cal
H}(\chi),\in,<^*_\chi)$ of cardinality $\lambda$ where $<^*_\chi$ is
any well ordering of ${\Cal H}(\chi)$
such that $\lambda +1
\subseteq M$ and let ${\Cal P} = M \cap [\lambda]^{< \theta}$ and we
shall prove that ${\Cal P}$ exemplifies $\lambda =
\lambda^{[\sigma,<\theta]}$.

Clearly ${\Cal P} \subseteq [\lambda]^{< \theta}$ has cardinality
$\lambda$ so let $u \in [\lambda]^{< \theta}$ and as $2^{< \theta} \le
\lambda$ it is enough to show that
$u$ is included in a union of $< \sigma$ members of ${\Cal P}$ thus
finishing.

Let $f$ be a one-to-one function from $\kappa =: |u|$ onto $u$ so
$\kappa < \theta$.  By
induction on $n$ we choose $f_n,\bar v_n$ such that
\mr
\item "{$\circledast$}"  $(a) \quad f_n$ is a function from $\kappa$
to $\lambda +1$
\sn
\item "{${{}}$}"  $(b) \quad \bar v_n = \langle
v_{n,\varepsilon}:\varepsilon < \varepsilon_n \rangle$ is a partition
of $\kappa$ which satisfies \nl

\hskip35pt $\varepsilon_n < \sigma$ and $\kappa =
\cup\{v_{n,\varepsilon}:\varepsilon < \varepsilon_n\}$
\sn
\item "{${{}}$}"  $(c) \quad f_0(i) = \lambda$ for every $i < \kappa$
\sn
\item "{${{}}$}"  $(d) \quad f_{n+1}(i) \le
f_n(i)$ for $i < \kappa$
\sn
\item "{${{}}$}"  $(e) \quad f(i) \le f_n(i)$
and if $f(i) < f_n(i)$ then $f_{n+1}(i) < f_n(i)$
\sn
\item "{${{}}$}"  $(f) \quad f_n \restriction v_{n,\varepsilon} \in
M$.
\ermn
This is sufficient: $\{\text{Rang}(f_n \restriction
v_{n,\varepsilon}):n < \omega,\varepsilon < \varepsilon_n\}$ is a family of
$< \sigma$ sets (as $\sigma = \text{ cf}(\sigma) > \aleph_0,\sigma >
\varepsilon_n$) each belonging to ${\Cal P}$ (as $f_n \restriction
v_{n,\varepsilon} \in M$) and their union includes $u$ because for
every $i < \kappa,f_n(i) = f(i)$ for
every $n$ large enough (by clauses (d) + (e) of $\circledast$).

For $n=0,f_n$ is constantly $\lambda$.  So let $n=m+1,f_m$ be given,
let 

$$
u_{n,0} = \{i < \kappa:f_m(i) = f(i)\}
$$

$$
u_{n,1} = \{i < \kappa:f_m(i) > f(i)
\text{ and is a successor ordinal or just has cofinality } < \theta\}
$$

$$
u_{n,2} = \kappa \backslash u_{n,0} \backslash u_{n,1}.
$$
\mn
As $2^\kappa \le 2^{< \theta} \le \lambda$, clearly 
the partition $\langle u_{n,0},u_{n,1},u_{n,2}\rangle$ of
$\kappa$ belongs to $M$, so it is enough to choose $f_{n+1}
\restriction u_{n,\ell}$ separately for $\ell=0,1,2$.
\nl
Let $f_n \restriction u_{n,0} = f_m \restriction u_{n,0}$.
\nl
Let $\bar C = \langle C_\alpha:\alpha \le \lambda \rangle \in M$ be
such that $C_0 = \emptyset,C_{\alpha +1} = \{\alpha\},C_\delta$ is a
club of $\delta$ of order type cf$(\delta)$ for limit ordinal $\delta
\le \lambda$.  Let $f_n \restriction u_{n,1}$ be defined by
$f_n(i) = \text{ Min}(C_{f_m(i)} \backslash f(i))$.
For each $\varepsilon < \varepsilon_m$ the function $f_n \restriction
(u_{n,1} \cap v_{m,\varepsilon})$ belongs to $M$ hence $\langle
C_{f_m(i)}:i \in u_{n,1} \cap v_{m,\varepsilon} \rangle$ belongs to $M$,
but the product $\dsize \prod_{i \in u_{n,1} \cap v_{m,\varepsilon}}
C_{f_m(i)}$ has cardinality $\le \theta^\kappa \le 2^{< \theta} \le
\lambda$ if cf$(\theta) \ge \kappa$ and \wilog \, $\theta_1 =
\sup\{|C_{f_m(i)}|:i \in u_{n,1} \cap v_{m,\varepsilon}\} < \theta$
hence has cardinality $\le 2^{\theta_1 + \kappa} \le 2^{< \theta} \le
\lambda$ if cf$(\theta) < \sigma$.

Lastly, it is enough to define $f_n \restriction (v_{m,\varepsilon}
\cap u_{n,2})$ for each $\varepsilon < \varepsilon_m$.  Let
$\lambda_{n,i} = \text{ cf}(f_m(i))$, so $\langle
\lambda_{n,\zeta}:\zeta \in v_{m,\varepsilon} \cap u_{n,2} \rangle \in
M$ hence there is a sequence $\langle h_{n,\zeta}:\zeta \in v_{n,\varepsilon}
\cap u_{n,2} \rangle \in M$ where $h_{n,\zeta}$ is increasing
continuous function from $\lambda_{n,\zeta}$ onto some club of $f_m(i)$. 
\nl
Let ${\frak a} = \{\lambda_{n,\zeta}:\zeta \in v_{m,\varepsilon} 
\cap u_{n,2}\}$.
Applying assumption (b), Definition \scite{g.5.Z}(1) it is easy to finish.
\hfill$\square_{\scite{g.5}}$\margincite{g.5}
\enddemo
\bigskip

\proclaim{\stag{g.6} Claim}  There is $\bar \lambda = \langle
\lambda_i:i < \kappa \rangle$ and $\sigma$-complete ideal $J$ on
$\kappa$ such that $T_J(\bar \lambda)\ge \lambda$ and $i < \kappa
\Rightarrow 2^\kappa < \lambda_i < \lambda$ when
\mr
\item "{$\circledast$}"  $(a) \quad 2^\kappa < \lambda,\aleph_0 <
\sigma = { \text{\rm cf\/}}(\sigma) \le \kappa$
\sn
\item "{${{}}$}"  $(b) \quad {\frak a} \subseteq { \text{\rm Reg\/}} \cap
\lambda \backslash (2^\kappa)^+$ has cardinality $\le \kappa$ and
${\frak a} \notin J^\sigma_\lambda[{\frak a}]$.
\endroster  
\endproclaim
\bigskip

\demo{Proof}  Let $\bar \lambda = 
\langle \lambda_i:i < \kappa \rangle$ list ${\frak a}$ and let
$J = J^\sigma_\lambda[{\frak a}]$ and by induction on $\alpha <
\lambda$ a function $f_\alpha \in \pi({\frak a})$ such that $\beta <
\alpha \Rightarrow f_\beta <_J f_\alpha$.  Arriving to $\alpha$ for
every ${\frak b} \subseteq {\frak a}$ let ${\Cal F}^\alpha_{\frak b} =
\{f_\beta \restriction {\frak b}:\beta < \alpha\}$, so by the
definition of $J^\sigma_\lambda[{\frak a}]$, for every ${\frak b} \in
J^+ := {\Cal P}({\frak a}) \backslash J,{\Cal F}^\alpha_{\frak b}$ does not
witness ${\frak b} \in J^\sigma_\lambda[{\frak a}]$ hence there is
$g^\alpha_{\frak b} \in \pi({\frak b})$ witnessing it.  Let $f_\alpha \in
\pi({\frak a})$ be defined by $f_\alpha(\theta) = \sup\{g^\alpha_{\frak
b}(\theta):{\frak b} \in J^+$ and $\theta \in {\frak b}\}$.  Now
$f_\alpha \in \Pi {\frak a}$ as $\theta \in {\frak a} \Rightarrow
f_\alpha(\theta) < \theta$ which holds as $|J^+| \le 2^{|{\frak a}|}
\le 2^\kappa < \theta$.  Also if $\beta < \alpha$ and we let 
${\frak b}^\alpha_\beta =: \{\theta \in {\frak a}:f_\beta(\theta) \ge
f_\alpha(\theta)\}$ then ${\frak b}^\alpha_\beta \in J^+$ implies easy
contradiction to the choice of $g^\alpha_{{\frak b}^\alpha_\beta}$
(and $f_\alpha$).  So we can carry the induction and so $\langle
f_\alpha:\alpha < \lambda \rangle,f_\alpha \in \pi(\bar \lambda)$
where $f'_\alpha(i)= f_\alpha(\lambda_i)$ exemplify $T_J(\bar \lambda)
\ge \lambda$ as required.  \hfill$\square_{\scite{g.6}}$\margincite{g.6}
\enddemo
\bigskip

\remark{Remark}  This is the case 
Min$({\frak a}) > 2^{|{\frak a}|}$ from \cite[XIII]{Sh:b}).
\endremark
\bigskip

\proclaim{\stag{g.7} Claim}  If $\circledast$ from \scite{g.6} then 
there is $\bar \lambda' = \langle
\lambda'_i:i < \kappa \rangle$ such that
\mr
\item "{$(\alpha)$}"  $2^\kappa < \lambda'_i \le \lambda_i$
\sn
\item "{$(\beta)$}"  if $f \in \dsize \prod_{i < \kappa} \lambda_i$ then
$T_J(f) < \lambda$
\sn
\item "{$(\gamma)$}"  $T_J(\bar \lambda') = \lambda$ 
\ermn
\ub{where}
\mr
\item "{$\circledast$}"  $(a) \quad 2^\kappa < \lambda,\aleph_0 < \sigma = {
\text{\rm cf\/}}(\sigma) \le \kappa$
\sn
\item "{${{}}$}"  $(b) \quad 2^\kappa < \lambda_i < \lambda$
\sn
\item "{${{}}$}"  $(c) \quad J$ is a $\sigma$-complete ideal on $\kappa$
\sn
\item "{${{}}$}"  $(d) \quad T_J(\bar \lambda) \ge \lambda$.
\endroster  
\endproclaim
\bigskip

\demo{Proof}  Clearly $\{i:\lambda_i \le (2^\kappa)^{+n}\} \in J$ for $n
< \omega$ (as $((2^\kappa)^{+n})^\kappa = (2^\kappa)^{+n}$) so by
\scite{0.B}(1) \wilog \, $i < \kappa \Rightarrow \lambda_i > (2^\kappa)^{+2}$.

As $(\dsize \prod_{i < \kappa} \lambda^+_i,<_J)$ is well founded 
and there is $f \in \dsize \prod_{i < \kappa} (\lambda_i+1)$
satisfying $T_J(f) \ge \lambda$ (i.e. $\bar \lambda$ itself) clearly there is
$f \in \dsize \prod_{i < \kappa} (\lambda_i +1)$ for which $T_J(f)
\ge \lambda$ satisfying $g \in \dsize \prod_{i < \kappa} (\lambda_i+1),g <_J
f$ implies $T_D(g) < \lambda$.  Now easily $\{i< \kappa:f(i) \le
(2^\kappa)^{+2}\} \in J$, so \wilog \, $i < \kappa \Rightarrow
f(i) > (2^\kappa)^{+2}$.  Let
$\lambda'_i = |f(i)|$, hence $\bar \lambda'$ satisfies demands
$(\alpha) + (\beta)$ of the desired conclusion, i.e., and 
$T_J(\bar \lambda') = T_J(f) \ge \lambda$.
So assume toward contradiction that it fails clause $(\gamma)$, so by
the last sentence we have $T_J(\bar \lambda') > \lambda$ and we shall derive a
contradiction thus finishing.  So there is $\{f_\alpha:\alpha <
\lambda^+\} \subseteq \dsize \prod_{i < \kappa} \lambda'_i$ such that
$\alpha \ne \beta \Rightarrow f_\alpha \ne_J f_\beta$ and let
$u_\alpha =: \{\beta:f_\beta <_J f_\alpha\}$.  If for some $\alpha <
\lambda,|u_\alpha| \ge \lambda$ then $\{f_\beta:\beta \in u_\alpha\}$
exemplifies that $T_D(f_\alpha) \ge \lambda$ and clearly $f_\alpha <_J
\bar \lambda' \le f$, contradiction to the choice of $f$.  So $\alpha
< \lambda^+ \Rightarrow |u_\alpha| < \lambda$.  Hence there is $S
\subseteq \lambda$ such that $(\forall \alpha \ne \beta \in S)(\beta
\notin u_\alpha)$ contradicting \scite{0.K}.  \hfill$\square_{\scite{g.7}}$\margincite{g.7}
\enddemo
\bigskip

\proclaim{\stag{g.8} The revised GCH Theorem}   If 
$\theta$ is strong limit singular \ub{then} for 
every $\lambda \ge \theta$ for some $\sigma < \theta$ we 
have $\lambda = \lambda^{[\sigma,\theta]}$.
\endproclaim
\bigskip

\remark{Remark}  1) Hence for every $\lambda \ge \theta$ for some $n < \omega$
and $\kappa_\ell < \theta(\ell \le n),
\aleph_0 = \kappa_0 < \kappa_1 < \ldots < \kappa_n = \theta$, 
for each $\ell < n,2^{\kappa_\ell} \ge \kappa_{\ell+1}$ or
$\lambda = \lambda^{[\kappa'_\ell,<\kappa_{\ell +1}]}$ where
$\kappa'_\ell = (2^{\kappa_\ell})^+$. 
\nl
2) If $\sigma \in (\text{cf}(\theta),\theta)$ and $\lambda \ge \theta$
then $\lambda^{[\sigma,\theta]} = \lambda^{[\sigma,< \theta]} =
\Sigma\{\lambda^{[\sigma,\theta']}:\theta' \in [\sigma,\theta)\}$.
\endremark
\bigskip

\demo{Proof}  We prove by induction on $\lambda \ge \theta$.

Let $\sigma =: (\text{cf}(\theta))^+ < \theta$.
\enddemo
\bn
\ub{Case 0}:  $\lambda = \theta$.

Let ${\Cal P}$ be the family of bounded subsets of $\theta$, so
$|{\Cal P}| = \theta$ and every $u \in [\theta]^{< \theta}$ is the
union of cf$(\theta)$ members of ${\Cal P}$.
\bn
\ub{Case 1}:  For every ${\frak a} \subseteq \text{ Reg } \cap
\lambda^+ \backslash \theta$ of cardinality $< \theta$ we have
pcf$_{\sigma\text{-comp}}({\frak a}) \subseteq \lambda^+$.  By
\scite{g.5}, we have $\lambda^{[\sigma,<\theta]} = \lambda$ 
(recalling cf$(\theta) < \sigma$ and \scite{g.3}).
\bn
\ub{Case 2}:  Neither Case 0 nor Case 1.

As not Case 1, the assumption of Claim \scite{g.6} holds for some
$\kappa$ for which $\sigma \le \kappa < \theta$, hence its conclusion
holds for some $\bar \lambda = \langle \lambda_i:i < \kappa \rangle$
and $J$, i.e.,  we have $2^\kappa < \lambda_i < \lambda$ and 
$T_J(\bar \lambda) \ge \lambda$ where $J$ is a $\sigma$-complete ideal
on $\kappa$.
So the assumptions of Claim \scite{g.7} holds, hence its conclusion
so by \scite{g.7}
\mr
\item "{$\circledast$}"  $(i) \quad J$ is a $\sigma$-complete ideal on
$\kappa$
\sn
\item "{${{}}$}"  $(ii) \quad \bar \lambda' = \langle \lambda'_i:i <
\kappa \rangle$
\sn
\item "{${{}}$}"  $(iii) \quad 2^\kappa < \lambda'_i < \lambda$
\sn
\item "{${{}}$}"  $(iv) \quad T_J(\bar \lambda') = \lambda$
\sn
\item "{${{}}$}"  $(v) \quad T_J(f) < \lambda$ if $f \in \dsize
\prod_{i < \kappa} \lambda'_i$.
\ermn
We can find an increasing sequence $\langle
\theta_\varepsilon:\varepsilon < \text{ cf}(\theta)\rangle$ of regular
cardinals from $(\sigma,\theta)$ with limit $\theta$.  By the
induction hypothesis for each $i < \lambda$ there is $\varepsilon(i)$
such that $\lambda_i = \lambda^{[\theta_{\varepsilon(i)},< \theta]}_i 
\ge \theta$ or $\lambda_i \le \theta_{\varepsilon(i)}$.  For 
$\zeta < \text{ cf}(\theta)$ define $A_\zeta = \{i < \kappa:\lambda_i \ge
\theta,\varepsilon(i) = \zeta\}$ and $A_{\text{cf}(\theta)+\zeta} = \{i <
\kappa:\lambda_i < \theta$ and $\varepsilon(i) = \zeta\}$.  So $\langle
A_\varepsilon:\varepsilon < \text{ cf}(\theta) + \text{ cf}(\theta)
\rangle$ is a partition of $\kappa$ to $< \sigma$ sets hence by \scite{0.B}(4)
we know that

$$
T^1_J(\bar \lambda') = \text{ Min}\{T^1_{J \restriction
A_\varepsilon}(\bar \lambda' \restriction A_\varepsilon):\varepsilon <
\text{ cf}(\theta) + \text{ cf}(\theta) \text{ and } A_\varepsilon \in
J^+\}
$$
\mn
hence for some $\zeta < \text{ cf}(\theta) + \text{
cf}(\theta),T_J(\bar \lambda') =
T_{J \restriction A_\zeta}(\bar \lambda' \restriction A_\zeta)$ and
$A_\zeta \in J^+$, so by renaming without loss of generality 
$A_\zeta = \kappa$.  If $\zeta \ge \text{ cf}(\theta)$ as
$\kappa < \theta,\theta$ strong limit we get $T_J(\bar \lambda') \le
\dsize \prod_{i < \kappa} \lambda'_i < (\theta_\zeta)^\kappa <
\theta$, a contradiction so
$\zeta < \text{ cf}(\theta)$.  Now apply \scite{g.4} with $J,\bar
\lambda',\sigma,\kappa,(2^\kappa)^+,\theta_\zeta$ here standing for
$J,\bar \lambda,\sigma,\kappa,\partial,\theta$
and get the desired result. \nl
${{}}$  \hfill$\square_{\scite{g.8}}$\margincite{g.8} 
\bigskip

\remark{\stag{g.9} Concluding Remarks}  1) We can in \scite{g.4} assume
less.  Instead $\theta = \theta^\kappa$, it is enough (which follows
to \cite[\S3]{Sh:506}, see \scite{0.Q})
\mr
\item "{$\circledast$}"  for every $\lambda' < \lambda$ we can find
${\Cal F} \subseteq \dsize \prod_{i < \kappa} \lambda_i$ of
cardinality $\lambda'$ such that $f \ne g \in {\Cal F} \Rightarrow f
\ne_J g$.
\ermn
This is seemingly a gain, but in the induction the case $(\forall
{\frak a} \subseteq \text{ Reg } \cap \lambda^+ \backslash \theta)
(|{\frak a}| \le \kappa \Rightarrow 
\text{ pcf}_{\aleph_1\text{-comp}}({\frak a})
\subset \lambda^+)$ is problematic.
\endremark
\newpage

\head {\S2 The finitely many exceptions} \endhead  \resetall \sectno=2
 \spuriousreset
\bigskip

We can prove a relative of \scite{g.4} assuming $i < \kappa
\Rightarrow \lambda^{[\sigma,\kappa,\theta]}_i = \lambda_i$, replacing
``$\theta = \theta^\kappa + 2^\theta \le \lambda$ by $2^\kappa
\le \lambda$ and getting $\lambda^{[\sigma,\kappa,\theta]} =
\lambda$.  But so far it has no conclusion parallel to \scite{g.8}.
\bn
In full:
\proclaim{\stag{h.1} Claim}  Assume
\mr
\item "{$(a)$}"  $\aleph_0 < \sigma = { \text{\rm cf\/}}(\sigma) \le
\theta$ and $\mu \le \theta$
\sn
\item "{$(b)$}"  $J$ is a $\sigma$-complete ideal on $\kappa$
\sn
\item "{$(c)$}"  $\bar \lambda = \langle \lambda_i:i < \kappa \rangle$
\sn
\item "{$(d)$}"  $T_J(\bar \lambda) = \lambda$, moreover this is
exemplified by a $\mu^{++}$-free family 
\sn
\item "{$(e)$}"  $\lambda^{< \partial,\mu,\theta>}_i =
\lambda_i$ for $i < \kappa$ 
\sn
\item "{$(f)$}"   if $\partial_i < \partial$
for $i < \kappa$ then $\dsize \prod_{i < \kappa} \partial_i < \partial$
\sn
\item "{$(g)$}"  $\theta = \theta^\kappa$ and $2^\kappa \le \lambda$.
\ermn
\ub{Then} $\lambda^{<\partial,\mu,\theta>} = \lambda$.
\endproclaim
\bigskip

\demo{Proof}  Let $\bar f = \langle f_\alpha:\alpha < \lambda \rangle$
be $\theta^+$-free, $f_\alpha \in \dsize \prod_{i < \kappa} \lambda_i$
pairwise $J$-different, (i.e., $\alpha \ne \beta \Rightarrow \{i:f_\alpha(i) =
f_\beta(i)\} \in J$) exists by clause (d) of the assumption.

For each $i < \kappa$ let ${\Cal P}_i \subseteq [\lambda_i]^{\le \mu}$ 
be of cardinality $\lambda_i$ and witness
$\lambda^{<\sigma,\mu,\theta>}_i = \lambda_i$ that is:  every $u \in
[\lambda_i]^{\le \theta}$ is the union of $< \sigma$ members of 

$$
\text{set}_{\theta,\mu}({\Cal P}_i) =: \{v:v \in [\lambda_i]^{\le \mu}
\text{ and every } w \in [v]^{\le \mu} \text{ is included in some
member of } {\Cal P}_i\};
$$
\mn
such family exists by assumption (e).  
Let $M \prec ({\Cal H}(\chi),\in)$ be
of cardinality $\lambda$ such that $\lambda +1 \subseteq M$ and
$\langle \lambda_i:i < \kappa \rangle,\langle {\Cal P}_i:i <
\kappa \rangle,J,{\Cal P}(\kappa)$ belong to $M$.
\nl
Let ${\Cal P} = M \cap [\lambda]^{\le \mu}$.  We shall show that
${\Cal P}$ exemplifies the desired conclusion.  Now ${\Cal P}$ is a
family of $\le \|M\| = \lambda$ of subsets of $\lambda$ each of
cardinality $\le \mu$ hence it is enough to show
\mr
\item "{$(*)$}"  if $u \in [\lambda]^{\le \mu}$ then $u$ is
included in the union of $< \partial$ sets $v \in 
\text{ set}_{\theta,\mu}({\Cal P}_i)$.
\endroster
\enddemo
\bigskip

\demo{Proof of $(*)$}  Let $u_i = \{f_\alpha(i):\alpha \in u\}$; so
$u_i \in [\lambda_i]^{\le \theta}$, hence we can find $\langle
v_{i,j}:j < j_i \rangle$ such that $[v_{i,j}]^{\le \mu} 
\in \text{ set}_{\theta,\mu}({\Cal P}_i)$ and $u_i =
\cup\{u_{i,j}:j < j_i\}$ and $0 < j_i < \partial$.  For each $\eta \in
\dsize \prod_{i < \kappa} j_i$ let 

$$
w_\eta = \{\alpha \in u:i <\kappa \Rightarrow f_\alpha(i) \in
v_{i,\eta(i)}\}.
$$
\mn
Clearly $u = \cup\{w_\eta:\eta \in \dsize \prod_{i<\kappa} j_i\}$ as
for any $\alpha \in u$ for each $i < \kappa$ we can choose
$\varepsilon_i(\alpha) < j$ such that $f_\alpha(i) \in
v_{i,\varepsilon_i(\alpha)}$ and let $\eta_\alpha = \langle
\varepsilon_i(\alpha):i < \alpha \rangle$ clearly $\eta_\alpha \in
\dsize \prod_{i < \kappa} j_i$ and $\alpha \in w_{\eta_\alpha}$.  By
the assumption (f) as $i < \kappa \Rightarrow j_i < \partial$ clearly
$|\dsize \prod_{i < \kappa} j_i| < \sigma$ hence it is enough to
prove that $\eta \in \dsize \prod_{i < \kappa} j_i \Rightarrow w_\eta
\in \text{ set}_{\theta,\mu}({\Cal P})$.  So it is
enough to prove for $\eta\in \dsize \prod_{i < \kappa} j_i$ and $w \in
[w_\eta]^{\le \mu}$ that
\mr
\item "{$\circledast$}"  $w$ is included in some $w' \in M \cap 
[\lambda]^{\le \mu}$.
\ermn
\ub{Proof of $\circledast$}: As $i < \kappa \Rightarrow |{\Cal P}_i|
= \lambda_i$ and $T_J(\bar \lambda) = \lambda$ there is ${\Cal G}
\subseteq \{g:g \in \dsize \prod_{i < \kappa} {\Cal P}_i\}$ satisfying
$|{\Cal G}| = \lambda$ and $(\forall g \in \dsize \prod_{i < \kappa} {\Cal
P}_i)(\exists g' \in {\Cal G})(\{i:g(i) = g'(i)\} \in J^+\}$.  As
$\langle {\Cal P}_i:i < \kappa \rangle \in M$ \wilog \, ${\Cal G}
\in M$ and as $\lambda +1 \subseteq M$ we have ${\Cal G}
\subseteq M$.  For each $i < \kappa$ we have $A_i =
\{f_\alpha(i):\alpha \in w\}$ is a subset of some $A'_i,{\Cal P}_i$.
Apply the choice of ${\Cal G}$ to $\langle A'_i:i < \kappa 
\rangle \in \dsize \prod_{i < \kappa} {\Cal
P}_i$, so for some $g \in {\Cal G} \subseteq M$ the
set $B =: \{i:A'_i = g(i)\}$ belongs to $J^+$.  Clearly $B \in M$ (as
$2^\kappa \le \theta^\kappa \le \lambda$) 
hence $\langle A'_i:i \in B \rangle \in M$ hence $w' = 
\{\alpha < \lambda$: for some $Y \in J^+$ for every $i \in Y$ we
have $f_\alpha(i) \in A'_i\}$ belongs to $M$.  Now $|w| 
\le \mu^\kappa$; as $\alpha < \beta < \lambda
\Rightarrow f_\alpha \ne_J f_\beta$ but as $\bar f$ is $\mu^{++}$-free
we even have $|w| \le \mu$.  Lastly $w \subseteq {\Cal F}$
as $\alpha \in w \and i < \kappa \Rightarrow f(i) \in A_i$, so 
we are done.  \hfill$\square_{\scite{h.1}}$\margincite{h.1}
\enddemo
\bigskip

\proclaim{\stag{h.3} Claim}  
If $\theta > \sigma = { \text{\rm cf\/}}(\sigma) > \aleph_0$ and 
$\lambda > \theta_* = 2^{< \theta}$ \ub{then}
there is $\langle (\lambda_\eta,{\Cal D}_\eta,D_\eta,\kappa_\eta):
\eta \in {\Cal T} \rangle$ such that
\mr
\item "{$(a)$}"  ${\Cal T}$ is a subtree of ${}^{\omega >} \theta$
(i.e. $<> \in {\Cal T} \subseteq {}^{\omega >} \theta,{\Cal T}$ is 
closed under initial segments) with no $\omega$-branch,
let {\rm max}$_{\Cal T}$ be the set of maximal nodes of ${\Cal T}$
\sn
\item "{$(b)$}"  $\lambda_\eta$ is a cardinal $\in (2^{<
\theta},\lambda]$ and $\nu \triangleleft \eta \Rightarrow 
\lambda_\nu > \lambda_\eta$ and $\lambda_{<>} = \lambda$
\sn
\item "{$(c)$}"  $\kappa_\eta$ is a regular cardinal $ \in
[\sigma,\theta)$ if $\eta \in {\Cal T} \backslash { \text{\rm
max\/}}_{\Cal T}$ and $\kappa_\eta$ 
is zero or undefined if $\eta \in { \text{\rm max\/}}_{\Cal T}$ and
$\eta {}^\frown \langle \alpha \rangle \in {\Cal T} \Leftrightarrow
\alpha <\kappa_\eta$
\sn
\item "{$(d)$}"  if $\eta \in { \text{ \rm max\/}}_{\Cal T}$ \ub{then}
for no $\kappa < \theta$ and $\sigma$-complete filter ${\Cal D}$ on $\kappa$
and cardinals $\lambda_i \in (2^{< \theta},\lambda_\eta)$ for $i <
\kappa$ do we have $T_{\Cal D}(\langle \lambda_i:i < \kappa \rangle) 
\ge \lambda$
\sn
\item "{$(e)$}"  ${\Cal D}_\eta$ is a $\sigma$-complete filter on
$\kappa_\eta$ when $\eta \in {\Cal T} \backslash { \text{\rm max\/}}_{\Cal T}$
\sn
\item "{$(f)$}"  $T_{{\Cal D}_\eta}(\langle 
\lambda_{\eta {}^\frown <\alpha>}:\alpha < \kappa_\eta \rangle) =
\lambda_\eta$ if $\eta \in {\Cal T} \backslash \max_{\Cal T}$
\sn
\item "{$(g)$}"  if $f \in \dsize \prod_{\alpha < \kappa_\eta}
\lambda_{\eta {}^\frown <\alpha>}$ then $T_{{\Cal D}_\eta}(f) <
\lambda_\eta$ 
\sn
\item "{$(h)$}"  $D_\eta$ is the $\sigma$-complete filter on {\rm
max}$_{{\Cal T},\eta} = \{\nu \in {\text{\rm max\/}}_{\Cal T}:\eta
\trianglelefteq \nu\}$ such that
{\roster
\itemitem{ $(\alpha)$ }  if $\eta \in { \text{\rm max\/}}_{\Cal T},
D_\eta = \{\{\eta\}\}$
\sn
\itemitem{ $(\beta)$ }  if $\eta \in {\Cal T} \backslash
{ \text{\rm max\/}}_{\Cal T}$ then
$$
\align
D_\eta = \{A \subseteq { \text{\rm max\/}}_{{\Cal T},\eta}:&\text{ the
following set belongs to } {\Cal D}_\eta \\
  &\{\alpha < \kappa_\eta:A \cap { \text{\rm max\/}}_{{\Cal T},\eta
{}^\frown<\alpha>} \in D_{\eta \char 94<\alpha>}\}\}
\endalign
$$
\endroster}
\sn
\item "{$(i)$}"  $(\alpha) \quad$ if {\rm cf}$(\lambda) > \theta_*$ then
$\eta \in {\Cal T} \Rightarrow {\text{\rm cf\/}}(\lambda_\eta) >
\theta_*$
\sn
\item "{${{}}$}"  $(\beta) \quad$ we can replace $\theta_* = 2^{<
\theta}$ by any cardinal such that \nl

\hskip40pt {\rm cf}$(\theta_*) \ge \theta \wedge
(\forall \sigma < \sigma)(\forall \alpha < \theta_*)|\alpha|^\sigma < \theta_*$
\sn
\item "{${{}}$}"  $(\gamma)$ if $\lambda$ is 
regular then every $\lambda_\eta$ is regular [FILL proof].
\endroster
\endproclaim
\bigskip

\demo{Proof}
\mn
\ub{Case 1}:  Ignoring clause (i).

We prove this by induction $\lambda > 2^{<\theta}$.  
If $\lambda$ satisfies the requirement in
clause (d) let ${\Cal T} = \{<>\},\lambda_\eta = \lambda$ and
$\kappa_{<>},D_{<>}$ are trivial.  If $\lambda$ fails that demand use
claims \scite{g.6} +\scite{g.7} to find ${\Cal D},\kappa,\bar \lambda$
such that
\mr
\item "{$(*)$}"  $\kappa \in [\sigma,2^{< \theta}],{\Cal D}$ is a
$\sigma$-complete filter on $\kappa,\bar \lambda = \langle
\lambda_\alpha:\alpha < \kappa \rangle$ and $\lambda_\alpha \in
[\kappa,2^{< \theta}]$, a cardinal $T_{\Cal D}(\bar \lambda) =
\lambda$ but $f \in \dsize \prod_{\alpha < \kappa} \lambda_\alpha
\Rightarrow T_{\Cal D}(f) < \lambda$.
\ermn
Now for each $\alpha < \kappa$ we can use the induction hypothesis to
find $\langle (\lambda^\alpha_\eta:{\Cal D}^\alpha_\eta,{\Cal
D}^\alpha_\eta):\eta \in {\Cal T}_\alpha \rangle$ as required in the
claim for $\lambda_\alpha$.  Now we let:
\mr
\item "{$\circledast$}" $(a) \quad {\Cal T} = \{<>\} \cup \{\langle
\alpha \rangle {}^\frown \eta:\eta \in {\Cal T}_\alpha\}$
\sn
\item "{${{}}$}"  $(b) \quad \lambda^0_{<>} = \lambda,\kappa_{<>} =
\kappa$
\sn
\item "{${{}}$}"  $(c) \quad \lambda_{\langle \alpha \rangle {}^\frown
\eta} = \lambda^\alpha_\eta$ for $\alpha < \kappa,\eta \in {\Cal
T}_\alpha$
\sn
\item "{${{}}$}"  $(d) \quad {\Cal D}_{<>} = {\Cal D}$ 
\sn
\item "{${{}}$}"  $(e) \quad {\Cal D}_{\langle \alpha \rangle
{}^\frown \eta} = {\Cal D}^\alpha_\eta$ for $\alpha < \kappa,\eta \in
{\Cal T}_\alpha$
\sn
\item "{${{}}$}"  $(f) \quad D_{<>} = \{A:A \subseteq \max_{{\Cal
T},<>}$ and $\alpha < \kappa:\{\eta:\langle \alpha \rangle {}^\frown
\eta \in A \cap \max_{{\Cal T}_\alpha,<>}\}$ belongs to ${\Cal D}\}$
\sn
\item "{${{}}$}"  $(g) \quad D_{\langle \alpha \rangle {}^\frown \eta}
= \{\langle \alpha \rangle {}^\frown \nu:\nu \in B\}:B \in
D^\alpha_\eta\}$.
\ermn
Easily they are as required.
\mn
\ub{Case 2}:  Proving the claim with $(i)(\alpha)$, so dealing with
$\lambda$ satisfying cf$(\lambda) > \theta_x$.

If $\lambda$ satisfies the requirement in clause (a) we finish as
above.  Otherwise, we can find $\kappa \in
[\sigma,\theta),D,\bar\lambda$ such that
\mr
\item "{$(*)$}"  $\sigma$-complete filter ${\Cal D}$ on
$\kappa,\bar\lambda = \langle \lambda_\alpha:\alpha < \kappa \rangle$
and $\lambda_\alpha \in (2^{< \theta},\lambda)$ such that $\lambda =
T_D(\langle \lambda_\alpha:\alpha < \kappa \rangle)$ and $f \in \dsize
\prod_{\alpha < \kappa} \lambda_\alpha \Rightarrow T_D(f) < \lambda$.
\ermn
Let $B = \{\alpha:\text{cf}(\lambda_\alpha) > \theta_*\}$.  
If $B \in {\Cal D}^+$ and $T_{{\Cal D} \restriction B}(f \restriction B) 
< \lambda$ for $f \in \Pi \bar \lambda$ (hence $T_{{\Cal D}
\restriction B}(\bar \lambda \restriction B)(\lambda_{\eta {}^\frown
<\varepsilon>}) = \lambda)$, then we can use $\bar \lambda \restriction
B,{\Cal D} \restriction B$ (and renaming), so we are done, so
assume that this fails hence $\bar \lambda(\kappa \backslash B),
{\Cal D} \restriction (\kappa \backslash B)$ are as required in $(*)$, so by
renaming \wilog \, $B = \emptyset$.  For each $\alpha < \kappa$
let $\langle \lambda_{\alpha,\varepsilon}:\varepsilon <
\text{ cf}(\lambda_\alpha)\rangle$ be increasing continuous
with limit $\lambda_\alpha$ and let $\bar f =
\langle f_\alpha:\alpha < \lambda \rangle$ witness $T_{{\Cal
D}_\eta}(\bar \lambda) \ge \lambda$.  For each $\alpha < \lambda$ for some
$h_\alpha \in \dsize \prod_{\alpha < \kappa_\eta}$
cf$(\lambda_\alpha),f_\alpha < \langle \lambda_{\alpha,h(\alpha)}:
\alpha < \kappa \rangle$ so for some
$h,|\{\alpha:h_\alpha = h\}| = \lambda_\eta$, contradiction.

We finish as in case (1).
\mn
\ub{Case 3}: Restricting ourselves to a regular $\lambda$.
Similar to case to using \scite{h.5}(3) below.
\hfill$\square_{\scite{h.3}}$\margincite{h.3}
\enddemo
\bigskip

\definition{\stag{h.4} Definition}  1) We say that $\bar \lambda =
\langle \lambda_i:i < \kappa \rangle$ is a $D$-representation of
$\lambda$ \ub{when}: 
\mr
\item "{$(a)$}"  $D$ is a filter on $\kappa$
\sn
\item "{$(b)$}"  $T_D(\bar \lambda) = \lambda$
\sn
\item "{$(c)$}"  if $f \in \dsize \prod_{i < \kappa} \lambda_i$ then
$T_D(f) < \lambda$.
\ermn
2) We say that $\bar \lambda$ is an exact $D$-representation of
$\lambda$ when
\mr
\item "{$(a)$}"  $D$ is a filter on $\kappa$
\sn
\item "{$(b)$}"  $T_{D+A}(\bar \lambda) = \lambda$ for $A \in D^+$
\sn
\item "{$(c)$}"  if $f \in \dsize \prod_{i < \lambda} \lambda_i$ and
$A \in D^+$ then $T_{D+A}(f) < \lambda$.
\ermn
3) We say that the representation is true \ub{when}:
\mr
\item "{$(d)$}"  cf$(\lambda) = \text{ tcf}(\Pi \bar\lambda,<_D)$.
\ermn
4) We can replace the filter by the dual ideal.
\enddefinition
\bigskip

\definition{\stag{h.4A} Definition}  1) We say
$\langle(\lambda_\eta,{\Cal D}_\eta,D_\eta,\kappa_\eta):\eta \in {\Cal
T} \rangle$ is a representation if the conditions in Definition
\scite{h.3} holds.
\nl
2) We say it is an exact/true representation when each $\langle 
\lambda_{\eta {}^\frown<\alpha>}:\alpha < \kappa_n\rangle$ is exact/true a
$D_\eta$-representation of $\lambda_\eta$.
\enddefinition
\bigskip

\proclaim{\stag{h.5} Claim}  1) Assume
\mr
\item "{$\circledast$}"  $(a) \quad 
\bar\lambda^* = \langle \lambda_i:i < \kappa \rangle$ is a
$D_*$-representation of $\lambda$ 
\sn
\item "{${{}}$}"  $(b) \quad \bar \lambda^i = 
\langle \lambda_{i,j}:j < \kappa_i \rangle$ is a
$D_i$-representation of $\lambda_i$
\sn
\item "{${{}}$}"  $(c) \quad D$ is 
$\Sigma_{D_*} \langle D_i:i < \kappa \rangle$,
i.e., the filter on $u = \{(i,j):i < \kappa,j < \kappa_i\}$ \nl

\hskip20pt defined by
$D = \{A^* \subseteq S:\{i:\{j < \kappa_i:(i,j) \in A\} \in D_i\} \in D_*\}$
\sn
\item "{${{}}$}"  $(d) \quad$ {\rm cf}$(\lambda)$, \
{\rm cf}$(\lambda_i)$ are $> |u|$
and $\lambda,\lambda_i,\lambda_{i,j}$ are $> 2^{|u|}$.
\ermn
\ub{Then} $\bar\lambda = \langle \lambda_{i,j}:(i,j) \in u \rangle$ is a
$D$-representation of $\lambda$.
\nl
2) Similarly for exact representations, i.e., if in $\circledast(a),(b)$
we further assume that the representations are exact \ub{then} also
$\bar\lambda$ an exact $D$-representation of $\lambda$.
\nl
3) Similarly for true representations:
if $\lambda_i = { \text{\rm tcf\/}}(\dsize \prod_{j < \kappa} 
\lambda_{i,j},<_{D_i}),\lambda = { \text{\rm tcf\/}}(\dsize
\prod_{i < \kappa} \lambda_i,<_D)$ \ub{then} $\lambda = 
{ \text{\rm tcf\/}}(\dsize \prod_{(i,j)} \lambda_{i,j},<_D)$.  
Similarly for {\rm min-cf}, etc.
\nl
4) If $D$ is an $\aleph_1$-complete filter on $\kappa,\bar \lambda =
\langle \lambda_i:i < \kappa \rangle$ and $T_D(\bar \lambda) > \lambda
> 2^\kappa$ and $i < \kappa \Rightarrow \lambda_i > 2^\kappa$.
\ub{Then} we can find $\bar\lambda'$ such that $i < \kappa \Rightarrow
2^\kappa \le \lambda'_i < \lambda_i$ and $\bar \lambda'$ is a
$D$-representation of $\lambda$.  If we demand only $T_D(\bar\lambda)
\ge \lambda$ then we know only $\lambda'_i \le \lambda_i$.
\endproclaim
\bigskip

\demo{Proof}  1) 
\mr
\item "{$(*)_1$}"  $\lambda = T_D(\langle \lambda_{i,j}:(i,j) \in u
\rangle)$
\nl
[why?  Let ${\Cal G}^i = \{g^i_\alpha:\alpha < \lambda_i\}$ witness
that $T_{D_i}(\bar \lambda^i) = \lambda_i$ and let
${\Cal G}^*= \{g^*_\alpha:\alpha < \lambda\}$ witness 
that $T_{D_*}(\bar \lambda^*) = \lambda$.  We now define
${\Cal G} = \{g_\alpha:\alpha < \lambda\}$ where $g_\alpha \in \dsize
\prod_{(i,j) \in u} \lambda_{i,j}$ is defined by $g_\alpha((i,j)) =
g^i_{g^*_\alpha(i)}(j)$ and we can easily check that $\alpha < \beta < \lambda
\Rightarrow g_\alpha \ne g_\beta$ mod $D$.  Now if $g \in \dsize
\prod_{(i,j) \in u} \lambda_{i,j}$ then for each $i$ the function
(i.e. sequence) $\langle g((i,j)):j < \kappa_{i,j} 
\rangle$ belongs to $\dsize \prod_{j
< \kappa_i} \lambda_{i,j}$ so for some $\gamma_i < \lambda_i$ we have
$\{j:g((i,j)) = g^i_{\gamma_i}(j)\} \in D^+_i$.
Similarly for some $\beta < \lambda,\{i < \kappa:\gamma_i =
g^*_\beta(i)\} \in D^+_*$.  Easily $\{(i,j) \in u:g_\beta(i,j) = g(i,j)\} \in
D^+$, so ${\Cal G}$ witness that $T_D(\bar \lambda) = \lambda$ is as required.]
\sn
\item "{$(*)_2$}"  if $g \in \Pi\{\lambda_{i,j}:(i,j) \in u\}$ then
$T_D(g) < \lambda$
\nl
[Why?  Without loss of generality $g((i,j)) > 0$ for every $(i,j)
\in u$.  For each $i < \kappa$, let $g_i \in \dsize \prod_{j<\kappa_i}
\lambda_{i,j}$ be defined by $g_i(j) = g((i,j))$.  So $g_i \in 
\dsize \prod_{j<\kappa_i} \lambda_{i,j}$ hence $\mu_i =: T_{D_i}(g_i)
< \lambda_i$ hence there is a sequence $\langle h^i_\alpha:i < \mu_i
\rangle$ such that $h^i_\alpha \in \dsize \prod_{j<\kappa_i} g_i(j)$
and $(\forall h \in \dsize \prod_{j<\kappa_i} g_i(j))(\exists \alpha <
\mu_i)(\neg(h \ne_{D_i} h'_\alpha))$.  Clearly $\bar \mu = \langle
\mu_i:i < \kappa \rangle \in \dsize \prod_{j<\kappa} \lambda_i$ hence
$\mu_* =: T_{D_*}(\bar \mu) < \lambda$ and let $\langle
g^{**}_\alpha:\alpha < \mu \rangle$ exemplifies this.  We now define
$f^{**}_\alpha \in \dsize \prod_{(i,j) \in u} g((i,j))$ by
$f^{**}_\alpha((i,j)) = h^i_{g^{**}_\alpha(i)}(j)$ and it suffices to
show that $T_D(g) \le \mu (< \lambda)$ is exemplified by
$\{f^{**}_\alpha:\alpha < \mu\}$ which is proved as in $(*)_1$, the
second half of the proof.]
\ermn
So we are done. 
\nl
2) Similarly.
\nl
3) By \cite[I]{Sh:e}.  
\nl
4) Easy (and proved above).   \hfill$\square_{\scite{h.5}}$\margincite{h.5}
\enddemo
\bigskip

\remark{\stag{h.6} Remark}  1) So if $D$ is defined from $D^*,\langle
D_i:i < \kappa \rangle$, as in \scite{h.5} and $\bar \lambda =
\langle \lambda_{i,j}:(i,j) \in u \rangle,\lambda_i = T_{D_i}(\langle
\lambda_{i,j}:j < \kappa_i \rangle),\lambda = T_{D_*}(\langle \lambda_i:i
< \kappa \rangle)$ then $\lambda = T_D(\bar \lambda)$.
\endremark
\bn
We may wonder \nl
\margintag{h.7}\ub{\stag{h.7} Question}:   Can we add in Claim \scite{h.3}:
\nl
1) If $\lambda$ is regular can we add: each $\lambda_\eta$ is
regular. 
\nl
2) Can add the case of nice filters.
\bn
See below.
\proclaim{\stag{h.6.1} Claim}  If $\bar\lambda = \langle \lambda_i:i <
\kappa \rangle$ is a $J$-representation of $\lambda,\lambda = 
{ \text{\rm cf\/}}(\lambda) > 2^\kappa$ and $J$ is an
$\aleph_1$-complete ideal on $\kappa$ \ub{then} for some $J' \in J^+$,
the sequence $\bar \lambda$ is a $J'$-representation of $\lambda$ and
$\dsize \prod_{i < \kappa} \lambda_i/J'$ has true cofinality $\lambda$
(hence $\{i:\lambda_i$ singular$\} \in J$).
\endproclaim
\bigskip

\demo{Proof}  By the pcf theorem there is $u^* \subseteq \kappa$ such
that $\lambda \notin \text{ pcf}\{\lambda_i:i \in \kappa(u^*)\}$ and
$\lambda \ge \text{ cf}(\dsize \prod_{i \in u^*} \lambda_i)$.  If
$\bar\lambda$ is a $(J+u^*)$-representation of $\lambda$ then $\lambda
= T_{J+u^*}(\bar \lambda)$ but this implies that for some $u \in (J+
u^*)^+$ we have $\dsize \prod_{i \in u} \lambda_i/(J +u^*)
\restriction u$ has true cofinality $\lambda$ by \cite[1.1]{Sh:589},
contradiction.  Hence $\bar \lambda \restriction u^*$ is a $(J
\restriction u^*)$-representation of $\lambda$ (by \sciteu{xxX}), so
\wilog \, $u^* = \kappa$, so $\lambda = \text{ max pcf}\{\lambda_i:i <
\kappa\}$.  Let $J_1 = \{u \subseteq \kappa:u \in J$ or $u \notin J$
and ${\Cal P}(u) \cap J_2 \subseteq J_1$ where 
$J_\chi = \{u \subseteq \kappa:u \in J$ or for some $v \in J$ we have
$\lambda > \text{ max pcf}\{\lambda_i:i \in u \backslash v\}$.
Clearly $J_1,J_2$ are ideals on $\kappa$ extending $J$ adn by the pcf
theorem we hvae $J_1 \cap J_2 = J$.  So by \sciteu{xxX} for some $\ell
\in \{1,2\},\bar\lambda$ is a $J_\ell$-representation of $\lambda$.
\mn
\ub{Case 1}:  $\ell=1$.

So $\lambda = T_{J_1}(\bar \lambda)$ hence by \cite[1.1]{Sh:589} for some $v
\in (J_1)^+,\dsize \prod_{i \in v} \lambda_i(J_1 \restriction v)$ has
true cofinality $\lambda$.  So if $u \in J_2 \backslash J$ then for some
$u' \subseteq u,\lambda > \text{ max pcf}(\{\lambda_i:i \in u \backslash u'\})$
and $J_2 \restriction u = J \restriction u$ hence $v \cap u \in J$.
But this means $v \cap u \in J$ for every $u \in J_2 \backslash J$
(and, of course, for $u \in J$) hence $v \in J_1$, contradiction.
\mn
\ub{Case 2}:  $\ell=2$.

By the pcf theorem $\dsize \prod_{i < \kappa} \lambda_i/J_2$ has true
cofinality $\lambda$.  \hfill$\square_{\scite{h.6.1}}$\margincite{h.6.1}
\enddemo
\bigskip

\demo{\stag{h.7.2} Conclusion}  In \scite{h.3} we can add:
\mr
\item "{$(j)$}"  if $\lambda$ is regular then every $\lambda_\eta$ is
regular and for $\eta \in {\Cal T} \backslash \max_{\Cal T}$ we have
$\lambda_\eta = \text{ tcf}(\dsize \prod_{\alpha < \kappa_\eta}
\lambda_{\eta {}^\frown \langle \alpha \rangle}/{\Cal D}_\eta)$.
\endroster
\enddemo
\bn
Now \scite{h.7}(2) (and also \scite{h.7}(1)) are answered by:
\nl
\ub{Saharon}:  Change \scite{h.7.2} to cover singleton $\lambda$ as
needed in \scite{md.1}!
\bigskip

\demo{\stag{0.23.1} Observation}  1) Assume that
\mr
\item "{$(a)$}"  $J_1,J_2$ are ideals on $\kappa$ with intersection
$J_0$
\sn
\item "{$(b)$}"  $f_i \in {}^\kappa(\text{Ord} \backslash \omega)$.
\ermn
Then $T_J(f) = \text{ Min}\{T_{J_1}(f),T_{J_2}(f)\}$.
\nl
2) If (a) above holds and $\bar \lambda$ is a $J$-representation of
$\lambda$, \ub{then} for some $\ell \in \{1,2\},\bar\lambda$ is a
$J_\ell$-representaiton of $\lambda$.
\nl
3) See $D_1$ in the proof of \scite{md.1}.
\enddemo
\bigskip

\proclaim{\stag{h.8} Claim}  Assume \footnote{without this assumption
much more follows, see \cite[V]{Sh:g}.}
that the pair $(\bold K[S],\bold V)$ fails the
covering lemma for every $S \subseteq \beth_2(\kappa)$ (or less).
\ub{Then} in \scite{h.3} we can add:
\nl
1) If ${\frak a} \subseteq { \text{\rm Reg\/}} \cap \lambda \backslash
(2^{< \theta})^+$ and $|{\frak a}| < \theta$ and
{\rm pcf}$_{\sigma\text{-comp}}({\frak a}) \nsubseteq \lambda,\lambda
> 2^{< \theta}$ \ub{then} for some $\kappa = {
\text{\rm cf\/}}(\kappa) \in [\sigma,\theta)$ and $\kappa$-complete
ideal $J$ on $\kappa$ and $\bar \lambda = \langle \lambda_i:i < \kappa
\rangle$ we have:
\mr
\item "{$(a)$}"  {\rm cf}$(\lambda) > 2^{< \theta} \Rightarrow$ 
{\rm cf}$(\lambda_i) > 2^{< \theta}$
\sn
\item "{$(b)$}"  $\langle \lambda_i:i < \kappa \rangle$ is an exact
true $J$-representation of $\lambda$
\sn
\item "{$(c)$}"  if $\lambda$ is regular then every $\lambda_i$ is
regular.
\ermn
2) For any normal filter $D$ on $\kappa$ we can further demand in part (1) 
that for some
$\iota:\kappa \rightarrow \kappa,(J,\iota)$ is nice and $A \in D
\Rightarrow \iota^{-1}(\kappa \backslash A) \in J$.   
\nl
2A) If $\sigma \ge \partial = { \text{\rm cf\/}}(\partial) > \aleph_0$
and $D$ is a normal filter on $\partial$ we can add in part (1) that
$(J,\iota)$ is nice and $A \in D = \iota^{-1}(\kappa \backslash A)$.
Similarly for normal filters on $[\sigma_1]^{< \partial}$.
\nl
3) So in \scite{h.3}, we can strengthen clauses (f),(g) to
\mr
\item "{$(f)^+$}"  if $A \in {\Cal D}^+_\eta,\eta \in {\Cal T}
\backslash { \text{\rm max\/}}_{\Cal T}$ and $f \in \dsize
\prod_{\alpha < \kappa_\eta} \lambda_{\eta {}^\frown <\alpha>}$
then $T_{{\Cal D}_\eta +A}(\langle \lambda_{\eta {}^\frown<\alpha>}:
\alpha < \kappa_\eta) = \lambda_\eta > T_{{\Cal D}_\eta +A}(f)$ 
(hence the parallel result for $D_\eta$)
\sn
\item "{$(g)^+$}"  if $\eta \in {\Cal T} \backslash
{ \text{\rm max\/}}_{\Cal T},A \in {\Cal D}^+_\eta$ and $f \in \dsize
\prod_{\alpha < \kappa_\eta} \lambda_{\eta {}^\frown<\alpha>}$
\ub{then} $T_{{\Cal D}_\eta +A}(f) < \lambda_\eta$ (hence 
the parallel result for $D_\eta$)
\sn
\item "{$(h)^+$}"  for each $\eta \in {\Cal T} \backslash \max_I$ for some
$\iota_\eta:\kappa_\eta \rightarrow \kappa_\eta$ the pair
$({\Cal D}_\eta,\iota_\eta)$ is nice
\sn
\item "{$(i)$}"  $(\gamma) \quad$ if $\lambda$ is regular \ub{then} every 
$\lambda_\eta$ is regular.
\endroster
\endproclaim
\bigskip

\demo{Proof}  By \cite[\S3]{Sh:420}, better \cite{Sh:835}, very close
to \cite{Sh:386}. \nl
1) There are $\kappa$ a $\kappa$-complete filter on $\kappa$ and
$\lambda_i < \lambda$ such that $T_D(\langle \lambda_i:i < \kappa
\rangle) \ge \lambda$.  By the results quoted above \wilog \, $D$ is 
a normal filter on $\kappa \times \kappa$ 
with $\iota(\alpha,\beta) = \alpha$.  Now we can choose
$(D,\bar \lambda)$ such that $D$ is a nice filter on $\kappa \times
\kappa,T_D(\bar \lambda) \ge \lambda$ and rk$^3_D(\bar \lambda)$ is
minimal.  As $D_1 \subseteq D_2 \Rightarrow T_{D_1}(\bar \lambda) \le
T_{D_2}(\bar \lambda)$ \wilog \, rk$^3_D(\bar \lambda) = 
\text{ rk}^2_D(\bar \lambda)$ and so $A \in D^+ \Rightarrow \text{
rk}^3_{D+A}(\lambda) = \text{ rk}^2_{D+A}(\lambda) = \text{
rk}^3_D(\bar \lambda)$ and $T_{D+A}(\bar \lambda) \ge T_D(\bar
\lambda)$.  If $T_{D+A}(\bar \lambda)$ then for some $f \in \pi
\bar\lambda,T_{D+A}(f) \ge \lambda$, let $\bar \lambda' = \langle
f(i):i < \kappa \rangle$, so $\bar \lambda' <_D \bar \lambda$ hence
rk$^3_{D+A}(\bar \lambda') < \text{ rk}^3_{D+A}(\bar \lambda)$ and we
get a contradiction).
\nl
2), 2A), 3)  Left to the reader.  \hfill$\square_{\scite{h.8}}$\margincite{h.8}
\enddemo
\bigskip

\definition{\stag{h.9} Definition}  Assume $\sigma < \theta < \lambda$.
\nl
1) Let ${\frak d}_0(\lambda) = {\frak d}^0_{\sigma,\theta}(\lambda) =
\{\kappa:\kappa = \text{ cf}(\kappa) \in
\text{ Reg } \cap \theta \backslash \sigma$ such that 
we cannot find $\langle (\lambda_\eta,{\Cal D}_\eta,
D_\eta,\kappa_\eta):\eta \in {\Cal T} \rangle$ satisfying
$[\eta \in {\Cal T} \Rightarrow |{\Cal T}| < \lambda_\eta]$ and
$\kappa \notin \{\kappa_\eta:\eta \in {\Cal T}\}\}$ (so finite!) 
If $\sigma = \aleph_1$ we may omit it.  If $\sigma = \aleph_1,\theta =
\lambda$ we may omit both.
\nl
2) Let ${\frak d}_1(\lambda) = {\frak d}^1_{\sigma,\theta} = 
\{\kappa:\kappa = \text{ cf}(\kappa) < \lambda$
and for arbitrarily large $\alpha < \lambda$ we have $\kappa \in
{\frak d}_0(|\alpha|)\}$; note that if cf$(\lambda) > \aleph_0$ we can
deduce the finiteness of ${\frak d}_1(\lambda)$ from the finiteness
of ${\frak d}_0(\lambda)$.
\nl
3) Let ${\frak d}'_\ell(\lambda) = {\frak d}_\ell(\lambda) \cup
\{\aleph_0\}$ for $\ell=0,1$; similarly ${\frak
d}'_{\ell,\theta,\sigma}$ or ${\frak d}'_{\ell,\theta}$.
\enddefinition
\newpage

\head {\S3 The main results $(\text{Pr}_\ell,\text{Ps}_\ell)$} \endhead  \resetall \sectno=3
 \spuriousreset
\bigskip

\proclaim{\stag{md.1} Lemma}  Assume that $\mu > \aleph_0$ is
 strong limit singular, $\lambda \ge { \text{\rm cf\/}}(\lambda) >
\mu$ and $h:{\text{\rm cf\/}}(\lambda) \rightarrow \lambda$ is 
increasing continuous with $\lambda = \cup\{h(\alpha):\alpha < \lambda\}$
 and satisfies {\rm cf}$(\lambda) < 
\lambda \Rightarrow { \text{\rm Rang\/}}(h) \subseteq$ {\rm Card} and
$\lambda  = { \text{\rm cf\/}}(\lambda) \Rightarrow h(\alpha) =
 \alpha$.  \ub{Then} for some $\kappa < \mu$ 
and finite ${\frak d} \subseteq { \text{\rm Reg\/}} \cap \mu$ there is
 $\bar{\Cal P}$ such that
\mr
\item "{$(*)_{\lambda,\bar{\Cal P}}$}"  $\bar {\Cal P} = 
\langle {\Cal P}_\alpha:
\alpha < { \text{\rm cf\/}}(\lambda) \rangle,
{\Cal P}_\alpha \subseteq [h(\alpha)]^{< \mu},|{\Cal P}_\alpha| < \lambda,
{\Cal P}_\alpha$ increasing,
\sn
\item "{$(*)^{{\frak d},\kappa}_{\lambda,\bar{\Cal P}}$}"  for 
every $A$ satisfying $A \subseteq { \text{\rm cf\/}}(\lambda)$ or
(more generally) $A \subseteq \lambda \and (\forall \alpha \in
A)[{\text{\rm Min\/}}(A \backslash (\alpha + 1) < { \text{\rm
Min\/}}({\text{\rm Rang\/}}(h)) \backslash (\alpha+1))]$ and satisfying
$|A| < \mu$ \ub{there is} $\bold c:[A]^2 \rightarrow
\kappa$ such that: \nl
if $B \subseteq A$ has no last element
and $\bold c \restriction [B]^2$ is constant and
$\delta = \cup\{{\text{\rm Min\/}}(\alpha +1):\gamma <
h(\alpha)\}:\gamma \in B\}$ has cofinality $\in 
({\text{\rm Reg\/}} \, \cap \mu \backslash {\frak d})$ so $B
\subseteq h(\delta)$ then $B \in {\Cal P}_\delta$.
\endroster
\endproclaim
\bigskip

\remark{\stag{md.1R} Remark}  1) The proof of \scite{md.1} is 
simpler if $\lambda$ is regular.
\nl
2) The conclusion of \scite{md.1} means that for $\lambda > \mu$, for
all but finitely many $\kappa = \text{ cf}(\kappa) < \mu$,
Pr$(\lambda$,cf$(\lambda),\kappa)$ holds (see Definition \scite{md.6}(b)).
\endremark
\bn
Similarly
\proclaim{\stag{md.1S} Claim}  In fact in \scite{md.1} we can 
choose ${\frak d} = {\frak d}'_{0,\mu}(\lambda)$.
\endproclaim
\bigskip

\demo{Proof}  Without loss of generality cf$(\mu) = \aleph_0$ 
(this is no loss by Fodor's Lemma; otherwise we may use 
$\mu > \text{ cf}(\mu)$ or replace
$\aleph_1$ by (cf$(\mu))^+$.)  We first prove the desired conclusions
for cardinal $\lambda$ such that
\mr
\item "{$\boxtimes_\lambda$}"   ${\frak a} \subseteq \lambda \cap 
\text{ Reg} \backslash \mu \and |{\frak a}| < \mu 
\Rightarrow \text{ pcf}_{\aleph_1\text{-comp}}({\frak a}) 
\subseteq \lambda$.
\ermn
Let $\bar M = \langle M_\alpha:
\alpha < \text{ cf}(\lambda)\rangle$ be such that
\mr
\item "{$\circledast_1$}"   $(a) \quad M_\alpha \prec ({\Cal
H}(\chi),\in)$ is increasing continuous
\sn
\item "{${{}}$}"   $(b) \quad \lambda \in M_\alpha,\|M_\alpha\| <
\lambda,h(\alpha) \subseteq M_\alpha$
\sn
\item "{${{}}$}"   $(c) \quad \langle M_\alpha:\alpha \le \beta
\rangle \in M_{\beta +1}$
\sn
\item "{${{}}$}"   $(d) \quad (\alpha) \quad$ if $\lambda$ is regular then 
$M_\alpha \cap \lambda \in \lambda$
\sn
\item "{${{}}$}"   $\quad \quad\,\, (\beta) \quad$ if $\lambda$ is singular and
$\alpha < \text{ cf}(\lambda)$ then $\mu_\alpha +1
\subseteq M_{\alpha + 1}$ \nl

\hskip43pt where $\mu_\alpha = \text{ Min}\{\mu$: if ${\frak a} \subseteq
\mu \cap \text{ Reg} \backslash \mu$ and $|{\frak a}| < \mu$ \nl

\hskip43pt then pcf$_{\aleph_1-\text{com}}({\frak a}) 
\subseteq \mu$ and $\mu \ge \|M_\alpha\|\}$.
\ermn
We let ${\Cal P}_\alpha =: M_{\alpha +1} \cap [h(\alpha)]^{< \mu}$ and
${\frak d} = \{\aleph_0\}$ and $\kappa = \aleph_0$ and
will show that $\langle {\Cal P}_\alpha:\alpha < \text{ cf}
(\lambda)\rangle,{\frak d}$ are as required.  
Now $(*)_{\lambda,\bar{\Cal P}}$ of the claim holds trivially.  To prove
$(*)^{{\frak d},\kappa}_{\lambda,\bar{\Cal P}}$ let 
$A \subseteq \lambda$, otp$(A) \le \mu$ be as there and let
$\{\alpha_\varepsilon:\varepsilon < \varepsilon(*)\}$ list $A$ in
increasing order.  Hence there is $\langle \beta_\varepsilon:\varepsilon <
\varepsilon(*) \rangle$ increasing continuous such that $\beta_\varepsilon <
\text{ cf}(\lambda),h(\beta_\varepsilon) \le \alpha_\varepsilon <
h(\beta_{\varepsilon +1})$.  
By the assumption (and \scite{g.5}, i.e., \cite[II,5.4]{Sh:g}) if
$\lambda$ is regular then for 
each $\varepsilon < \varepsilon(*)$ there is a set 
${\Cal P}^\varepsilon \subseteq [h(\beta_\varepsilon)]^{< \mu}$ 
of cardinality $< \lambda$ such that every $a \in
[h(\beta_\varepsilon)]^{< \mu}$ is equal to the union of $\le \kappa$ of 
them (by the choice of $\kappa$) hence \wilog \, 
${\Cal P}^\varepsilon \in M_{\beta_\varepsilon +1}$ 
hence ${\Cal P}^\varepsilon \subseteq M_{\beta_\varepsilon +1} \cap
[h(\beta_\varepsilon)]^{< \mu} = {\Cal P}_{\beta_\varepsilon}$.  
If $\lambda$ is singular using clause $(d)(\beta)$ we get the same conclusion.
So there is a sequence $\langle A_{\varepsilon,i}:i < \kappa \rangle$ such
that $A_{\varepsilon,i} \in {\Cal P}_{\beta_\varepsilon},A \cap
\alpha_\varepsilon = \cup\{A_{\varepsilon,i}:i < \kappa\}$.  We
defined $\bold c:[A]^2 \rightarrow \kappa$ by: for $\varepsilon <
\zeta < \varepsilon(*),\bold c\{\alpha_\varepsilon,\alpha_\zeta\} =
\text{ Min}\{i:\alpha_\varepsilon \in A_{\zeta,i}\}$.  So assume $B
\subseteq A$ and $\bold c \restriction [B]^2$ is constantly 
$i < \kappa$ and $\delta =
\sup(B)$ has cofinality $\theta \in \text{ Reg } \cap \mu \backslash
{\frak d}$.  Clearly
$\alpha_\varepsilon \in B \Rightarrow \alpha_\varepsilon \cap B
\subseteq \{\alpha_\zeta:\zeta < \varepsilon,\bold
c\{\alpha_\zeta,\alpha_\varepsilon\} =i\} \subseteq A_{\varepsilon,i}
\in {\Cal P}_{\beta_\varepsilon}$.  But ${\Cal P}_\alpha = M_{\alpha
+1} \cap [h(\alpha)]^{< \mu}$ is closed under subsets hence
$\alpha_\varepsilon \in B \Rightarrow \alpha_\varepsilon \cap B \in
{\Cal P}_{\beta_\varepsilon}$.

Now in $M_{\delta +1}$ we can define a tree ${\Cal T}$; it has
otp$(B)$ levels; 
\sn
\block
the $i$-level is $\{a \in M_\delta:a \subseteq
\delta$ and otp$(a) = i\}$
\endblock
\sn
the order is $\triangleleft$, being initial segments.

So by the assumptions (and \cite[\S2]{Sh:589}), as 
$\aleph_1 \le \text{ cf}(\delta) <
\mu$ the number of $\delta$-branches of ${\Cal T}$ is $< \lambda$,
so as ${\Cal T} \in M_{\delta +1}$, every $\delta$-branch of 
${\Cal T}$ belongs to $M_{\delta +1}$, hence $B \in M_{\delta +1}$,
which implies that $B \in {\Cal P}_\delta$, as required.

Now we prove the statement in general.
\bn
We prove this by induction on $\lambda$.  For $\lambda = \mu^+$ this 
is trivial by the first part of the proof.  So
assume $\lambda > \mu^+$ and the conclusion fails,
but the first part does not apply.

In particular for some ${\frak a} \subseteq \text{ Reg } \cap
\lambda^+,|{\frak a}| < \mu$ and
pcf$_{\aleph_1\text{-comp}}({\frak a}) \nsubseteq \lambda$.
Hence by \scite{g.6} + \scite{g.7} + \scite{h.3} + \scite{h.6.1}
(alternatively + \scite{h.8} as in the proof
of \scite{h.3}; if the hypothesis of \scite{h.8} fails that as in
\cite{Sh:386} Jensen Dodd covering lemma gives much more), for some 
$\kappa = \text{ cf}(\kappa) \in [\aleph_1,\mu)$, 
\mr
\item "{$(*)_1$}"  there is a sequence $\langle \lambda_i:i < \kappa
\rangle$ and an $\aleph_1$-complete filter $D$ on $\kappa$ such that
{\roster 
\itemitem{ $(a)$ }  $T_D(\dsize \prod_{i < \kappa} \lambda_i) = \lambda$
\sn
\itemitem{ $(b)$ }  $\text{cf}(\lambda_i) > \mu$
\sn
\itemitem{ $(c)$ }  if $\lambda'_i < \lambda_i$ for $i < \kappa$, 
\ub{then} $T_D(\langle \lambda'_i:i < \kappa \rangle) < \lambda$
\sn
\itemitem{ $(d)$ }  tcf$(\dsize \prod_{i < \kappa} \lambda_i,<_D) = 
\text{ cf}(\lambda)$. 
\endroster}
\ermn
Clearly we can find $\langle h_i:i < \kappa \rangle$ such that
\mr
\item "{$(*)_2$}"  $h_i$ is an increasing
continuous function from cf$(\lambda_i)$ to $\lambda_i$.
\ermn
Let

$$
\align
D_1 = \{A:&A \in D \text{ or } A \notin D,A \in D^+ \text{ and} \\
  &T_{D+(\kappa \backslash A)}(\bar \lambda') \ge \lambda \text{ for
  some } \bar \lambda' \in \dsize \prod_{i < \kappa} \lambda_i\}.
\endalign
$$
\mn
Clearly
\mr
\item "{$(*)_3$}"  $D_1$ is a $\kappa$-complete filter on $\kappa$
extending $D$ and we can replace $D$ by $D+A$ whenever $A \in D^+_1$.
\ermn
By the induction hypothesis applied to $\lambda_i$, as $\lambda_i >
\mu$ there is a pair $(\kappa_i,{\frak d}_i)$ as in the conclusion.
Without loss of generality $\kappa^\kappa_i = \kappa_i$.
So for some $m(*) < \omega$ and 
$\kappa(*) < \mu$ the set $\{i < \kappa:|{\frak d}_i| =
m(*),\kappa(i) \le \kappa(*)\} \in D^+_1$ so \wilog
\mr
\item "{$(*)_4$}"  $i < \kappa \Rightarrow |{\frak d}_i| = m(*) \and
\kappa_i = \kappa(*)$.
\ermn
By (d) of $(*)_1$ there is $\bar f$ such that
\mr
\item "{$(*)_5$}"  $\bar f = \langle f_\alpha:\alpha < \text{ cf}(\lambda)
\rangle$ is $<_D$-increasing and cofinal in $\dsize \prod_{i < \kappa} \,
\lambda_i$ and if $\delta < \text{ cf}(\lambda),\text{cf}(\delta) < \mu$;
and $\bar f \restriction \delta$ has a $<_D$-eub then $f_\delta$ is
such a $<_D$-eub and let $f'_\alpha \in \dsize \prod_{i < \kappa}
\lambda_i$ be $f'_\alpha(i) = \text{ Min}(\text{Rang}(h_i)
\backslash f_\alpha(i))$ and $f''_\alpha(i) = h^{-1}_i(f'_\alpha(i))$.
\ermn
For each $i$ let $\bar{\Cal P}^i = \langle {\Cal P}^i_\alpha:\alpha <
\text{ cf}(\lambda_i)\rangle$ be such that 
$(*)^{{\frak d}_i,\kappa}_{\lambda_i,\bar{\Cal P}^i,h_i}$ holds.  We 
now choose $M_\alpha$ for $\alpha < \text{ cf}(\lambda)$ such that
\mr
\item "{$\circledast_2$}" $(a) \quad M_\alpha \prec ({\Cal H}(\chi),\in)$
\sn 
\item "{${{}}$}"  $(b) \quad \|M_\alpha\| <
\lambda,M_\alpha$ is increasing continuous, $\beta < \alpha 
\Rightarrow h(\beta) \subseteq M_{\alpha +1}$ and \nl

\hskip22pt $\beta < \alpha \Rightarrow \langle
M_\gamma:\gamma \le \beta \rangle \in M_{\alpha +1}$
\sn
\item "{${{}}$}"  $(c) \quad$ the following objects belong to
$M_\alpha:\langle \bar{\Cal P}^i:i < \kappa \rangle$, \nl

\hskip22pt $\langle \lambda_i,h_i:i < \kappa \rangle,\bar f,D$ and
$\mu$
\sn
\item "{${{}}$}"  $(d) \quad$ if $A \in D^+_1$, so $T_{D+A}(\langle
|{\Cal P}^i_{f_\alpha(i)}|:i < \kappa \rangle) < \lambda$ then
$T_D(f_\alpha) +1 \subseteq M_{\alpha +1}$ \nl

\hskip22pt (remember cf$(\lambda) > \mu > 2^\kappa$).
\ermn
Let ${\frak d}^* = \{\theta:\theta = \kappa$ or 
$\{i < \kappa:\theta \notin {\frak d}_i\} =
\emptyset$ mod $D_1\}$; it should be clear that $|{\frak d}^*| \le m(*)+1$.
\nl
Let ${\Cal P}_\alpha = M_{\alpha +1} \cap [h(\alpha)]^{< \mu}$ and
$\bar{\Cal P} = \langle {\Cal P}_\alpha:\alpha < \text{ cf}(\lambda)\rangle$.

It is enough now to prove that 
$(*)^{{\frak d}^*,\kappa(*)}_{\lambda,\bar{\Cal P},h}$ holds.   
\nl
Let $A \subseteq \lambda,|A| < \mu$ be as in the assumption and 
we should find $\bold c:[A] \rightarrow \kappa(*)$ as required.  
For $i < \kappa$ let
$A_i = \{f_\alpha(i):\alpha \in A\}$, so $A_i \in [\lambda_i]^{< \mu}$ 
hence there is $\bold c_i:[A_i]^2 \rightarrow \kappa(*)$ as required.  
Recalling that $\kappa(*)^\kappa = \kappa(*)$, we can 
define $\bold c:[A]^2 \rightarrow \kappa(*)$ such that
\mr
\item "{$\circledast_3$}" if $\alpha_1 < \beta_1,\alpha_2 < \beta_2$ are
from $A$ and $\bold c\{\alpha_1,\beta_1\} = \bold
c\{\alpha_2,\beta_2\}$ \ub{then}
{\roster
\itemitem{ $(i)$ }  if  $i < \kappa$ then $f_{\alpha_1}(i) <
f_{\beta_1}(i) \equiv f_{\alpha_2}(i) < f_{\beta_2}(i)$
\sn
\itemitem{ $(ii)$ }  if  $i < \kappa$ then $f_{\alpha_1}(i) >
f_{\beta_1}(i) \equiv f_{\alpha_2}(i) > f_{\beta_2}(i)$
\sn
\itemitem{ $(iii)$ }  if  $i < \kappa$ and $f_{\alpha_1}(i) <
f_{\beta_1}(i)$ then $\bold c_i\{f_{\alpha_1}(i),f_{\beta_1}(i)\} =
\bold c_i\{f_{\alpha_2}(i),f_{\beta_2}(i)\}$.
\endroster}
\ermn
Let $\theta \in \text{ Reg } \cap \mu \backslash {\frak d}^*$ and let
$\delta < \text{ cf}(\lambda)$ and $B \subseteq A \cap h(\delta)$ be 
such that $\bold c \restriction [B]^2$ is
constantly $i$ and $\theta = \text{ cf}(\delta)$ and $\delta =
\sup(B)$.  We can replace $D$ by $D + \{i < \kappa:\theta \notin
{\frak d}_i\}$.  So for some set $a \subseteq \theta$ we have
\mr
\item "{$\circledast_4$}"  if $\alpha < \beta$ are from $B$ then $a =
\{i < \kappa:f_\alpha(i) < f_\beta(i)\}$.
\ermn
Clearly $a \in D$.  Note that by $\circledast_3$ for 
each $i \in a$ the sequence $\langle
f_\alpha(i):\alpha \in B \rangle$ is increasing and let $B_i =
\{f_\alpha(i):\alpha \in B\}$ so $\delta_i =: \sup(B_i)$ has
cofinality $\theta$ and $\bold c_i \restriction [B_i]^2$ is
constant.  Hence by the choice of $\bar{\Cal P}^i$ and $\circledast(iii)$
clearly $B_i \in {\Cal P}^i_{\delta_i}$.  Also by the choice of 
$a$ and $\circledast$ above (and \cite[II,\S1]{Sh:g}) $\bar f \restriction
\delta$ has a $\le_D$-eub $f',f'(i) =: \cup\{f_\alpha(i):i \in B\}$,
hence \footnote{Note that here we use $\theta \ne \kappa$, in fact the
only point that we use it, if we could avoid it, then ${\frak d}$ could be
chosen as $\{\aleph_0\}$}
$a' = \{i \in a:f_\delta(i)=\delta_i\}$
belongs to $D$.  Now $|{\Cal P}^i_{f_\delta(i)}| < \lambda_i$, hence
$T_D(\langle {\Cal P}^i_{f_\delta(i)}:i < \kappa \rangle) <
\lambda$, so there is ${\Cal F} \subseteq \dsize \prod_{i < \kappa}
{\Cal P}^i_{f_\delta(i)},|{\Cal F}| < \lambda$ such that for every $g
\in \dsize \prod_{i < \kappa} {\Cal P}^i_{f_\delta(i)}$ there is $g'
\in {\Cal F}$ such that $\{i:g(i)=g'(i)\} \in D^+$.  By the choice 
of $M_{\delta +1}$, i.e., clause (d) of $\circledast_2$ 
\wilog \, ${\Cal F} \in M_{\delta +1}$, hence ${\Cal F}
\subseteq M_{\delta +1}$.  We can define $g \in 
\dsize \prod_{i < \kappa} {\Cal P}^i_{f_\delta(i)}$ 
by letting $i \in a' \Rightarrow g(i) =
B_i$.  So there is $g' \in {\Cal F} \subseteq M_{\delta +1}$ such that $b =
\{i:g(i) = g'(i)\} \in D^+$ hence $b \cap a' \in D^+$.  That is $b' =:
\{i \in a':g'(i) = B_i\} \in D^+$.  Clearly $b' \in M_{\delta
+1}$ (as $\mu \in M_{\delta +1}$ hence ${\Cal H}(\mu) \subseteq
M_{\delta +1})$ and $g' \in M_{\delta +1}$, hence $g' \restriction b' \in
M_{\delta +1}$, hence also the set $B^*$ belongs to $M_{\delta +1}$ where

$$
B^* =: \{\gamma < \lambda:\{i \in b':f_\gamma(i) \in g'(i) = g(i) = B_i\} \in
D^+\}.
$$
\mn
Now $|B^*| \le \dsize \prod_{i < \kappa} B_i < \mu$ and
$\alpha \in B \Rightarrow \alpha \in B^*$.  But as $B^* \in M_{\delta
+1}$ every subset of $B^*$ belongs to $M_{\delta +1}$ hence $B \in
M_{\delta +1}$ so $B \in {\Cal P}_\delta$, as required.
\hfill$\square_{\scite{md.1}}$\margincite{md.1}
\enddemo
\bn
\margintag{md.1Q}\ub{\stag{md.1Q} Discussion}:  1) Note that in a sense what was done in
\cite{Sh:108}, i.e., $I[x]$ large for $\lambda = \mu^+$ is done in
\scite{md.1} for any $\lambda$ with cf$(\lambda) > \mu$.
\nl
2) We may consider replacing 
${\frak d}_\lambda$ by $\{\aleph_0\}$ in \scite{md.1}.  The base of the
induction is clear (pcf$_{\aleph_1}$-inaccessibility).  So eventually we
have $f_\delta$ for it as above $\langle f_\alpha:\alpha \in B
\rangle$, the hard case is cf$(\text{otp}(B)) = \kappa$, we have the induced $h_* \in {}^\kappa
\kappa$ such that $\alpha < \kappa \Rightarrow \{i:d < h_*(i)\} \in
D$, but $(\forall^D i)[\text{cf}(h_*(i)) = \aleph_0]$ (otherwise using
niceness of the filter, etc., we are done).

Note that is problem appear even in the problem ``assume $\mu$ is
strong limit of cofinality $\aleph_1$ (or $\kappa \in [\aleph_1,\mu))$
and $2^\mu = \mu^+$, does it follows that
$\diamondsuit_{S^{\mu^+}_{\text{cf}(\mu)}}$ holds?"  See
\cite{Sh:186}, Cummings Dzamonja Shelah \cite{CDSh:571}, Dzamonja
Shelah \cite{DjSh:692}.

But if $\kappa = \text{ cf}(\kappa) > \aleph_0$ and in \scite{h.9} we
use $D = D_\kappa + S^\kappa_{\aleph_1}$.  For each $\alpha < \kappa$
we should consider $\iota(t)$; if $D$-positively we have $\iota(t) \le
h_*(t)$ we are done.  But if $\iota(t) > h_*(t),D$-positively 
then on some $A \in D^+,h_* \restriction A$ is constant.
\bigskip

\demo{\stag{md.2} Conclusion}  Assume $\mu < \lambda,\mu$ is strong 
limit $> \aleph_0,\lambda$ is regular (or just cf$(\lambda) > \mu$).  
\ub{Then} for some $\kappa <
\mu$ and finite ${\frak d} \subseteq \text{ Reg } \cap \mu$ to which
$\aleph_0$ belongs (in fact $({\frak d}^0_\mu(\lambda)
\cup\{\aleph_0\})$ is O.K.), there is $\bar{\Cal F}$ such that
\mr
\item "{$\circledast^{\mu,{\frak d},\kappa}_{\lambda,\bar{\Cal F}}$}"
$(a) \quad \bar{\Cal F} = \langle {\Cal F}_\alpha:\alpha <
\lambda \rangle,|{\Cal F}_\alpha| < \lambda$ for $\alpha < \lambda$
\sn
\item "{${{}}$}" $(b) \quad {\Cal F}_\alpha \subseteq \{f:f$ a partial
function from $\alpha$ to $\alpha,|\text{Dom}(f)| < \mu\},
{\Cal F}_\alpha$ closed \nl

\hskip25pt under restriction 
\sn
\item "{${{}}$}" $(c) \quad$ for every $A \subseteq \lambda,|A| < \mu$
and $f:A \rightarrow \lambda$ for some $\bold c:[A]^2 \rightarrow
\kappa$ we have
{\roster
\itemitem{ $\boxdot_1$ }   if $B \subseteq A,\delta = \sup(B) \in E,\bold c
\restriction [B]^2$ is constant, $[\alpha \in B \Rightarrow f(\alpha) <
\delta]$ and cf$(\delta) \notin {\frak d}$ then $f \restriction B \in
{\Cal F}_\delta$ and $\alpha \in B \Rightarrow f \restriction (B \cap
\alpha) \in {\Cal F}_\alpha$.
\endroster}
\endroster
\enddemo
\bigskip

\demo{Proof}   We use the result of \scite{md.1}.

For clause (c) we use pairing function pr on $\lambda$ such that
pr$(\alpha,\beta) <$ 
Max$\{\omega,\alpha + |\alpha|,\beta + |\beta|\}$ to
replace the function $f$ in clause (c) by the set
$\{\text{pr}(\alpha,f(\alpha)):\alpha \in A\}$ and first we restrict
ourselves to $\delta$ in some club $E$ of $\lambda$ (the range of $h$
in \scite{md.1}'s notation) such that $\delta \in E \Rightarrow |\delta|$ 
divides $\delta$ (hence $\delta$ is closed under pr); 
so if $B \subseteq \lambda,\sup(B) \in E$ we are done.  The other
cases are easier as \wilog \, if $\alpha < \delta \in E$, then $\alpha +
\text{ Min}\{\chi:\mu \ge |\alpha|$ and if ${\frak a} \subseteq \text{
Reg } \cap \chi^+,|{\frak a}| < \mu$,
pcf$_{\kappa(*)^+\text{-comp}}){\frak a}) \subseteq \mu^+\} < \delta$,
and easy to finish.   \hfill$\square_{\scite{md.2}}$\margincite{md.2}
\enddemo
\bigskip

\demo{\stag{md.3} Conclusion}   Assume that $\mu$ is strong limit
singular, $\lambda = \lambda^{< \mu}$ (equivalently $\lambda = \lambda^\mu$)
 and $c \ell:[\lambda]^{< \mu}
\rightarrow [\lambda]^{< \mu}$ satisfies for notational simplicity $c
\ell(B) = \cup\{c \ell(B \cap (\alpha +1)):\alpha \in B\}$ and $B_1
\trianglelefteq B_2 \Rightarrow B_1 \subseteq c \ell(B_1) \subseteq c
\ell(B_2)$.

\ub{Then} in \scite{md.2} we can add to (a),(b) and (c) also
\mr
\item "{$(d)$}"  $\bold g$ is a function from $\{f \restriction
u:f \in {}^\lambda \lambda$ and $u \in [\lambda]^{< \mu}\}$ to $\lambda$
\sn
\item "{$(e)$}"  for every $f:\lambda \rightarrow \lambda$ for some
$g_f:[\lambda]^{< \mu} \rightarrow \lambda$ (in fact $g_f(u) = \bold
g(f \restriction c \ell(u))$ we have
{\roster
\itemitem{ $\boxtimes$ }  for every $A \subseteq \lambda$ of
cardinality $< \mu$ such that $\alpha \in A \Rightarrow g_f(A \cap
\alpha) < \alpha$, for some function $\bold c:[A]^2 \rightarrow
\kappa$ we have
\sn
\itemitem{ ${{}}$ }  $\otimes \quad$ if $B \subseteq A,\bold c
\restriction [B]$ is constant and $B$ has no last element, \nl

\hskip45pt $\delta =
\sup(B)$ has cofinality $\notin {\frak d}$ \ub{then} $f \restriction 
 c \ell(B)$ \nl

\hskip45pt belong to ${\Cal F}_\delta$ and 
$\alpha \in B \Rightarrow f \restriction c \ell(B \cap \alpha) \in
\cup\{{\Cal F}_\beta:\beta < \delta\}$
\endroster}
\sn
\item "{$(f)$}"  if $\lambda$ is regular then there is a sequence
$\bar C = \langle C_\delta:\delta \in S \rangle$ is
{\roster
\itemitem{ $(\alpha)$ }  $S \subseteq S^* = \{\delta < \lambda$:
cf$(\delta) \in [\aleph_1,\mu)\}$,
\sn
\itemitem{ $(\beta)$ }  $C_\delta$ a club of $\delta$ of
order type cf$(\delta)$ \nl
and in clause (e) we can add:
\sn
\itemitem{ $(\gamma)$ }    $f \restriction c \ell(C_\delta) 
\in {\Cal F}_\delta$ and 
\sn
\itemitem{ $(\delta)$ }  $\alpha \in C_\delta
\Rightarrow f \restriction c \ell(C_\delta \cap \alpha) \in
\dbcu_{\beta < \delta} {\Cal F}_\beta$ and
\sn
\itemitem{ $(\varepsilon)$ }   $\alpha \in \text{ nacc}
(C_{\delta_1}) \cap \text{ nacc}(C_{\delta_2}) \Rightarrow
C_{\delta_1} \cap \alpha = C_{\delta_2} \cap \alpha$
\sn
\itemitem{ $(\zeta)$ }  if $\alpha < \text{ cf}(\lambda)$ is limit,
cf$(\alpha) \notin {\frak d}$ then $\{C_\delta \cap \alpha:\alpha \in
\text{ acc}(C_\delta)\}$ has coordinates $< \lambda$
\sn
\itemitem{ $(\eta)$ }  if $\beta \le \lambda,|B| < \mu$ then for some
$\bold c:[B]^2 \rightarrow \kappa$ if $B' \subseteq B$ has no last
member and $\bold c \restriction B'$ is constant and cf$(\sup(B'))
\notin {\frak d}$ then sup$(B') \in S$.
\endroster}
\endroster
\enddemo
\bigskip

\demo{Proof}  We repeat the proof of \scite{md.1}.

Choose $h:\text{cf}(\lambda) \rightarrow \lambda$ and
$\langle M_\alpha:\alpha < \text{ cf}(\lambda)\rangle$ as in 
the proof of \scite{md.1} but add the requirement that 
$c \ell \in M_0$ and still
use ${\Cal F}_\alpha = M_{\alpha +1} \cap \{f:f$ a partial function
from $\alpha$ to $\alpha$ with domain of cardinality $< \mu\}$.

Choose $\bold g$ such that
\mr
\item "{$\boxtimes$}"  $(a) \quad \bold g$ is a function from $\{f
\restriction u:f \in {}^\lambda \lambda$ and 
$u \in [\lambda]^{< \mu}\}$ onto $\lambda$
\sn
\item "{${{}}$}"  $(b) \quad f_1 \subseteq f_2 \in \text{ Dom}(\bold g)
\Rightarrow \bold g(f_1) \le \bold g(f_2)$ and
\sn
\item "{${{}}$}"  $(c) \quad$ for each $\alpha < \lambda$ 
for some $f \in \text{ Dom}(\bold g)$ we have $\bold g(f) = \alpha
\and$ \nl

\hskip20pt $(\forall f')[\bold g(f') = \alpha \Rightarrow f' \subseteq f]$
\sn
\item "{${{}}$}"  $(d) \quad$ if $f:B_2 \rightarrow \lambda$ and $B_1
\triangleleft B_2$ then $\bold g(B_1) < g(B_2)$
\sn
\item "{${{}}$}"  $(e) \quad \bold g(f) = \alpha \Rightarrow \text{
Dom}(f) \subseteq \alpha$.
\ermn
Without loss of generality $\bold g \in M_0$ so clause (d) (of the
conclusion of \scite{md.3}) holds
trivially and let us prove clause $(e)$.  So, as $\bold g$ has already
been chosen, we
are given $A \subseteq \text{ cf}(\lambda)$ of cardinality $< \mu$ and $f:A
\rightarrow \lambda$ such that $\alpha \in A \Rightarrow 
\bold g(f \restriction (A \cap \alpha)) < \alpha$.

Now $\alpha \mapsto \bold g(f \restriction (A \cap \alpha))$ is an
increasing function from $A$ to $\lambda$ and let $A' = \{\bold g(f
\restriction (A \cap \alpha)):\alpha \in A\}$ and let $\bold c':[A']^2
\rightarrow \kappa$ be as proved to exist in \scite{md.1} and by
$\bold c:[A]^2 \rightarrow \kappa$ be defined by $\bold
c\{\alpha,\beta\} = \bold c'\{\bold g(f \restriction (A \cap
\alpha)),\bold g(f \restriction (A \cap \beta))\}$.  

It is easy to check that $\bold c$ is as required.  We turn to
proving clause (f) of the claim.  Now there is a function
$F:{}^\omega \lambda \rightarrow \lambda$
such that for any $\bar \alpha \in {}^\omega \lambda$ for every large
enough $n < \omega$ there are $m_0 < m_1 < m_2 < \ldots < \omega$ which are $>
n$ and $\alpha_n = F(\alpha_{m_0},\alpha_{m_1},\ldots)$; by \cite{EH}.
For any $u \in [\lambda]^{< \mu}$ we define $c \ell_*(u)$ as follows:
let $u^{+ \bold g} = u \cup \{\bold g(1_v):v \subseteq u \cap \alpha$
and $\alpha \in u\}$ and let $c \ell_*(u)$ be the minimal set
$v$ such that $u^{+ \bold g} \subseteq v$ and $[\delta= \sup(v \cap
\delta) < \sup(u^{+ \bold g}) \and \text{ cf}(\delta) \le |u|
\Rightarrow \delta \in v]$ and $[\bold g(1_w) \in v \and |w| \le |u|
\Rightarrow w \subseteq v]$  and
$\bar\alpha \in {}^\omega v \Rightarrow F(\alpha) \in v$; so $|c
\ell_*(u)| \le (|u|^+ +2)^{\aleph_0}$.

In the proof above we can replace $c \ell$ by $c \ell_* \circ c \ell$.
Now if $\delta < \lambda,\aleph_0 < \text{ cf}(\delta) < \mu$
for some club $C^*_\delta$ of $\delta$ of order type cf$(\delta)$ we
have: if $C \subseteq C^*_\delta$ is a club of $\delta$ then $c
\ell_* \circ c \ell(C) = c \ell_* \circ c \ell(C^*_\delta)$ 
(exists by the choice of $F$).  Alternatively, let 
$C'_\delta = \cap \{c \ell_*(C):C$ a club of
$\delta\}$; however, $C'_\delta$ seemingly has order type just
$<(\text{cf}(\delta)^{\aleph_0})^+$.  Now if $C^*_\delta$ satisfies
$(\forall \alpha \in C^*_\delta)(\bold g(1_{C_\delta \cap \alpha}) <
\delta)$ then we can find $C^{**}_\delta,C_\delta$ such that:
\mr
\item "{$\circledast_1$}"  $C^{**}_\delta \subseteq c \ell_* \circ c
\ell(C^*_\delta)$ is a
club of $\delta$ of order type cf$(\delta)$ such that $\alpha \in
\text{ nacc}(C^{**}_\delta) \Rightarrow \sup((C^{**}_\delta
\cup \{0\}) \cap \alpha) < \bold g(1_{((C^{**}_\delta \cup\{0\}) \cap
\alpha)}) < \alpha$
\sn
\item "{$\circledast_2$}"  $C_\delta$ is $\{\bold g
(1_{((C^{**}_\delta \cup\{0\}) \cap \alpha)}):\alpha \in \text{
nacc}(C^{**}_\delta)\} \cup \text{ acc}(C^{**}_\delta)$.
\ermn
Clearly
\mr
\item "{$\circledast_3$}"  $C_\delta \subseteq c \ell_*(B)$ whenever
$B \subseteq \delta = \sup(B)$
\sn
\item "{$\circledast_4$}"  if $\alpha \in \text{ nacc}(C_{\delta_1})
\cap \text{ nacc}(C_{\delta_2})$ then $C_{\delta_1} \cap \alpha =
C_{\delta_2} \cap \alpha$.
\ermn
We are done as we have used $c \ell_* \circ c \ell$ and
\mr
\item "{$(*)$}"  if $\delta < \lambda,\aleph_0 < \text{ cf}(\delta) <
\mu$ and $B$ is an unbounded subset of $\delta$ then $C_\delta
\subseteq c \ell_*(B)$.
\nl
${{}}$  \hfill$\square_{\scite{md.3}}$\margincite{md.3}
\endroster
\enddemo
\bigskip

\remark{\stag{md.3A} Remark}  1) In \scite{md.1}, \scite{md.2},
\scite{md.3} if $\lambda$ is regular; then

$$
\align
A_{\bar M} = \{\delta:&\delta < \lambda,\text{cf}(\delta) < \delta
\text{ and there is} \\
  &u \subseteq \delta = \sup(u),\text{otp}(u) < \delta,
(\forall \alpha < \delta)(u \cap \alpha \in M_\alpha)\}
\endalign
$$
\mn
belongs to $I[\lambda]$ and the $\delta$ mentioned in $(*)^{{\frak
d},\kappa}_{\lambda,\bar{\Cal P}}$ of \scite{md.1},(c) of \scite{md.2}
necessarily belongs to $A_{\bar M}$.  So $A_{\bar M}$, for ordinals of
cofinality $\in \text{ Reg} \cap \mu \backslash {\frak d}$ contains
``almost all of them" in the appropriate sense.
\nl
2) We can use them to upgrade if $\{\delta < \omega_2:
S^{\beth_{\delta +1}}_\kappa \in I(\beth^+_\delta)\}$ 
then $S^{\beth_{\omega +1}}_\kappa \in 
I[\beth^+_{\omega_1+1}]$ when $\kappa = \text{ cf}(\kappa) >
\aleph_1$, see \cite{Sh:589}.
\endremark
\bigskip

\demo{\stag{md.5} Main Conclusion}  1) If $\mu$ is strong limit and
$\lambda = \lambda^{< \mu}$ \ub{then} for all but finitely many
regular $\kappa < \mu$ (actually $\kappa \notin {\frak d}^0_\mu
(\lambda) \cup \{\aleph_0\}$ is enough) we have 
Ps$_1(\lambda,\kappa)$, see Definition \scite{md.6} below.
\nl
2) We also get Ps$_1(\text{cf}(\lambda),\lambda,\kappa)$ when $\kappa
 > \aleph_0$.
\enddemo
\bigskip

\demo{Proof}  1) By \scite{md.2}, \scite{md.3}.
\enddemo
\bigskip

\definition{\stag{md.6} Definition}  1) Ps$_1(\lambda,\kappa)$ means
that Ps$_2(\lambda,S)$ for some stationary $S \subseteq S^\lambda_\kappa$.
\nl
2) Ps$_2(\lambda,S)$ means that for some $\bar C = \langle
C_\delta:\delta \in S \rangle$ we have Ps$_3(\lambda,\bar C)$.
\nl
3) Ps$_3(\lambda,\bar C)$ means that for some $\bar{\Cal F}$ we have
$\text{Ps}_4(\lambda,\bar C,\bar{\Cal F})$.
\nl
4) Ps$_4(\lambda,\bar C,\bar{\Cal F})$ means that for some $S$:
\mr
\item "{$(a)$}"  $S$ is a stationary subset of $\lambda$
\sn
\item "{$(b)$}"  $\bar C$ has the form $\langle C_\delta:\delta \in S
\rangle$
\sn
\item "{$(c)$}"  $\bar{\Cal F}$ has the form 
$\bar{\Cal F} = \langle {\Cal F}_\alpha:\alpha \in S \rangle$
\sn
\item "{$(d)$}"  $C_\delta$ is a club of $\delta$ of order type
cf$(\delta)$ and $\alpha \in \text{ nacc}(C_{\delta_1}) \cap
\text{ nacc}(C_{\delta_2} \Rightarrow \alpha \cap C_{\delta_1} =
\alpha \cap C_{\delta_2}$
\sn
\item "{$(e)$}"  ${\Cal F}_\delta$ is a set of functions from
$C_\delta$ to $\delta$ of cardinality $< \lambda$
\sn
\item "{$(f)$}"  if $f:\lambda \rightarrow \lambda$ then for
stationarily many $\delta \in S$ we have $f \restriction C_\delta \in
{\Cal F}_\delta$.
\ermn
5) Ps$_4(\lambda,\mu,h,\bar C,\bar{\Cal F})$ is defined similarly (so
$\lambda$ is regular) except that
\mr
\item "{$(e)'_1$}"  $h$ is an increasing continuous function from
$\lambda$ to $\mu$ with limit $\mu$
\sn
\item "{$(e)_2$}"  ${\Cal F}_\delta$ is a set of functions from
$\delta$ to $h(\delta)$ of cardinality $< \mu$
\sn
\item "{$(f)$}"  if $f:\lambda \rightarrow \mu$ then for stationarily
many $\delta \in S$ we have $f \restriction C_\delta \in {\Cal
F}_\delta$.
\ermn
6) If in (5) we omit $h$ we mean some $h$.
\nl
7) Ps$_1(\lambda,\mu,\kappa)$, Ps$_2(\lambda,\mu,S)$,
Ps$_3(\lambda,\mu,\bar C)$ are defined parallely. 
\enddefinition
\bigskip

\definition{\stag{md.6B} Definition}  Pr$_\ell$ are defined similarly
except not using $\bar C$ and ${\Cal F}_\delta$ is a set of functions
from some unbounded subset of $\delta$ into $\delta$ (or $h(\delta))$
that is:
\nl
1) Pr$_1(\lambda,\kappa)$ means that Pr$_2(\lambda,S)$ for some
stationary $S \subseteq S^\lambda_\kappa$.
\nl
2) Pr$_4(\lambda,\bar{\Cal F})$ means that for some $S$:
\mr
\item "{$(a)$}"  $S$ is a stationary subset of $\lambda$
\sn
\item "{$(b)$}"   $\bar{\Cal F}$ has the form 
$\bar{\Cal F} = \langle {\Cal F}_\alpha:\alpha \in S \rangle$
\sn
\item "{$(c)$}"  ${\Cal F}_\delta$ is a set of cardinality $< \lambda$ 
of functions from some
unbounded subset of $\delta$ to $\delta$
\sn
\item "{$(d)$}"  if $f:\lambda \rightarrow \lambda$ then for
stationarily many $\delta \in S$ we have $f \restriction C \in
{\Cal F}_\delta$ for some $A \subseteq \delta = \sup(C)$.
\ermn
3) Pr$_4(\lambda,\mu,h,\bar{\Cal F})$ is defined similarly to
except that
\mr
\item "{$(c)'_1$}"  $h$ is an increasing continuous function from 
$\lambda$ to $\mu$ with limit $\mu$
\sn
\item "{$(c)_2$}"  ${\Cal F}_\delta$ is a set of cardinality $<
\lambda$ of functions from some
unbounded subset of $\delta$ to $h(\delta)$
\sn
\item "{$(d)$}"  if $f:\lambda \rightarrow \mu$ then for stationarily
many $\delta \in S$ we have $f \restriction A \in {\Cal F}_\delta$ for
some $A \subseteq \delta = \sup(A)$.
\ermn
4) If in (5) we omit $h$ we mean some $h$.
\enddefinition
\bigskip

\demo{\stag{bb.6U} Observation}:  If Ps$_4(\lambda,\bar C,{\Cal F}),
\lambda_1 = \text{ cf}(\lambda) < \lambda,\bar C = \langle
C_\delta:\delta \in S \rangle,(\forall \delta \in S)[\text{cf}(\delta) >
\aleph_0],h:\lambda_1 \rightarrow \lambda$ is increasing continuous
with limit $\lambda,S' = \{\delta < \lambda_1:h(\delta) \in
S\},C'_\delta = \{\alpha < \delta:h(\alpha) \in C_\delta\},\bar C' =
\langle C'_\delta:\delta \in S'\},{\Cal F}'_\delta = \{h \circ f:f \in
{\Cal F}_\delta\}$ \ub{then} Ps$_4(\lambda_1,\lambda,h,\bar C',\bar{\Cal F}')$.
\enddemo
\bigskip

\demo{\stag{md.4} Conclusion}  1) If $\lambda = \text{ cf}(\lambda) >
\mu > \aleph_0,\mu$ strong limit singular \ub{then} for some $A \in
I[\lambda],\kappa < \mu$ and finite ${\frak d} \subseteq \text{ Reg }
\cap \mu$ (in fact ${\frak d} = {\frak d}'_{0,\mu}(\lambda)$ we have:
\mr
\item "{$(*)$}"   for every $\kappa(2) = \kappa(2)^{\kappa(1)} <
\mu,\kappa(1) > \kappa$ and increasing continuous sequence $\langle
\alpha_\varepsilon:\varepsilon < \kappa(2)^+ \rangle$ we have:
there is a club $C$ of $\kappa(2)^+$ such that $\{\alpha \in
C:\text{cf}(\alpha) \notin {\frak d}$ and cf$(\alpha) \le
\kappa(1)^+\} \subseteq A$.
\ermn
2) If above $\lambda = \lambda^{< \lambda}$ we can add: $\kappa \in
\text{ Reg } \cap \mu \backslash {\frak a} \Rightarrow (D
\ell)_{S^\lambda_\kappa}$ (and is even $(D \ell)_S$ for any $S
\subseteq S^\lambda_\kappa$ which is $\ne \emptyset$ modulo for a 
suitable filter similarly in (3)). 
\enddemo
\bn
On diamond from instances of GCH see \cite{Sh:460}.  
Whereas $\lambda = \mu^+$ successor of regular cardinals has
strong partial squares (\cite[\S4]{Sh:351}) for successor of singular
we have much less.  If $\lambda = \mu^+,\mu^\theta = \mu$ for
cofinalities $\le \theta$ we still have this.
\demo{\stag{md.4A} Conclusion}  Assume $\lambda = { \text{\rm
cf\/}}(\lambda) > \mu > \aleph_0,\mu$ strong limit and ${\frak d} =
{\frak d}'_{1,\mu}(\lambda) = \{\kappa:\kappa = { \text{\rm cf\/}}(\kappa)$
and $\kappa \in {\frak d}'_0(|\alpha|)$ for arbitrarily large $\alpha <
\lambda\}$ which is finite.
If $\lambda = \chi^+ = 2^\chi$ and $\kappa \in { \text{\rm Reg\/}} \,
\cap \mu \backslash {\frak d}$ then $\diamondsuit_{S^\lambda_\kappa}$. 
\enddemo
\newpage

\head {\S4 Middle Diamonds and Black Boxes} \endhead  \resetall \sectno=4
 \spuriousreset
\bigskip

We use \S3 to improve the main results of \cite{Sh:775} on the middle
diamond.  The point is that there we use \cite{Sh:460}, here we use
\S3 instead.  For our aim we quote some results and definitions.  See
\scite{md.12} and \scite{md.13}.  The conclusions are
\bn
In full
\proclaim{\stag{md.10} The Middle Diamond Claim}  Assume
\mr
\item "{$(a)$}"   $\lambda = {\text {\rm cf\/}}(2^\mu),D$ is a
$\mu^+$-complete filter on $\lambda$ extending the club filter
\sn
\item "{$(b)$}"   $\kappa = {\text {\rm cf\/}}(\kappa) < \lambda$ and
$S \subseteq S^\lambda_\kappa$
\sn
\item "{$(c)$}"   $\bar C = \langle C_\delta:\delta \in S
\rangle,C_\delta$ a club of $\delta$ of order type $\kappa$ and
$\lambda = 2^\mu \and \delta \in S \Rightarrow \lambda > 
|\{C_\delta \cap \alpha:\alpha \in { \text{\rm nacc\/}}(C_\delta)\}$
and $S \in D$
\sn
\item "{$(d)$}"   $2^{< \chi} \le 2^\mu$ and $\theta \le \mu$ 
\sn
\item "{$(e)$}"   {\rm Ps}$_1(\lambda,2^\mu,\bar C)$ (see Definition
\scite{md.6}) 
\sn
\item "{$(f)$}"  {\rm Sep}$(\mu,\theta)$ (see Definition \scite{bb.6}
below and \scite{bb.6P} on sufficient conditions).
\ermn
\ub{Then} $\bar C$ has the $(D,2^\mu,\theta)$-{\rm Md}-property
{\text{\rm (\/}}recall that this
means that the number of colours is $\theta$ not just 2, see
Definition \scite{md.13} below.
\endproclaim
\bigskip

\demo{Proof}  By the proof 
of \cite[1.10]{Sh:775} or see the more general claim on the squared
middle diamond (here or in subsequent work).
\enddemo
\bigskip

\definition{\stag{bb.6} Definition}   1) Sep$(\mu,\theta)$ means that
for some $\bar f$ and $\Upsilon$:
\mr
\item "{$(a)$}"  $\bar f = \langle f_\varepsilon:\varepsilon < \mu
\rangle$
\sn
\item "{$(b)$}"  $f_\varepsilon$ is a function from ${}^{\mu} \theta$
to $\theta$
\sn
\item "{$(c)$}"  for every $\partial \in {}^\chi \theta$ the set
$\{\nu \in {}^\mu \theta$: for every $\varepsilon < \chi$ we have
$f_\varepsilon(\nu) \ne \partial(\varepsilon)\}$ has cardinality $<
\Upsilon$
\sn
\item "{$(d)$}"  $\Upsilon = \text{ cf}(\Upsilon) \le 2^\mu$.
\ermn
2) Sep$_\sigma(\mu,\theta)$ means that for some $\bar f,R$ and
   $\Upsilon$ we have
\mr
\item "{$(a)$}"  $\bar f = \langle f^i_\varepsilon:\varepsilon < \mu
\text{ and } i < \sigma \rangle$
\sn
\item "{$(b)$}"  $f^i_\alpha$ is a function from ${}^M \theta$ to
 ${}^\mu \theta$
\sn
\item "{$(c)$}"  $R \subseteq {}^\mu \theta;|R| = 2^\mu$ (if $R =
{}^\mu \theta$ we may omit it)
\sn
\item "{$(d)$}"  $\bar{\Cal I} = \langle {\Cal I}_i:i < \sigma
\rangle,{\Cal I}_i \subseteq {\Cal P}(\mu)$ and if $A_j \in {\Cal
I}_j$ for $j<j^* < \sigma$ then $\mu \ne \cup\{A_j:j < j^*\}$
(e.g. ${\Cal I}_i = {\Cal I}$ a $\sigma$-complete ideal on $\mu$
\sn
\item "{$(e)$}"  if $\eta \in {}^\mu \theta$ and $i < \sigma$ then
$\Upsilon > |\text{Sol}_\eta|$ where
$$
\text{Sol}_\eta =\{\rho \in R:\text{ the set } \{\varepsilon <
\mu:(f^i_\alpha(\eta))(\varepsilon) \ne \eta(\varepsilon)\} \text{
belong to } I\}.
$$ 
\endroster
\enddefinition
\bn
We may wonder if clause (f) of the assumption is reasonable; the
following Claim gives some sufficient conditions for clause (f) to hold.
\bigskip

\proclaim{\stag{bb.6P} Claim}  Clause (f) of  \scite{md.10} holds, 
i.e., {\rm Sep}$(\mu,\theta)$ holds, \ub{if} at least one of the 
following holds:
\mr
\item "{$(a)$}"  $\mu = \mu^\theta$
\sn
\item "{$(b)$}"   $\bold U_\theta(\mu) = \mu + 2^\theta \le \mu$,
\sn
\item "{$(c)$}"   $\bold U_J(\mu) = \mu$ where for some $\sigma$ we
have $J = [\sigma]^{< \theta},\theta
\le \sigma,2^{< \sigma} < \mu$
\sn
\item "{$(d)$}"   $\mu$ is a strong limit of cofinality $\ne \theta =
{\text {\rm cf\/}}(\theta) < \mu$
\sn
\item "{$(e)$}"  $\mu \ge \beth_\omega(\theta)$.
\endroster
\endproclaim
\bigskip

\demo{Proof}  This is \cite[1.11]{Sh:775}.
\enddemo
\bigskip

\definition{\stag{md.12} Definition}  1) We say that $\lambda$ has the
$(\kappa,\theta)$-MD$^+$ property \ub{when} there are $\bar C^i = \langle
C_\delta:\delta \in S_i \rangle$ for $i < \lambda$ such that
\mr
\item "{$\boxdot^{\lambda,\kappa}_{\bar C}$}"  $(a) \quad S_i$
are pairwise disjoint stationary subsets of $\lambda$
\sn
\item "{${{}}$}"  $(b) \quad \delta \in S_i \Rightarrow \text{
cf}(\delta) =  \kappa$
\sn 
\item "{${{}}$}"  $(c) \quad C_\delta$ is a club of $\delta$ of order
type $\kappa$ and every $\alpha \in \text{ nacc}(C_\delta)$ is a \nl

\hskip30pt successor ordinal,
\sn
\item "{${{}}$}"  $(d) \quad$ if $\alpha \in \text{
nacc}(C_{\delta_1}) \cap \text{ nacc}(C_{\delta_2})$ then
$C_{\delta_1} \cap \alpha = C_{\delta_2} \cap \alpha$
\sn
\item "{${{}}$}"  $(e) \quad \bar C^i$ has the $\theta$-MD property which means
that there is $\bar f = \langle f_\delta:\delta \in S_i \rangle$ \nl

\hskip30pt such that
$f_\delta:{}^{\omega >}(C_\delta) \rightarrow \theta$ and for every $f \in
{}^{\omega >}\lambda \rightarrow \theta$ for \nl

\hskip30pt stationarily many $\delta
\in S_i$ we have $f_\delta = f \restriction C_\delta$.
\ermn
2) We write MD instead of MD$^+$ if we omit clause (d); we write
 MD$^\pm$ if we replace ``$C_\delta$ a club of $\delta$" by ``$C_\delta
 \subseteq \delta = \sup(C_\delta)$ and MD$^-$ if we do both changes.
\enddefinition
\bigskip

\definition{\stag{md.13} Definition}  1) We say that $\bar C$
exemplifies Md$^+(\lambda,\kappa,\theta,\Upsilon,D)$ when
\mr
\item "{$(a)$}"  $\lambda > \kappa$ are regular cardinals, $\Upsilon$
an ordinal (or a function with domain $\lambda$ or ${}^{\omega >}\lambda$ 
in this case a function $f$ from $X$ to $\Upsilon$ means that $f$ is a
function with domain $X$ and $f(x) \in \Upsilon(x)$, so
${}^C \Upsilon = \{f:f \text{ is a function with Dom}(f) = C$ and
$\alpha \in C \Rightarrow f(\alpha) \in \Upsilon(\alpha)\}$)
\sn
\item "{$(b)$}"  $C = \langle C_\delta:\delta \in S \rangle,S$ a
stationary subset of $\lambda$ such that $\delta \in S \Rightarrow
\text{ cf}(\delta) = \kappa$
\sn
\item "{$(c)^+$}"  $C_\delta$ is a club of $\delta$
disjoint to $S$ and $\alpha \in \text{ nacc}(C_{\delta_1}) \cap \text{
nacc}(C_{\delta_2}) \Rightarrow 
C_{\delta_2} \cap \alpha = C_{\delta_2} \cap \alpha$ so we
may define $C_\alpha = C_\delta \cap \alpha$ when $\alpha \in \text{ nacc}
(C_\delta)$
\sn
\item "{$(d)$}"  if $\bold F$ is a function from
$\dbcu_{\delta \in S} \{f:f$ is a function from ${}^{\omega >}(C_\delta)$ to
$\Upsilon\}$ to $\theta$ \ub{then}
for some $\bold c \in {}^S \theta$ for every $f \in {}^\lambda
\Upsilon$ the set $\{\delta \in S:\bold F(f \restriction C_\delta)) =
\bold c(\delta)\} \in D^+$. 
\ermn
2) We write Md instead Md$^+$ if we weaken (c)$^+$ to
\mr
\item "{$(c)$}"  $C_\delta$ is an unbounded subset of $\delta$.
\ermn
3) We say $\bar C$ has the $(D,\Upsilon,\theta)$-Md property \ub{when}
clauses (a),(b),(c),(d) above holds; we 
say $\lambda$ has this property if some
$\bar C = \langle C_\delta:\delta \in S \rangle$ has it, $S \subseteq
S^\lambda_\theta$ stationary.
\enddefinition
\bigskip

\remark{\stag{md.13K} Remark}  1) How strong is the demand that 
$S$ can be divided to $\lambda$ sets $S_i$ with the property?  
It is hard not to have it.
\nl
2) In \scite{md.11} to have more than one exception is a heavy demand
 on ${\Cal H}(\mu)$.
\nl
3) We can improve \scite{md.11} including the case cf$(\mu_*) =
 \aleph_0$, even $\mu_* = \beth_{\alpha + \omega}$.  Then probably
 in part (2) we have to distinguish $\lambda$ successor of regular
 (easy), success of singular (harder), rest (hardest).
\endremark
\bigskip

\demo{\stag{md.11} Conclusion}  1) If $\mu_*$ is strong limit $>
\aleph_0,\mu \ge \mu_*,\lambda = \text{ cf}(2^\mu)$ and $\Upsilon =
2^\mu$ \ub{then} for all but finitely many $\kappa \in 
\text{ Reg } \cap \mu_*$ (even every $\kappa \in \text{ Reg } \cap \mu_*
\backslash {\frak d}'_{1,\mu_*}(2^\mu))$, there is $\bar C 
= \langle C_\delta:\delta \in S \rangle$ exemplifies
Md$^+(\lambda,\kappa,\theta,\Upsilon)$ hence $(\kappa,\theta)$-MD$^+$.
\nl
2) Assume $\mu_*$ is strong limit singular of uncountable cofinality and
$\lambda = \text{ cf}(\lambda) > \mu_*$ is not strongly inaccessible.
\ub{Then} for all but finitely many $\kappa \in \text{ Reg } \cap
\mu_*$ for every $\theta < \mu_*,\lambda$ has $(\kappa,\theta)$-MD
hence $(\kappa,\theta)$-MD$^+$ (moreover only one of the exceptions
depend on $\lambda$).
\enddemo
\bigskip

\demo{Proof}  1) Let ${\frak d} = {\frak d}'_{\mu_*,0}(\lambda)$.  
So by \S3 we have $\kappa \in \text{ Reg } \cap
\mu_* \backslash {\frak d} \Rightarrow \text{ Ps}_1(\lambda,2^\mu,\bar
C)$ for some $\bar C$ satisfying clause (c) of \scite{md.10}, moreover
clauses (c) (d) of \scite{md.12}(1).  So we apply \scite{md.10}.
\nl
2) Let $\langle \mu_i:i <\text{ cf}(\mu_*) \rangle$ be increasing
continuous with limit $\mu_*$ each $\mu_i$ is strong limit
singular.  For each $i < \text{ cf}(\mu_*)$ let ${\frak d}_i =
{\frak d}'_{\mu_i,0}(\text{cf}(2^{\mu_i}))$, so it is
finite and let ${\frak d} = \{\kappa:\kappa = \text{ cf}(\kappa) <
\mu_*$ and $\kappa \in {\frak d}_i$ for every $i < \text{ cf}(\mu_*)$
large enough$\}$.
\enddemo
\bn
\ub{Case 1}:   $(\forall \alpha < \lambda)[|\alpha|^{< \mu_*} <
\lambda]$.  

So we can find $\mu < \lambda \le 2^\lambda$, let $\mu_1 =
((\mu)^{< \mu_*})^{< \mu_*}$ this cardinal is $< \lambda$ and $\mu_1 =
(\mu_1)^{\mu_*}$.  \nl
Now use \cite[\S2]{Sh:775}.
\bn
\ub{Case 2}:  $(\exists \alpha < \lambda)[|\alpha|^{< \mu_*} \ge
\lambda]$.  

As $\lambda$ is regular for some $\kappa < \lambda,\mu <
\lambda$ we have $\mu^\kappa \ge \lambda$.  Let $\mu = \text{
Min}\{\mu:\mu^\kappa \ge \lambda$ for some $\kappa < \mu_*\}$.
\nl
NOTE:  Here getting $\lambda$ pairwise disjoint $S_i$ should be
done.  Again we use \cite[\S2]{Sh:775}.
\bigskip

\remark{Remark}   $\aleph_0 \in {\frak d}$ as we need $F:{}^\omega \lambda
\rightarrow \lambda$ as in \S3!!  \hfill$\square_{\scite{md.11}}$\margincite{md.11}
\endremark
\bigskip

\definition{\stag{bb.1} Definition}  We say that $\bar C$ exemplify
BB$_0(\lambda,\kappa,\theta)$ when
\mr
\item "{$(a)$}"  $\lambda > \kappa$ are regular
\sn
\item "{$(b)$}"  $\bar C = \langle C_\delta:\delta \in S \rangle,S$ a
stationary subset of $\lambda$ such that $\delta \in S \Rightarrow
\text{ cf}(\delta) = \kappa$
\sn
\item "{$(c)$}"  $C_\delta$ is an unbounded subset of $\delta$
disjoint to $S$ such that $\alpha \in C_{\delta_1} \cap C_{\delta_2}
\Rightarrow C_{\delta_1} \cap \alpha = C_{\delta_2} \cap \alpha$
\sn
\item "{$(d)$}"  assume $\tau_0 \subseteq \tau_1 \subseteq \tau_2$ are
vocabularies of cardinality $\le \theta,\tau_1 \backslash \tau_0$ has
only predicates, $\tau_2 \backslash \tau_1$ has only function symbols
(allowed to be partial) ${\frak B}$ is a $\tau_0$-model with universe
$\lambda$ (if $\lambda = \lambda^\theta$ also individual constants) 
\ub{then} we can find $\langle {\Cal M}_\delta:\delta \in S
\rangle$ such that
{\roster
\itemitem{ $(\alpha)$ }  every $M \in {\Cal M}_\delta$ is a
$\tau_2$-model of cardinality $\theta$ expanding ${\frak B}
\restriction |M_\delta|$
\sn
\itemitem{ $(\beta)$ }  if $M \in {\Cal M}_\delta,F \in \tau_2 \backslash
\tau_1$ then $F^M$ has domain $\subseteq C_\delta$ (i.e., arity $F(C_\delta))$
\sn
\itemitem{ $(\gamma)$ }  every $M \in {\Cal M}_\delta$ has universe
which includes $C_\delta$ and is included in $\delta$ and the universe
of $M$ is the ${\frak B}$-closure of $C_\delta \cup \{F(\bar \alpha):F
\in \tau_2 \backslash \tau_1$ and $\bar \alpha \in
{}^{\text{arity}(F)}(C_\delta)\}$ 
\sn
\itemitem{ $(\delta)$ }  if $M',M'' \in {\Cal M}_\delta$ then
$(M',\gamma)_{\gamma \in C_\delta},(M'',\gamma)_{\gamma \in C_\delta}$
are isomorphic
\sn
\itemitem{ $(\varepsilon)$ }  if ${\frak B}^+$ is a
$\tau_2$-expansion of ${\frak B}$ then for stationarily many
$\delta \in S$ for some $M \in {\Cal M}_\delta$ we have:
\mr
\itemitem{ ${{}}$ }  $(i) \quad F \in \tau_2 \backslash \tau_2 
\Rightarrow F^{{\frak B}^+} \restriction C_\delta = F^M \restriction
C_\delta$ ($= F^M$)
\sn 
\itemitem{ ${{}}$ }  $(ii) \quad M \restriction \tau_1 \subseteq
{\frak B}^+ \restriction \tau_1$. 
\endroster}
\endroster
\enddefinition
\bigskip

\demo{\stag{bb.2} Observation}  1) In \scite{bb.1} if the order $<$ on
$\lambda$ is a relation of ${\frak B}$ (which is no loss) \ub{then}
the isomorphism is unique as it is necessarily the unique order
preserving function from $|M'|$ onto $|M''|$.
\nl
2) In \scite{bb.1}, if the function $F_i$ where $\alpha < \beta \in
C_\delta,\alpha \in C_\delta$, otp$(C_\delta \cap \alpha) = i
\Rightarrow F_i(\beta) = \alpha$, then for any $M \in \cup\{{\Cal
M}_\delta:\delta \in S\}$ and $\delta,M \cap C_\delta$ is
an initial segment of $C_\delta$.
\enddemo
\bigskip

\definition{\stag{bb.3} Definition}  We say that $\bar C$ exemplifies
BB$_1(\lambda,\kappa,\theta)$ when (a),(b),(d),(e) from \scite{bb.1}
holds + $(\varepsilon)$ below.  BB$_2(\lambda,\kappa,\theta)$ holds when we add
$(\zeta)$ to clause (d) where
\mr
\item "{$(\epsilon)$}"  the isomorphism type of $(M,\gamma)_{\gamma \in
C_\delta}$ for $M \in {\Cal M}_\delta$ depend on
$\tau_0,\tau_1,\tau_1$ but not on ${\frak B}$
\sn
\item "{$(\zeta)$}"  if $M',M'' \in {\Cal M}_\delta$ and $\pi$
is an isomorphism from $M'$ onto $M'$ and $\delta',\delta'' \in
S,C_{\delta'} \subseteq M',C_{\delta''} \subseteq M'$ and $\pi$ maps
$C_{\delta'}$ onto $C_{\delta''}$, \ub{then} for any $N' \in {\Cal
M}_{\delta'},N'' \in {\Cal M}_{\delta''}$ we have
$(N_\delta,\gamma)_{\gamma \in C_{\delta'}} \cong
(N_{\delta''},\gamma)_{\gamma \in C_{\delta''}}$.
\endroster
\enddefinition
\bigskip

\proclaim{\stag{bb.4} Claim}  If $\mu > \aleph_0$ is strong limit and
$\lambda = { \text{\rm cf\/}}(2^\mu)$ or $\lambda > 2^{2^\mu}$
not strongly inaccessible \ub{then} for all but finitely many $\kappa
\in { \text{\rm Reg\/}} \, \cap \theta \, (\kappa \in { \text{\rm
Reg\/}} \, \cap \mu \backslash {\frak d}'_1(2^\mu))$ for every $\theta < \mu$,
{\rm BB}$_1(\lambda,\kappa,\theta)$ holds.
\endproclaim
\bigskip

\demo{Proof}  Use also \scite{bb.11} below.  
\enddemo
\bigskip

\demo{\stag{bb.5} Observation}  1) If $\bar C$ exemplies
BB$_\ell(\lambda,\kappa,\theta)$ \ub{then} for some pairwise disjoint
$\langle S_\varepsilon:\varepsilon < \lambda \rangle$ we have each
$\bar C \restriction S_\varepsilon$ exemplifies
BB$_\ell(\lambda,\kappa,\theta)$. 
\nl
2) If $\lambda = \lambda^\theta$ we can allow in $\tau_1 \backslash \tau_0$
individual constant.
\nl
We delay their proof as we first use them.
\enddemo
\bn
Now we turn to proving \scite{bb.4}, \scite{bb.5}.
\proclaim{\stag{bb.11} Claim}  1) If $\bar C$ exemplifies
{\rm MD}$(\lambda,\kappa,2^\theta,\lambda)$ \ub{then} $\bar C$ exemplifies
{\rm BB}$_1(\lambda,\kappa,\theta)$ [Rethink: if we use $C * \chi,\chi
= \beth_\kappa$ enough to have many guess.]
\nl
2) $\bar C$ exemplifies {\rm BB}$_1(\lambda,\kappa,\theta)$ when there
are $\lambda_1,\bar C^1$
\mr
\item "{$(a)$}"  $\bar C$ exemplifies {\rm
MD}$(\lambda,\kappa,2^\theta,\lambda)$ (hence $\bar C^1 = \langle
C^1_\delta:\delta \in S_1 \rangle$ exemplifies {\rm
BB}$_1(\lambda,\kappa,\theta)$ but apparently we need more
\sn
\item "{$(b)$}"  $\bar h = \langle h_\delta:\delta \in S_1 \rangle$
where $h_\delta$ is an increasing function from $C_\delta$ onto some
$\gamma = \gamma(\delta) \in S_1$
\sn
\item "{$(c)$}"  for every club $C$ of $\lambda$ there is an
increasing continuous function $g$ from $\lambda_1$ into $C$ such that
$\alpha \in S_1 \Rightarrow g(\alpha) \in S \and \gamma_{g(\alpha)} =
\alpha$.
\ermn
3) If $\bar C$ exemplifies {\rm MD}$(\lambda,\kappa,2^\theta)$ then $\bar C$
 exemplifies {\rm BB}$_2(\lambda,\kappa,\theta)$.
\endproclaim
\bigskip

\demo{Proof}  1)  $\bar C$ has the
$(D,2^\mu,\theta)$-Md-property (which is like the desired conclusion
except that we write $F_\delta(\nu \restriction C_\delta)$ instead of
$F(\nu \restriction C_\delta,\bar C \restriction \bar C_\delta)$.  But
let $\beta = \alpha/\theta$ mean that $\theta \beta \le \alpha <
\theta \beta +1$.  But define $F'_\delta(\nu) = F_\delta(\langle
\nu(\alpha)/\theta:\alpha \in C_\delta \rangle,\langle
\nu(\alpha)-\theta(\nu(\alpha)/\theta):\alpha \in C_\delta \rangle)$.
So for $\langle F'_\delta:\delta \in S \rangle$ we have $\bar c$ as
required in the original requirement; the same $\bar c$ is as required
for our $\bar F$. \nl
2), 3) Left to the reader.   \hfill$\square_{\scite{bb.11}}$\margincite{bb.11} 
\enddemo
\bigskip

\demo{\stag{bb.12.Y} Conclusion}  If 
$\lambda = \text{ cf}(\lambda) > \beth_{\omega
+3}$ not strongly inaccessible, \ub{then} for every regular $\kappa <
\beth_\omega$ except possibly finitely many we have:
\mr
\item "{$\circledast$}"  for some topological space $X$ and $\bar C =
\langle C_\delta:\delta \in S \rangle$ we have
{\roster
\itemitem{ $(a)$ }  $X$ is Hausdorff with a clopen basis set
\sn
\itemitem{ $(b)$ }  every $Y \subseteq X$ of cardinality $< \kappa$ is
closed
\sn
\itemitem{ $(c)$ }  every point has a neighborhood of cardinality $\le \kappa$
\sn
\itemitem{ $(d)$ }  there is $f:X \rightarrow \kappa$ such that:
\nl
if $X = \dbcu_{\alpha < \beta} X_\alpha,\beta < \alpha$ \ub{then} some
non-isolated point $x$ has a neighborhood included in $x_{f(x)}$ (so $f(x) <
\beta$).
\endroster}
\endroster
\enddemo
\newpage

     \shlhetal 

\nocite{ignore-this-bibtex-warning} 
\newpage
    
REFERENCES.  
\bibliographystyle{lit-plain}
\bibliography{lista,listb,listx,listf,liste}

\enddocument